\documentclass[sn-mathphys]{sn-jnl}

\jyear{2021}

\usepackage{multirow}
\usepackage{dsfont}
\usepackage{comment}
\theoremstyle{thmstyleone}
\newtheorem{definition}{Definition}
\newtheorem{theorem}{Theorem}
\newtheorem{lemma}{Lemma}
\newtheorem{proposition}{Proposition}

\newtheorem{remark}{Remark}

\def\N{\mathbb{N}}
\def\R{\mathbb{R}}
\def\C{\mathbb{C}}
\def\Sp{\mathbb{S}}

\def\spvard{\mathcal{V}_d}
\def\origrass{\widetilde{G}^n_d}
\def\orientgrass{\widetilde{G}^n_1}
\def\vol{\text{vol}}

\begin{document}
\title{Weight metamorphosis of varifolds and the LDDMM-Fisher-Rao metric}
%\author{Hsi-Wei Hsieh and Nicolas Charon}
%\date{} % Activate to display a given date or no date (if empty),
         % otherwise the current date is printed
         
\author[1]{\fnm{Hsi-Wei} \sur{Hsieh}}\email{hsiwei.hsieh@austin.utexas.edu}

\author[2]{\fnm{Nicolas} \sur{Charon}}\email{charon@cis.jhu.edu}

\affil[1]{\orgdiv{Oden Institute for Computational Engineering and Sciences}, \orgname{The University of Texas at Austin}}

\affil[2]{\orgdiv{Center of Imaging Sciences}, \orgname{Johns Hopkins University}}

%\affil[3]{\orgdiv{Department}, \orgname{Organization}, \orgaddress{\street{Street}, \city{City}, \postcode{610101}, \state{State}, \country{Country}}}

\abstract{This paper introduces and studies a metamorphosis framework for geometric measures known as varifolds, which extends the diffeomorphic registration model for objects such as curves, surfaces and measures by complementing diffeomorphic deformations with a transformation process on the varifold weights. We consider two classes of cost functionals to penalize those combined transformations, in particular the LDDMM-Fisher-Rao energy which, as we show, leads to a well-defined Riemannian metric on the space of varifolds with existence of corresponding geodesics. We further introduce relaxed formulations of the respective optimal control problems, study their well-posedness and derive optimality conditions for the solutions. From these, we propose a numerical approach to compute optimal metamorphoses between discrete varifolds and illustrate the interest of this model in the situation of partially missing data.}


\maketitle

\section{Introduction}
Diffeomorphic shape analysis has come a long way since its origins in the 90s and some seminal works such as \cite{Grenander1993,christensen1996deformable}. Besides being at the origin of a constant development of various new mathematical models and numerical methods, it has further shown its wide potential for applications most notably to domains such as computational anatomy or computer vision. Although shape analysis is typically concerned with the usual issues of standard statistics, for instance the quantitative comparison of objects, the estimation of the mean of a population or of the directions of principal morphological variability, it remains such an active field primarily because of the very intricate mathematical structure of shape spaces that makes the generalization of those notions particularly delicate and still a largely open problem. Taking for example the case of shapes such as curves, surfaces or submanifolds in an Euclidean space, one generally needs to consider those objects as elements of a quotient of an infinite dimensional space by an infinite dimensional group, specifically as the equivalence classes of parametrization functions modulo all their reparametrizations (and in some cases other additional invariance groups such as rigid motions). Thus, even the definition of an adequate notion of metric is highly non trivial and much work has been conducted towards the construction and study of so called \textit{intrinsic} Riemannian metrics on such quotient shape spaces, see e.g. \cite{Michor2007,Bauer2011b}.

The construction of metrics and by extension the statistical analysis of shapes can be addressed through a different approach however, which was pioneered by the works of Grenander \cite{Grenander1993}. In Grenander's shape space framework, one views shapes as objects being acted on by a certain (potentially infinite dimensional) group $G$ of \textit{extrinsic} deformations. For the above situation of submanifolds embedded in $\R^n$, it could be for instance the group of diffeomorphisms $\text{Diff}(\R^n)$ or a subgroup such as rigid, affine or projective transformations that acts by transporting submanifolds. Then the distance between two shapes can be technically induced from a right-invariant distance on the deformation group itself by looking for a minimal deformation to transform one shape to the other. It is thus quite naturally that this line of work triggered the study of Riemannian metrics on diffeomorphism groups, among which the model coined \textit{Large Deformation Diffeomorphic Metric Mapping} (LDDMM) in \cite{Beg2005} proved particularly prolific in good part because it allows to operate with large deformations of the space and provides a principled approach to deal with a variety of geometric shapes including landmarks, images, curves and surfaces or even tensor fields. Yet, this metric formulation involves solving a registration problem i.e. finding an optimal deformation in the group between two given shapes, which a fortiori assumes that those two shapes belong to the same orbit for the action of $G$. As this is often not a realistic setting when dealing with real data or because of the difficulty of actually solving such boundary value problems, it is very common to relax the exact matching constraint and only enforce that the deformation maps the two shapes approximately as measured by some data attachment (or fidelity) term. For shapes such as landmarks or images, this measure of similarity can be simply taken as the sum of squared differences between the landmark positions or pixel values \cite{Joshi2000,Beg2005}. But the case of curves and surfaces is typically more elusive because of the aforementioned reparametrization invariance that needs to be embedded within the fidelity term. From a discrete perspective, this means a notion of discrepancy between two curves or two surfaces that does not assume predefined point correspondences and that is robust to differences of sampling and mesh structure. 

This precise issue motivated, in particular, the adaptation of ideas from the field of \textit{geometric measure theory} in order to obtain convenient representations for the design of adequate data attachment metrics in the case of submanifold data. The underlying principle is to map shapes into certain spaces of generalized measures and compare them through those measure representations. This was first proposed based on the framework of mathematical currents in \cite{Glaunes2,durrleman2009statistical} and later extended to the representations of \textit{varifolds} \cite{Charon2013} and \textit{oriented varifolds} \cite{kaltenmark2017general}. In each case, the construction of kernel metrics on the corresponding measure spaces lead to simple fidelity terms that can be effectively used in the above inexact diffeomorphic registration framework. But interestingly, this approach does not need to limit itself to submanifolds as varifolds in fact encompass a much wider category of geometric objects which can be loosely described as spatial distributions of local orientation planes. Thus, in their recent works \cite{hsieh2020diffeomorphic,hsieh2021metrics}, the authors suggested to formulate a more general diffeomorphic registration problem directly on varifolds themselves by considering a proper notion of group action of $\text{Diff}(\R^n)$ on the varifold space. 

Nevertheless, even formulated in the more general setting of varifolds, diffeomorphic models can remain insufficient in adequately dealing with some of the geometric variability encountered in data. Indeed, it is common for two given shapes to exhibit differences which cannot be entirely represented by a diffeomorphism. An obvious situation is the presence of topological changes between the two shapes. It can be also the result of imbalances, such as different fiber densities when comparing two fiber bundles which is common for instance with white matter fiber tracts obtained from diffusion MRI \cite{Gori2016}. The attempt to complement diffeomorphic deformations with transformations of a different nature was first formalized through the fundamental concept of \textit{metamorphoses} in \cite{Trouve1}. With subsequent works that include \cite{Holm2009,Richardson2013,Richardson2015,berkels2015time,Charon2018}, it appeared that metamorphoses can provide an effective framework to extend the Riemannian metric setting of LDDMM by incorporating a richer class of shape transformations. Yet metamorphoses have so far been primarily studied and implemented for images and landmarks. With the exception of the measure metamorphosis model of \cite{Richardson2013}, there has been very little to no work done on trying to adapt this framework to submanifolds and let alone to varifolds.

The main goal of this paper is precisely to address that issue by defining and studying a new Riemannian metric on varifolds based on a specific model of varifold metamorphosis that generalizes the purely diffeomorphic approach of \cite{hsieh2021metrics}. Our model essentially augments the diffeomorphic transport of a varifold with a dynamical change of its \textit{weight} (or mass) at each point. In order to associate a metric to the varifold metamorphosis, we take inspiration from a related class of models in the context of optimal transport and specifically the unbalanced optimal transport framework that was introduced independently in \cite{Liero2016} and \cite{chizat2018interpolating}. In its dynamic formulation, the Wasserstein-Fisher-Rao (or Hellinger-Kantorovich) metric between measures of $\R^n$ (i.e. $0$-dimensional varifolds) that is defined in those works combine usual optimal transport with a Fisher-Rao metric to penalize weight changes in the measures. By analogy, we introduce and study the \textit{LDDMM-Fisher-Rao metric} between varifolds where the optimal transport component is here replaced by a metric induced from the right-invariant metric on diffeomorphisms of the LDDMM model through its action on varifolds. Besides the mathematical analysis of this novel varifold metamorphosis model, we also tackle its numerical implementation for which we focus on discrete varifolds, namely finite sums of generalized Dirac masses, and again introduce a relaxation of the matching problem based on the aforementioned kernel fidelity terms. The practical interest of this numerical framework for data applications is multifold. Through some of the presented simulations, we will show that it can provide robustness to different types of density imbalances in structured and unstructured geometric data. But we shall also illustrate its potential to deal with \textit{partially observed data} with curves and surfaces, which has been a recurrent and challenging issue for many different shape analysis models \cite{bronstein2009partial,robinson2012functional,rodola2017partial,antonsanti2021partial,sukurdeep2021new,attaiki2021dpfm}. 

\vskip2ex

\textbf{Relationship to other works.} The approach we introduce in this paper relates but differs from several previous works in the following way. In the special case of $0$-dimensional varifolds, it leads to a metamorphosis metric between classical measures of $\R^n$ which shares the same diffeomorphic component as the measure metamorphosis model of \cite{Richardson2013} but combined with the Fisher-Rao weight transformation metric of unbalanced optimal transport \cite{Liero2016,chizat2018interpolating}. Unlike these two models, the LDDMM-Fisher-Rao metric is only well-defined between measures belonging to the same orbit under the combined action of diffeomorphisms and reweighting functions. Compared to \cite{Richardson2013} however, our use of a more constrained model and metric for the non-diffeomorphic part of the metamorphosis allows us to circumvent the theoretical issues that were uncovered by the authors of \cite{Richardson2013}; in particular we recover the existence of geodesics that remain in the space of positive measures. In contrast to unbalanced optimal transport on the other hand, the diffeomorphic component of our model guarantees smooth bijective geometric transformations which is often desirable in registration problems. Moreover, as we shall see, the generalization from measures to higher-dimensional varifolds and the corresponding change in the transport action induces further significant differences with these two models. 
%We point out in addition that the weighted varifold metamorphosis of this paper is also fundamentally different than the functional shapes framework of \cite{Charon2018} which introduced metamorphosis of functions on deformable curves and surfaces although such functions are instead interpreted as separate information that is independent of the shape itself, in sharp contrast with the weight functions considered here which transformation is coupled with the one of the object's geometry. 
Finally, related to the aforementioned challenge of partial data registration, we should mention the two recent works \cite{antonsanti2021partial} and \cite{sukurdeep2021new} which both examine alternative approaches that also rely, to some degree, on the varifold representation. The key difference consists in the fact that \cite{antonsanti2021partial} rather modifies the fidelity term used in the registration problem into a pseudo-distance that allows for partial overlap of the matched shape and the target while our approach actively models and estimate local varifold weight change itself. On the other hand, the method of \cite{sukurdeep2021new}, while also based on the estimation of weight changes, focuses on the space of curves equipped with intrinsic Sobolev metrics and does not fall in the setting of Riemannian metrics between varifolds that we follow in the present work. 

\vskip2ex 

\textbf{Structure of the paper.} The paper is organized as follows. In Section \ref{sec:weight_metam}, we start by reviewing basic definitions and properties of varifolds, their relationship to submanifolds and the action of the diffeomorphism group. We then introduce a generalized action with varifold weight changes and proceed in the definition of induced metrics from the resulting transformation model. Although we are primarily interested in the LDDMM-FR metric mentioned earlier (and introduced in Section \ref{ssec:LDDMM_FR}), we first discuss in Section \ref{ssec:L2_model} a simplified model involving a static $L^2$ penalty on the weight change. This approach was only briefly considered in our preliminary work \cite{hsieh2021diffeomorphic} that also mainly focused on the special case of constant weight change functions, and we shall expand its analysis in this paper. In the case of LDDMM-FR, we prove that it leads to a well-defined distance between varifolds in a given orbit and show the existence of geodesics. We further derive, in Section \ref{ssec:geodes_single_Diracs}, the exact expression of those geodesics between two single Dirac $0$- or $1$- varifolds. In Section \ref{sec:relaxed_problems}, we introduce relaxed versions of the matching problems that rely on fidelity metrics derived from reproducing kernels on varifolds; this allows to extend the two models to the comparison of varifolds in different orbits. We also obtain the existence of solutions to the corresponding optimal control problems under adequate assumptions and derive optimality conditions for those solutions. Based on these, in Section \ref{sec:implementation_results}, we propose a numerical approach to estimate the optimal matching between discrete varifolds and illustrate it on several simple examples, in particular for the registration of shapes with partially missing data, which are meant to serve as proof-of-concept of the validity of the model. Our Python implementation is also made openly available on Github \footnote{\url{https://github.com/charoncode/Var_metamorph}}. For the purpose of readability and concision of the main text, we have grouped all proofs of the theorems and propositions in the Appendix.

\section{Weight metamorphoses on varifolds}
\label{sec:weight_metam}
\subsection{Diffeomorphic varifold transformation}
The model we shall study in this paper builds on the varifold diffeomorphic registration approach that the authors had introduced in their previous works \cite{hsieh2020diffeomorphic,hsieh2021metrics}. In the following paragraphs, we give a brief summary of the general framework of those papers and thereby introduce some notations and definitions that will be necessary for the upcoming sections. 

In all the paper, we shall consider the Euclidean space $\R^n$ with $n\geq 2$ to be the ambient space in which our ``shapes'' of interest live. The Euclidean inner product between two vectors $a,b \in \R^n$ will be written $a\cdot b$ or $a^Tb$. For any integer $1\leq r \leq \infty$, a $C^r$ \textit{diffeomorphism} of $\R^n$ is a bijective map $\R^n \rightarrow \R^n$ of class $C^r$ such that the inverse is also $C^r$. These will constitute our set of geometric deformations as we shall detail further below. In addition to the Euclidean space, another important set for the rest of the paper is the $d$-dimensional \textit{oriented Grassmannian} for $0\leq d \leq n$ which we denote by $\origrass$ and which is defined as the set of all $d$-dimensional oriented linear subspaces of $\R^n$. The oriented Grassmannian carries a natural manifold structure as it can be identified with the quotient of groups $SO(n)/(SO(d)\times SO(n-d))$. Alternatively, one can think of an element $U\in \origrass$ as the equivalence class of all the oriented frames $(u_1,\ldots,u_d) \in (\R^n)^d$ that span $U$ with the correct orientation. In particular, for $d=1$ or $d=n-1$, we can identify $\origrass$ with the unit sphere $\Sp^{n-1}$. Furthermore, there is a natural metric on $\origrass$ that is inherited from this frame representation. Given $U, U' \in \origrass$ and $(u_1,\ldots,u_d)$ and $(u_1',\ldots,u_d')$ representative frames of those oriented spaces, it is defined by:
\begin{equation}\label{eq:inner_product_Grass}
    \langle U,U' \rangle = \det(u_i \cdot u_j')_{i,j=1,\ldots,d} 
\end{equation}

We can now define the central mathematical object of this paper, namely \textit{oriented varifolds} which correspond to an oriented version of the classical notion of varifold introduced within the field of geometric measure theory in the seminal works of \cite{Almgren,Allard}. Specifically, 
\begin{definition}
 An oriented $d$-varifold $\mu$ on $\mathbb{R}^n$ is a nonnegative finite Radon measure on the space $\mathbb{R}^n \times \origrass$. We denote by $\spvard$ the space of all oriented $d$-varifolds.
\end{definition}
In the rest of the paper, with a slight abuse of vocabulary, we will use the word varifold instead of oriented varifold to keep the denomination short but note that usual varifolds in the sense of \cite{Allard,Simon} simply result in replacing the oriented Grassmannian in the above definition by its unoriented counterpart. From the Riesz representation theorem, we can equivalently view any varifold $\mu$ as a distribution, i.e. an element of the dual space $C_0(\mathbb{R}^n \times \widetilde{G}^n_d)^*$, where $C_0(\mathbb{R}^n \times \origrass)$ denotes the set of continuous functions vanishing at infinity on $\mathbb{R}^n \times \origrass$. It is defined for any test function $\omega \in C_0(\mathbb{R}^n \times \origrass)$ by:
\begin{equation}
\label{eq:var_mu_distribution}
    (\mu \vert \omega) \doteq \int_{\mathbb{R}^n \times \origrass} \omega(x,U) d\mu(x,U).
\end{equation}
The \textit{weight} of a varifold $\mu \in \spvard$ is the finite Radon measure $\vert\mu\vert$ on $\R^n$ defined by $\vert\mu\vert(A) := \mu(A\times \widetilde{G}^n_d)$ for all Borel subset $A$ of $\R^n$. As a consequence of the disintegration theorem for measures on product spaces (c.f. \cite{Ambrosio2000} Chap. 2), any varifold $\mu$ can be decomposed into its weight measure and a family of probability measures on the oriented Grassmannian, namely: 
\begin{proposition}
\label{prop:disintegration}
 Let $\mu \in \spvard$. For $\vert\mu\vert$-almost every $x$ in $\R^n$, there exists a probability measure $\nu_x$ on $\origrass$ such that $x \mapsto \nu_x$ is $\vert\mu\vert$-measurable and $\mu = \vert\mu\vert \otimes \nu_x$ meaning that for all $\omega \in C_0(\mathbb{R}^n \times \origrass)$
\begin{equation}
\label{eq:var_mu_disintegration}
    (\mu \vert \omega) = \int_{\R^n} \int_{\origrass} \omega(x,U) d\nu_x(U) d\vert\mu\vert(x).
\end{equation}
\end{proposition}

Varifolds can be transformed by the action of diffeomorphisms in particular via the notion of pushforward. If $\mu \in \spvard$ and $\phi$ is a $C^1$-diffeomorphism, the \textit{pushforward} of $\mu$ by $\phi$ is the varifold $\phi_\sharp \mu \in \spvard$ such that for all $\omega \in C_0(\mathbb{R}^n \times \origrass)$:
\begin{equation}
\label{eq:def_pushforward}
    (\phi_{\#} \mu \vert \omega) \doteq \int_{\mathbb{R}^d \times \origrass} \omega(\phi(x),d_{x} \phi \cdot U) J_U \phi(x) d \mu(x,U).
\end{equation}
In the above, $d_x \phi$ is the Jacobian of $\phi$ at $x$ and $d_{x} \phi \cdot U$ denotes the oriented subspace obtained by transporting $U$ by the linear map $d_x \phi$ i.e. if $(u_1,\ldots,u_d)$ is an oriented frame spanning $U$, $d_{x} \phi \cdot U$ is the oriented subspace spanned by the frame $(d_x \phi(u_1),\ldots, d_x\phi(u_d))$. Finally $J_U \phi(x)$ is the Jacobian determinant of $\phi$ at $x$ along the subspace $U$ i.e. the change of $d$-volume induced by $\phi$ along $U$ which is given precisely by $J_U \phi(x) =\sqrt{\det(d_x \phi(u_i) \cdot d_x \phi(u_j))_{i,j}}$ if $(u_1,\ldots,u_d)$ is any orthonormal frame of $U$. As we shall explain just below, this seemingly convoluted definition of varifold pushforward extends the classical diffeomorphic transformation of submanifolds. 

Varifolds provides a representation that embeds a very wide class of mathematical structures among which usual densities and discrete measures. In $\spvard$, we will write $\delta_{(x,U)}$ a Dirac mass located at $x \in \R^n$ with attached oriented subspace $U \in \origrass$. For such a Dirac varifold, we have as a particular case of the above definitions $\vert\delta_{(x,U)}\vert = \delta_x$ and for any diffeomorphism $\phi$, $\phi_\# \delta_{(x,U)} = J_U \phi(x) \delta_{(\phi(x),d_x\phi \cdot U)}$. But beyond densities and Dirac masses, varifolds further encompass geometric structures such as submanifolds or rectifiable subsets of $\R^n$. Indeed, let $X$ be a $d$-dimensional oriented rectifiable subset of $\R^n$. As this will not be of critical importance for the rest of the paper, we refer the reader to \cite{Simon} or \cite{hsieh2021metrics} for precise definitions of rectifiable sets; otherwise the reader unfamiliar with those notions can instead restrict $X$ to be an oriented $d$-dimensional smooth submanifold of $\R^n$. Then $X$ can be naturally represented by a varifold $\mu_X \in \spvard$ defined by:
\begin{equation}
\label{eq:def_rect_var}
    (\mu_X \vert \omega) \doteq \int_{X} \omega(x,T_x X)\, d\vol_X(x)
\end{equation}
where $T_x X \in \origrass$ denotes the oriented tangent space to $X$ at $x$ and $\vol_X$ the $d$-volume measure of $X$. We will refer to such varifolds associated to rectifiable subsets through \eqref{eq:def_rect_var} as \textit{rectifiable varifolds}. In this case, the disintegration given by Proposition \ref{prop:disintegration} is given more specifically by $\vert\mu_X\vert = \vol_X$ and $\nu_x = \delta_{T_x X}$. In addition, one can check by direct application of the area formula that for any diffeomorphism $\phi$, the pushforward varifold $\phi_{\#} \mu_X$ is nothing but the varifold $\mu_{\phi(X)}$ associated to the deformed set $\phi(X)$. Lastly, we conclude this brief review by pointing out that, beyond its interest for the shape analysis problems we consider here, this representation of rectifiable sets as varifolds can be also very useful in computational geometry, for example in the estimation of discrete curvatures \cite{buet2018discretization} and curvature flows \cite{buet2020mean}. 

\subsection{Weight change model}
%introduce the extended transformation and group action, maybe discuss the orbits under this action...
Although the model discussed above does allow transformation of mass through Jacobians of deformations, in many situations considerable  inconsistencies or density variations between measures cannot be fully described by diffeomorphic transformations. We will thus extend the diffeomorphic varifold transformation model by augmenting the pushforward action \eqref{eq:def_pushforward} with a weight or density changing process. To be concrete, we consider rescaling functions living in $\mathcal{B}(\R^n,\R_+)$, the space of positive Borel measurable functions defined on $\R^n$. $\mathcal{B}(\R^n,\R_+)$ is a group under pointwise multiplication, and each element $\alpha$ can be applied to varifolds $\mu \in \spvard$ via the action defined by:
\begin{align}\label{eq:resc_fun}
    (\alpha \mu \vert \omega) = (\mu\vert \alpha \omega) = \int_{\R^n \times \origrass} \alpha(x) \omega(x,T) d\mu(x,T).
\end{align}
Informally speaking, this action modifies the density of a varifold at each point in the space $\R^n$.

We can combine the applications of the diffeomorphism group $\textrm{Diff}(\R^n)$ and of $\mathcal{B}(\R^n,\R_+)$ on varifolds by semi-direct product. We define a homomorphism $\Psi: \textrm{Diff}(\R^n) \mapsto \textrm{Aut}(\mathcal{B}(\R^n,\R_+))$, which sends elements $\varphi$ to automorphisms $\Psi_{\varphi}(\alpha) \doteq \alpha \circ \varphi$ of $\mathcal{B}(\R^n,\R_+)$. The group $\textrm{Diff}(\R^n) \ltimes_{\Psi} \mathcal{B}(\R^n,\R_+)$ of semi-direct product between $\textrm{Diff}(\R^n)$ and $\mathcal{B}(\R^n,\R_+)$ can be defined by the following group law:
\begin{align*}
(\varphi_1, \alpha_1) \cdot (\varphi_2, \alpha_2) \doteq (\varphi_1 \circ \varphi_2, (\alpha_1 \circ \varphi_2) \alpha_2),
\end{align*}
 and it is straightforward to verify that the identity is $(\text{id},1)$ and the inverse of the element $(\varphi,\alpha)$ is $(\varphi^{-1},\frac{1}{\alpha \circ \varphi^{-1}})$. A natural left action of the group $\textrm{Diff}(\R^n) \ltimes_{\Psi} \mathcal{B}(\R^n,\R_+)$ on the space of varifolds $\spvard$ can be then defined as follows:
 \begin{align}\label{eq:action_meta}
     (\varphi,\alpha) \cdot \mu \doteq \varphi_{\#} (\alpha \mu) = (\alpha \circ \varphi^{-1}) \varphi_{\#} \mu,
 \end{align}
 which corresponds to first rescaling $\mu$ by $\alpha$ then deforming $\alpha\mu$ by $\varphi$, or equivalently first deforming $\mu$ by $\varphi$ then recaling the deformed varifold by $\alpha$ with the change of coordinate $\varphi^{-1}$.
 
This extended group action leads, as we can expect, to richer orbits than the sole diffeomorphic pushforward. For a given varifold $\mu$, one can see for instance that the reweighing function $\alpha$ may set some parts of $\mu$ to zero mass, which we will leverage for applications to partial matching or for the modelling of topological changes.

\subsection{A static $L^2$ energy model}
%mention the IPMI paper, recap the LDDMM metric on diffeomorphism, write the new cost and argue that it is not a metric
\label{ssec:L2_model}
In the spirit of Grenander's approach to construct metrics on shape spaces \cite{Grenander1993}, we may attempt to recover a notion of distance on the space of varifolds by introducing right-invariant Riemannian metrics on the transformation group and relying on the previous group action. For purely diffeomorphic transformations, this was addressed in particular through the Large Diffeomorphic Deformation Metric Mapping (LDDMM) model introduced in \cite{Beg2005} which we will recap briefly. 

Its key principle is to consider diffeomorphisms which are obtained as the flow of a vector field on $\R^n$. Specifically, let $V$ be a predefined Hilbert space of vector fields $\R^n \rightarrow \R^n$ with metric written $\|\cdot\|_V$. All throughout this paper, we will make the assumption that $V$ is continuously embedded into $C^2_0(\R^n,\R^n)$, the space of $C^2$ vector fields of $\R^n$ that vanish at infinity as well as all first and second order derivatives. In other words, there exists a constant $c_V >0$ such that for any $u \in V$, we have $\|u\|_{2,\infty} \leq c_V \|u\|_V$ where $\| \cdot \|_{2,\infty}$ denotes the sup-norm of $u$ and all derivatives up to order $2$. Under such assumptions, $V$ can be further shown (c.f. \cite{Younes2019} Chap. 8) to be a \textit{reproducing kernel Hilbert space} (RKHS) associated to a positive vector kernel which we shall write $K_V: \R^n \times \R^n \rightarrow \R^{n \times n}$. The Hilbert norm of $V$ can be also expressed based on this kernel operator, which we shall used later on. Now, considering the space $L^2([0,1],V)$ of time-dependent vector fields in $V$, for any $v \in L^2([0,1],V)$, the flow map $t\in[0,1] \mapsto \varphi_t^v$ is defined for any $t\in [0,1]$ and $x \in \R^n$ by the integral equation:
\begin{equation*}
 \varphi_t(x) = x + \int_0^t v_s \circ \varphi_s^v(x) ds
\end{equation*}
or equivalently by $\varphi_0^v =\text{Id}$ and $\partial_t \varphi_t^v = v_t \circ \varphi_t^v$. The results of \cite{Beg2005} (see also \cite{Younes2019} Chap. 7) guarantee that $G_V \doteq \{\varphi_1^v \ \vert \ v \in L^2([0,1],V) \}$ is a subgroup of the group of $C^1$-diffeomorphisms of $\R^n$ which can be equipped with the right-invariant distance such that:
\begin{equation*}
 d_{G_V}(\text{Id},\phi) = \inf \left \{ \left(\int_0^1 \|v_t\|_V^2 dt\right)^{1/2} \ \vert \ v \in L^2([0,1],V), \ \varphi_1^v = \phi \right\}.
\end{equation*}
Intuitively, $\|v_t\|_V^2$ represents the instantaneous cost of the deformation at time $t$ and the total cost of a deformation $\phi \in G_V$ is obtained by minimizing the full path energy $\int_0^1 \|v_t\|_V^2 dt$. This metric on deformations can in turn be used to measure the distance between two varifolds $\mu$ and $\mu'$ (provided they belong to the same orbit under the action of $G_V$) by finding an optimal deformation field $v$ minimizing $\int_0^1\|v_t\|_V^2 dt$ under the constraint that $(\varphi_1^v)_\# \mu = \mu'$. This is the idea underlying varifold diffeomorphic registration models considered in several earlier works such as \cite{Charon2013,kaltenmark2017general,hsieh2021metrics}. 

Now, for the weight change component, we shall first discuss a simple $L^2$ penalty on $\alpha$ which is consistent with the approach introduced in our previous work \cite{hsieh2021diffeomorphic} in the more restrictive setting of discrete varifolds. As we shall see however, this energy will not allow us to recover a real notion of metric on varifolds. Given a reference varifold $\mu$, we define $L^2(\vert \mu\vert)$ to be the space of real-valued functions on $\R^n$ which are square integrable with respect to the weight measure $\vert \mu \vert$. If $\mu'$ is in the orbit of $\mu$ for the action of $G_V \times L^2(\vert \mu \vert)$, we can define a deformation-$L^2$ discrepancy between $\mu$ and $\mu'$ as:
\begin{equation}
   \label{eq:def_LDDMM_L2_distance} 
   D_{V-\mathcal{L}^2}(\mu,\mu')^2 = \inf \left\{\frac{1}{2}\int_0^1 \|v_t\|_V^2 dt + \frac{\gamma}{2}\int_{\R^n} (\alpha(x) - 1)^2 d\vert \mu \vert(x) \right\}
\end{equation}
where the infimum is taken over all time dependent deformation fields $v \in L^2([0,1],V)$ and all weight rescaling function $\alpha \in L^2(\vert \mu \vert)$ under the constraint that $(\varphi,\alpha) \cdot \mu = \mu'$. The energy that is being minimized in \eqref{eq:def_LDDMM_L2_distance} combines the previous deformation cost with a second term measuring how $\alpha$ deviates from a base value of $1$ on $\text{supp}(\vert \mu \vert)$ with $\gamma >0$ being a balancing parameter between those two terms. Although $D_{\mathcal{D-L}^2}$ provides a relatively natural notion of discrepancy and the formulation allows for simple mathematical analysis and implementation as evidenced previously by \cite{hsieh2021diffeomorphic} and developed further in Section \ref{ssec:relaxed_L2_model}, it is quite clear that it does not define a real metric on the orbit of $\mu$. This is due to the energy in \eqref{eq:def_LDDMM_L2_distance} being fundamentally non symmetric since the second integral involves only the first of the two varifold and would thus differ if one instead goes from $\mu'$ to $\mu$. 

Note that the existence of optimal $v$ and $\alpha$ achieving the infimum in \eqref{eq:def_LDDMM_L2_distance} can be shown relatively easily but, since these are not linked to any notion of geodesic for an actual metric, we will skip that discussion and postpone the study of existence to the relaxed version of \eqref{eq:def_LDDMM_L2_distance} in Section \ref{ssec:relaxed_L2_model}. 

\subsection{The LDDMM-Fisher-Rao metric}
\label{ssec:LDDMM_FR}
%intruduce the model (in the original form and then with the change of variable), prove that we have a distance and existence of geodesics.
In order to recover a true  notion of Riemannian metric on each orbit, we actually need to also introduce a dynamical model for the weight change $\alpha$ as opposed to the static model of the previous section. This is also the idea behind the general concept of metamorphoses in shape spaces \cite{Trouve1}. From there on, similar to the deformation $\varphi_t^v$, we will thus consider a time-varying weight change function $t \mapsto \alpha_t \in \mathcal{B}(\R^n,\R_+)$ which is governed by the simple ODE $\partial_t \alpha_t = g_t \alpha_t $. Here $g_t$ can be interpreted as a growth factor and used as the control for the weight change function. The transformed varifold at $t\in [0,1]$ is now $\mu_t = (\alpha_t \circ (\varphi_t^v)^{-1}) (\varphi_t^v)_\# \mu_0$. 

To define an adequate energy on $\tilde{\eta}$, we can consider a classical and relatively natural Riemannian metric on spaces of densities: the Fisher-Rao metric studied e.g. in \cite{friedrich1991fisher}. Given a reference measure $\lambda_0$ on $\R^n$, $\lambda = \rho \lambda_0$ a measure with positive density $\rho$ with respect to $\lambda_0$ and $h:\R^n \rightarrow \R$, the Fisher-Rao metric at $\lambda$ is defined as:
\begin{equation*}
    G^{FR}_{\lambda}(h,h) = \int_{\R^n} \frac{h}{\rho}.\frac{h}{\rho} d\lambda = \int_{\R^n} \frac{h^2}{\rho} d\lambda_0.
\end{equation*}
The Fisher-Rao metric has been shown to satisfy very specific properties in particular when it comes to diffeomorphism invariance \cite{bauer2016uniqueness} and is also used as a penalty in several versions of unbalanced optimal transport \cite{Liero2016,chizat2018interpolating}. Inspired from those works, we propose to use a similar approach to measure and penalize the cost of the weight change process. Specifically, in our context, the infinitesimal variation of $\alpha_t$ can be quantified as $G^{FR}_{\vert\mu_t\vert}(\partial_t \alpha_t \circ (\varphi_t^v)^{-1},\partial_t \alpha_t \circ (\varphi_t^v)^{-1})$ and the total energy of the path $t\mapsto \alpha_t$ is then:
\begin{equation*}
    \int_0^1 G^{FR}_{\vert\mu_t\vert}(\partial_t \alpha_t \circ (\varphi_t^v)^{-1},\partial_t \alpha_t \circ (\varphi_t^v)^{-1}) dt = \int_0^1 \int_{\R^n} g_t^2 \circ (\varphi_t^v)^{-1}(x) \ d\vert \mu_t\vert(x).
\end{equation*}
Combined with the LDDMM deformation energy already described in Section \ref{ssec:L2_model}, we could then define the squared distance between two varifolds $\mu_0$ and $\mu'$ as the minimizer of the total cost:
\begin{equation}
\label{eq:ener_LDDMM_FR_original}
    \frac{1}{2} \int_0^1 \|v_t\|_V^2 dt + \frac{\gamma}{2} \int_0^1 \int_{\R^n} g_t^2 \circ (\varphi_t^v)^{-1}(x) \ d\vert \mu_t\vert(x).
\end{equation}
subject to $\mu' = \mu_1 = (\alpha_1 \circ (\varphi_1^v)^{-1}) (\varphi_1^v)_\# \mu_0$. The key difference with the static $L^2$ metric model of section \ref{ssec:L2_model} is that the energy of an instantaneous change in weight now evolves alongside the current transformed measure $\mu_t$ as opposed to freezing it to the initial $\mu_0$.

We shall in fact consider a formulation of the above distance that is formally equivalent but simpler to handle technically in view of the derivations that will follow. It is obtained by introducing the change of variable $\tilde{\alpha}_t = \sqrt{\alpha_t}$ which gives after differentiation $\partial_t \tilde{\alpha}_t = \frac{1}{2} \eta_t$ with $\eta_t = g_t \sqrt{\alpha_t}$ and \eqref{eq:ener_LDDMM_FR_original} can be rewritten with respect to $\eta$ as:
\begin{equation}
\label{eq:ener_LDDMM_FR}
    C_{\mu_0}(v,\eta) = \frac{1}{2} \int_0^1 \|v_t\|_V^2 dt + \frac{\gamma}{2} \int_0^1 \int_{\R^n} \eta_t^2 \circ (\varphi_t^v)^{-1}(x) \ d\vert (\varphi_t^v)_{\#} \mu_0 \vert(x) dt.
\end{equation}
Using the definition of $(\varphi_t^v)_{\#} \mu_0$ and the change of variable $x \mapsto \varphi_t^v(x)$ in the above integral, one also obtains the following equivalent expression of the cost:  
\begin{equation}
\label{eq:ener_LDDMM_FR2}
    C_{\mu_0}(v,\eta) = \frac{1}{2} \int_0^1 \|v_t\|_V^2 dt + \frac{\gamma}{2} \int_0^1 \int_{\R^n \times \origrass} \eta_t^2(x) J_U \varphi_t^v(x) \ d\mu_0(x,U) dt.
\end{equation}
Let us now write the precise definition of the LDDMM-Fisher-Rao (LDDMM-FR) metric between varifolds. For any $\mu \in \spvard$, define the \textit{orbit} of $\mu$ as:
\begin{align}
 \label{eq:orbit_LDDMM_FR}
 \Theta(\mu) = &\big\{((\tilde{\alpha}_1^{\eta})^2\circ (\varphi_1^v)^{-1}) (\varphi_1^v)_{\#} \mu \ \text{s.t } (v,\eta) \in L^2([0,1],V \times L^2(\vert \mu \vert) \nonumber\\ 
 &\tilde{\alpha}_t^{\eta}(x)>0 \ \text{for all } t\in[0,1] \ \text{and } \vert\mu\vert-a.e \ x \in \R^n \big\}
\end{align}
in which $\tilde{\alpha}_t^{\eta}$ is by definition the solution of $\partial_t \tilde{\alpha}_t^\eta = \frac{1}{2} \eta_t$ with $\tilde{\alpha}_0^\eta = 1$ on $\text{supp}(\vert \mu \vert)$ i.e. $\tilde{\alpha}_t^\eta(x) = 1 + \frac{1}{2} \int_0^t \eta_s(x) ds$ for all $x \in \text{supp}(\vert \mu \vert)$. Note that in the definition of $\Theta(\mu)$, we impose the constraint $\tilde{\alpha}_t^{\eta}(x)>0$ meaning that we do not allow cancellation of mass for varifolds in $\Theta(\mu)$. In contrast, we also define the \textit{extended orbit} of $\mu$ as:
\begin{align}
 \label{eq:orbit_LDDMM_FR2}
 \bar{\Theta}(\mu) = \big\{&((\tilde{\alpha}_1^{\eta})^2\circ (\varphi_1^v)^{-1}) (\varphi_1^v)_{\#} \mu \ \text{s.t } (v,\eta) \in L^2([0,1],V \times L^2(\vert \mu \vert) \nonumber\\ 
 &\tilde{\alpha}_t^{\eta}(x)\geq 0 \ \text{for all } t\in[0,1], \ \vert\mu\vert-a.e \ x \in \R^n \big\}.
\end{align}
Now, let $\mu_0 \in \spvard$ and $\mu_1 \in \bar{\Theta}(\mu)$. We define the LDDMM-FR distance between $\mu_0$ and $\mu_1$ through the following optimal control problem on $(v,\eta)$:
\begin{equation}
   \label{eq:def_LDDMM_FR_distance} 
   d_{V-\mathcal{FR}}(\mu_0,\mu_1)^2 = \inf_{(v,\eta) \in L^2([0,1],V \times L^2(\vert \mu_0 \vert)} C_{\mu_0}(v,\eta)
\end{equation}
where the cost $C_{\mu_0}(v,\eta)$ is given by \eqref{eq:ener_LDDMM_FR}, $\mu_t = (\tilde{\alpha}_t^{\eta})^2\circ (\varphi_t^v)^{-1}) (\varphi_t^v)_{\#} \mu$ and the control system is:
\begin{equation*}
 \left\lbrace\begin{aligned}
\partial_t \varphi_t^v &= v_t \circ \varphi_t^v \\
\partial_t \tilde{\alpha}_t^{\eta} &= \frac{1}{2} \eta_t
\end{aligned}\right.
\end{equation*}
In what follows, we shall call a path $\mu_t$ (or $(\varphi^v_t,\tilde{\alpha}^{\eta}_t)$) obtained from the above system of ODEs an \textbf{admissible} path.

Despite the fact that $d_{V-\mathcal{FR}}$ is built by combining two Riemannian metrics on the geometric deformation and the weight change function respectively, we emphasize that the two are coupled in particular through the presence of the Jacobian change in the second term of $C_{\mu_0}$. Therefore it is not yet clear that one recovers an actual distance between varifolds nor that optimal controls $v$ and $\eta$ (and thus geodesics) exist. Those questions are addressed by the following theorems.

\begin{theorem}
\label{thm:distance_LDDMM_FR}
For any $\mu_0 \in \spvard$, the function $d_{\mathcal{D-FR}}$ given by \eqref{eq:def_LDDMM_FR_distance} induces a distance on $\Theta(\mu_0)$, i.e. it is symmetric, satisfies the triangle inequality and $d_{\mathcal{D-FR}}(\mu',\mu'')=0$ if and only if $\mu'=\mu''$.
\end{theorem}
The proof is provided in Appendix. We point out that it is here essential to restrict to $\Theta(\mu_0)$ since the symmetry and triangle inequality would not hold anymore in the extended orbit $\bar{\Theta}(\mu_0)$. The second question is whether one can recover an optimal control and thereby a geodesic path $t \mapsto \mu_t$ for the LDDMM-FR distance between two varifolds. As a preliminary step, let us consider the problem of minimizing the cost with respect to $\eta$ when the deformation field $v$ is fixed which is the object of the following Lemma.
\begin{lemma}
\label{lemma:LDDMM_FR_fixed_deform}
 Let $\mu' \in \bar{\Theta}(\mu)$ and $v \in L^2([0,1],V)$ a fixed time-dependent vector field such that there exists $\eta_0 \in L^2([0,1],L^2(\vert \mu \vert))$ for which $\tilde{\alpha}_t^{\eta_0}\geq 0$ and $(\varphi_1^{v},(\tilde{\alpha}_1^{\eta_0})^2) \cdot \mu = \mu'$. Then the variational problem:
\begin{equation*}
 \inf_{\eta \in L^2([0,1],L^2(\vert \mu \vert))} \int_0^1 \int_{\R^n \times \origrass}\eta^2_t(x) J_U\varphi_t^v(x) \, d\mu(x,U) dt
\end{equation*}
subject to $\tilde{\alpha}_t^{\eta}(x)\geq 0$ and $(\varphi_1^{v},(\tilde{\alpha}_1^{\eta})^2) \cdot \mu = \mu'$ has a unique solution given by:
\begin{equation*}
 \bar{\eta}_t(x) = 2 \frac{\tilde{\alpha}_1^{\eta_0}(x)-1}{h_t(x)\int_0^1 1/h_s(x)ds}
\end{equation*}
%in which $\rho$ is the Radon-Nykodym derivative of the measure $|\mu'|$ with respect to the transported measure $|(\varphi_1^{v^0})_{\sharp}\mu|$ and 
in which $h_t(x) = \int_{\origrass} J_U\varphi_t^{v}(x) d\nu_x(U)$ where $\mu = \vert\mu\vert \otimes \nu_x$ denotes the disintegration of $\mu$ over $\R^n \times \origrass$. 
\end{lemma}
Interestingly, we see that for fixed $v$, one obtains a unique minimizer with respect to $\eta$ which can be expressed in closed form with respect to the resulting deformation $\varphi_t^{v}$. This can be used to prove the existence of solutions to the full optimal control problem. 
\begin{theorem}
\label{thm:geodesics_LDDMM_FR}
Let $\mu \in \spvard$ a varifold with compact support and $\mu'\in \bar{\Theta}(\mu)$. Then there exists $(v^*,\eta^*) \in L^2([0,1],V\times L^2(\vert \mu \vert))$ such that $d_{\mathcal{D-FR}}(\mu,\mu') = \sqrt{C_{\mu}(v^*,\eta^*)}$. 
\end{theorem}
The proofs of Lemma \ref{lemma:LDDMM_FR_fixed_deform} and Theorem \ref{thm:geodesics_LDDMM_FR} are provided in Appendix. Unlike the restricted problem of Lemma \ref{lemma:LDDMM_FR_fixed_deform}, the uniqueness of the minimizer for the joint optimization problem does not necessarily hold anymore.

\subsection{Explicit form of the geodesics between single Diracs} 
\label{ssec:geodes_single_Diracs}
Before we introduce a relaxed registration formulation in the next section that will be more amenable to the numerical estimation of the LDDMM-FR metric and its geodesics, it is fairly enlightening to first specify the above in the simplest situation of two single Dirac varifolds. We will focus specifically on $0$- and $1$-varifolds since, under certain assumptions, we will be able to recover explicit expression for the distance and geodesics and draw interesting comparisons between the different models. 
\vskip2ex

\textbf{$0$-varifolds.} Let us fix two single Dirac $0$-varifolds (i.e. usual measures) of $\R^n$, $\mu_0=r_0 \delta_{x_0}$ and $\mu_1 = r_1 \delta_{x_1}$ with $r_0>0$ and $r_1\geq 0$. In such a case, the action of the deformation and reweighting process are essentially acting independently on $x$ and $r$ respectively and it is quite easy to compute explicitly the optimal $v$ and $\eta$ together with the distance and the geodesic. We have specifically:
\begin{proposition}
\label{prop:geod_0_var}
The LDDMM-FR distance between the two Dirac measures $r_0 \delta_{x_0}$ and $r_1 \delta_{x_1}$ is:
\begin{equation*}
 d_{V-\mathcal{FR}}(r_0 \delta_{x_0},r_1 \delta_{x_1}) = \frac{1}{2} d_{K_V}(x_0,x_1)^2 + 2 \gamma(\sqrt{r_1}-\sqrt{r_0})^2.
\end{equation*}
where $d_{K_V}$ is the distance associated to the Riemannian metric on $\R^n$ given by $G_x(h,h)=h^T K_V(x,x)^{-1} h$. The geodesic is $r(t)\delta_{x(t)}$ with $x(t)$ being the geodesic between $x_0$ and $x_1$ for the metric $G$ and $r(t) = \left((1-t) \sqrt{r_0} + t \sqrt{r_1}\right)^2$. 
\end{proposition}
The proof is given in the appendix. Note that in the particular case of a radial scalar kernel of the form $K_V(x,y) = \rho\left(\frac{\vert x-y \vert^2}{\sigma^2} \right)$, the above reduces simply to $d_{\mathcal{D-FR}}(r_0 \delta_{x_0},r_1 \delta_{x_1}) = \frac{1}{2} \vert x_1-x_0\vert^2 + 2 \gamma(\sqrt{r_1}-\sqrt{r_0})^2$ and $x(t)=(1-t)x_0+t x_1$. In other words, the Dirac travels in straight line between its initial and final positions while the weight evolution follows the typical geodesic of the Fisher-Rao metric. 
\begin{remark}
 In comparison, the Wasserstein-Fisher-Rao metric of unbalanced optimal transport which was introduced independently in \cite{Liero2016} and \cite{chizat2018interpolating} involves a more intricate interaction between position and weight. Indeed, it is shown that the distance is given by:
 \begin{equation*}
  d_{WFR}(r_0 \delta_{x_0},r_1 \delta_{x_1}) = 2\gamma^2\left[r_0 + r_1 -2 \sqrt{r_0 r_1} \overline{\cos}\left(\frac{\vert x_1-x_0\vert}{2\gamma} \right) \right]
 \end{equation*}
where $\overline{\cos}(z) = \cos(|z|\wedge \frac{\pi}{2})$ and that there are in fact two types of geodesics: when $\vert x_1-x_0\vert>\pi \delta$, one obtains the geodesic $[t^2 r_1 \delta_{x_1} + (1-t)^2 r_0 \delta_{x_0}]$ in which no transport occur while for $\vert x_1-x_0\vert<\pi \delta$ the geodesic is a travelling Dirac $r(t) \delta_{x(t)}$ where $r(t)$ and $x(t)$ are determined by $r(t) = At^2 -2Bt + r_0$, $\dot{x}(t) =\omega_0/r(t)$ with 
\begin{equation*}
 \omega_0 \doteq 2\gamma \tau \sqrt{\frac{r_0r_1}{1+\tau^2}}, \ \ \tau \doteq \tan\left(\frac{x_1-x_0}{2\gamma} \right), \ \ A=r_1+r_0-2\sqrt{\frac{r_0 r_1}{1+\tau^2}}, \ \ B=r_0-\sqrt{\frac{r_0 r_1}{1+\tau^2}}.
\end{equation*}
\end{remark}

\vskip2ex

\begin{figure}[ht]
\begin{center}
    \begin{tabular}{cc}
        \includegraphics[width=4cm]{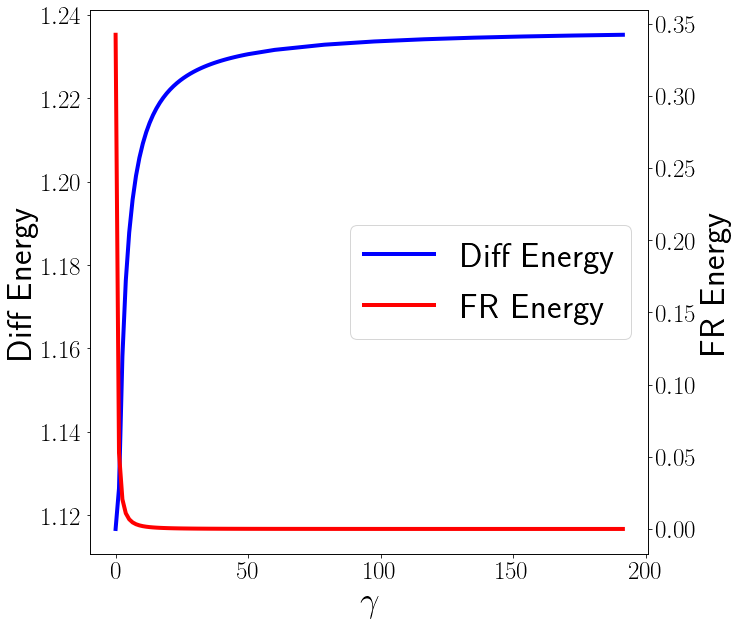}
         &\includegraphics[width=3.5cm]{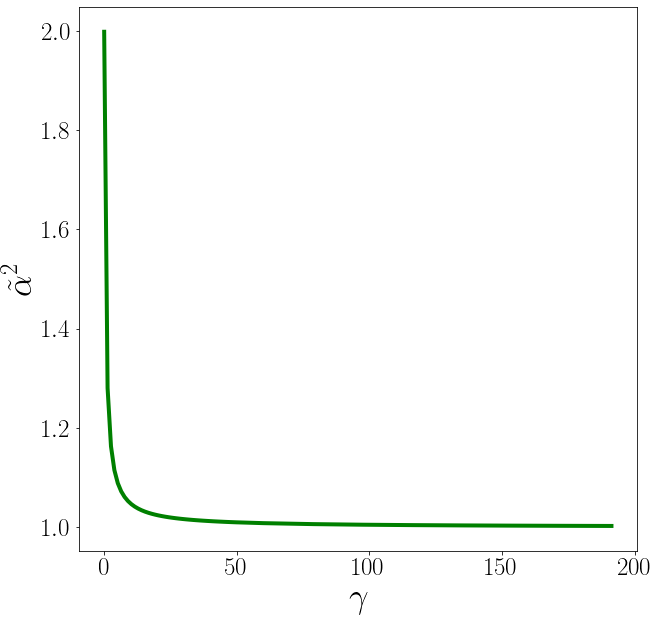}  
    \end{tabular}
    
     \begin{tabular}{ccccc}
     \rotatebox{90}{$\gamma=0.001$}
    &\includegraphics[width=2.4cm]{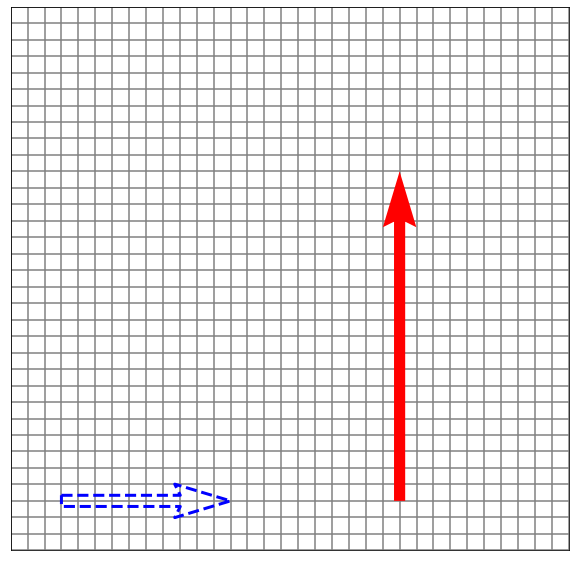}
    &\includegraphics[width=2.4cm]{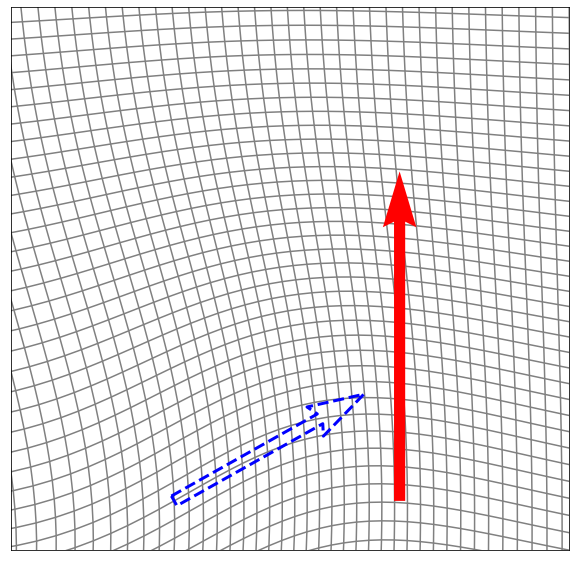}
    &\includegraphics[width=2.4cm]{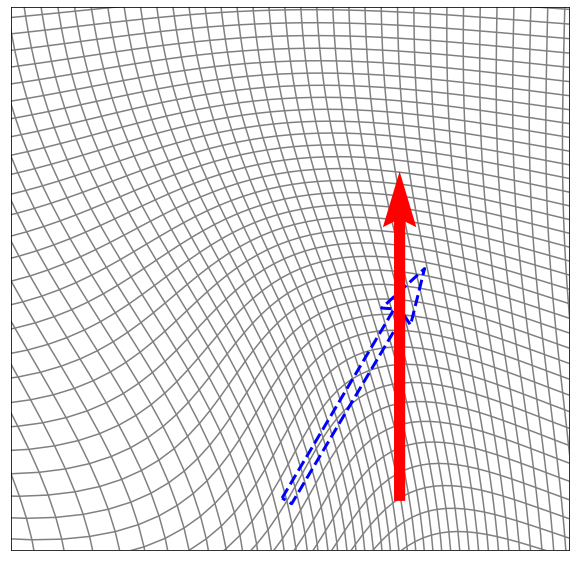}
    &\includegraphics[width=2.4cm]{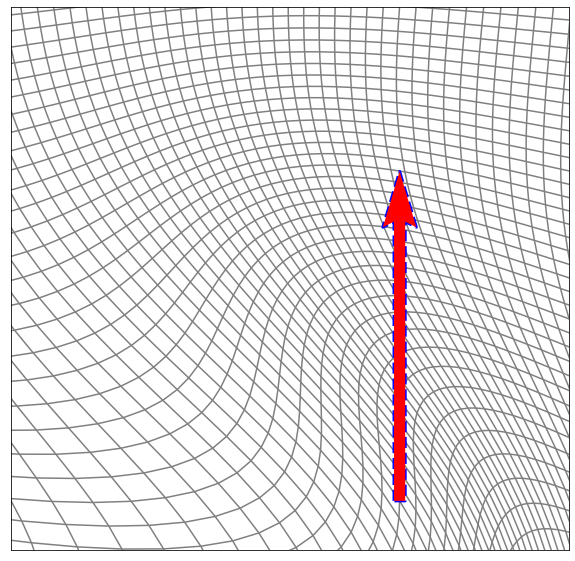}\\
    
    \rotatebox{90}{$\gamma=10$}
    &\includegraphics[width=2.4cm]{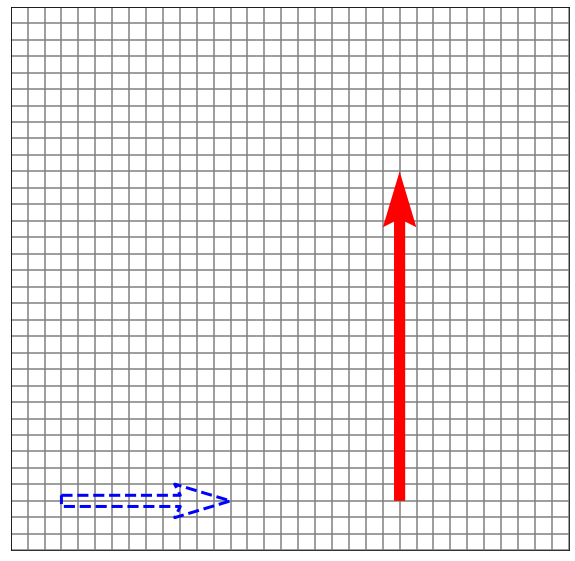}
    &\includegraphics[width=2.4cm]{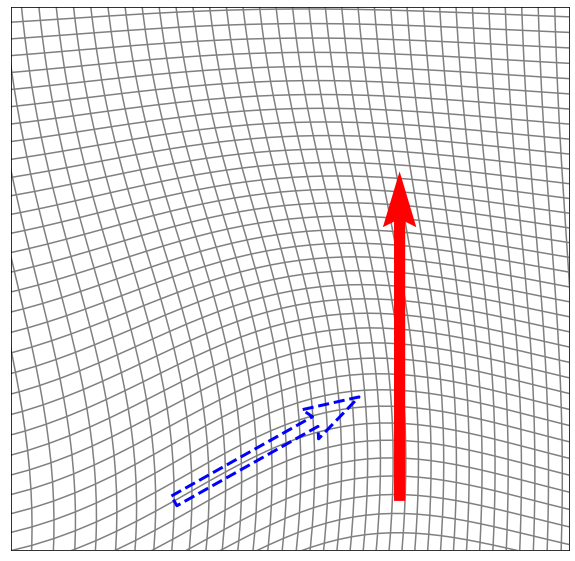}
    &\includegraphics[width=2.4cm]{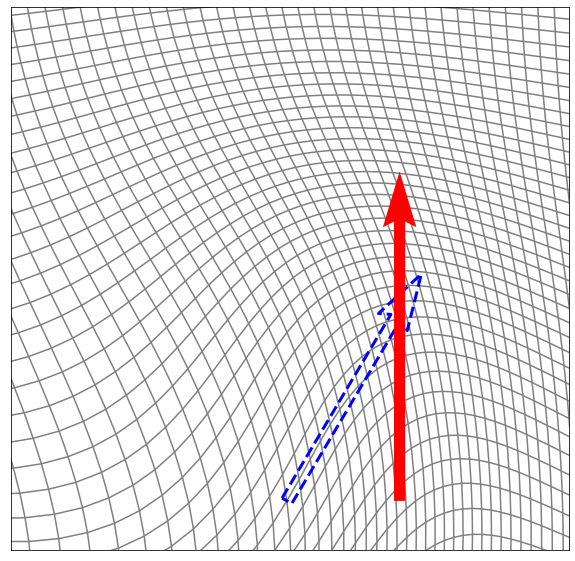}
    &\includegraphics[width=2.4cm]{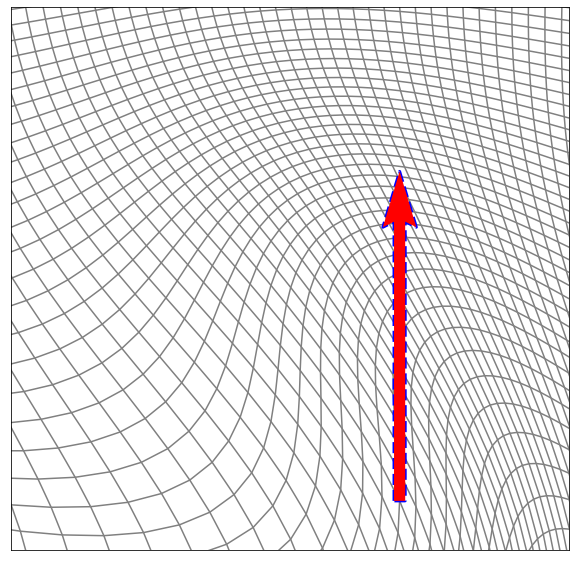}\\
    
    \rotatebox{90}{$\gamma=200$}
    &\includegraphics[width=2.4cm]{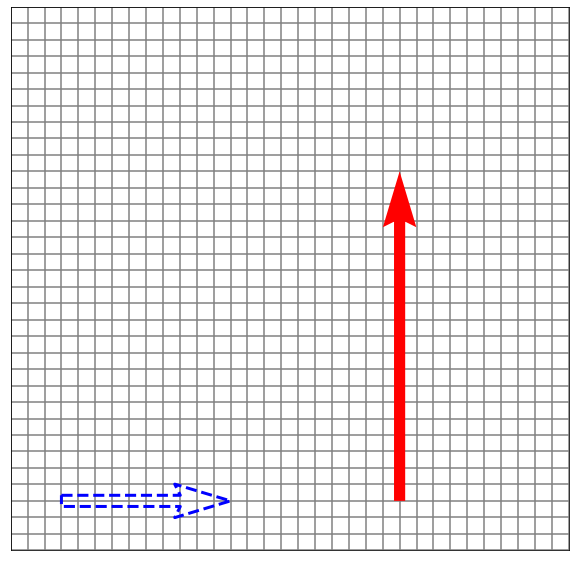}
    &\includegraphics[width=2.4cm]{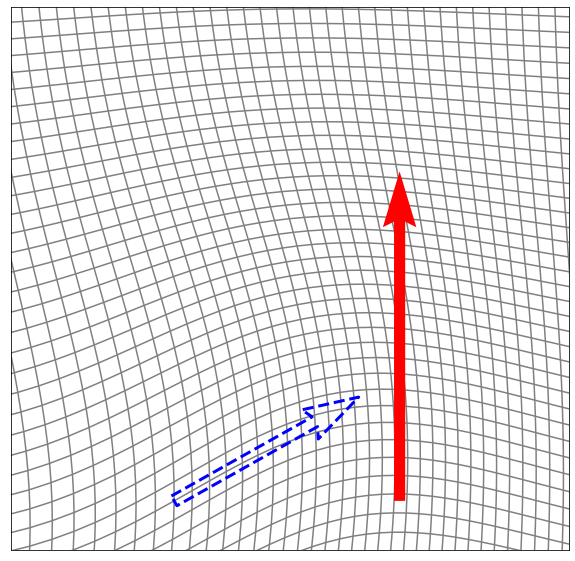}
    &\includegraphics[width=2.4cm]{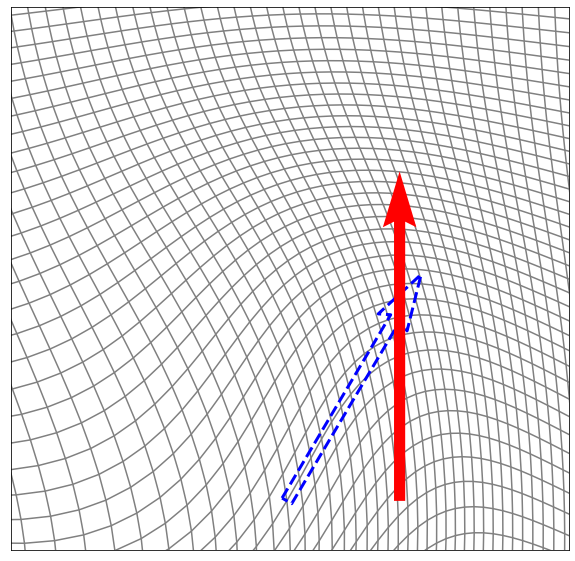}
    &\includegraphics[width=2.4cm]{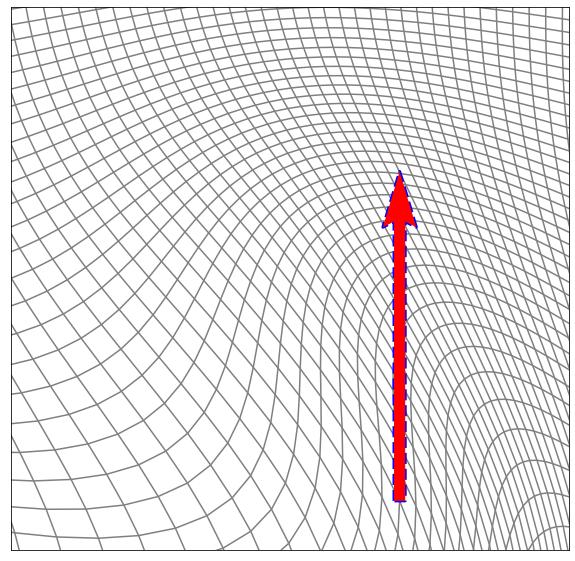}\\
    
    &$t=0$ &  $t=1/3$  &  $t=2/3$ & $t=1 $ 
    \end{tabular}   
    \caption{Exact matching between the two Diracs $\mu=r_0\delta_{(x_0,u_0)}$ and $\mu'=r_1\delta_{(x_1,u_1)}$ with $r_0=0.5$, $r_1=1$, $x_0 = (0,0),\ x_1 = (1,0) \in \R^2$ and $u_0=(1,0)$ and $u_1=(0,1)$. The figure illustrates the effect of $\gamma$ on the geodesic and optimal transformation calculated from the expression of Proposition \ref{prop:geod_1_var}. The first row shows the relative contributions of the deformation and Fisher-Rao energies as well as the estimated $(\tilde{\alpha}_1^{\eta})^2$ as functions of $\gamma$. The next rows showcase the geodesics obtained for three specific values $\gamma=0.001$, $\gamma=10$ and $\gamma=200$. A Dirac mass $r \delta_{(x,u)}$ is here represented as the vector $r u$ which foot point is at the position $x$. The transformed source varifold is plotted in blue, the target in red.} 
    \label{fig:single dirac}
\end{center}
\end{figure}

\textbf{$1$-varifolds.} As we shall see, the situation becomes quite a bit more complex in the case of 1-varifolds since unlike the 0-dimensional case, there is now a coupling between the diffeomorphic transformation and weight change through the Jacobian factor. In the following, we will directly restrict to a radial kernel $K_V(x,y) =\rho\left(\frac{\vert x-y\vert^2}{\sigma^2} \right)$ for simplicity for which we assume without loss of generality that $\rho(0)=1$. We will further assume for technical reasons that $\rho$ is $C^1$ and $\rho'(0)<0$. Let us again consider two single Dirac 1-varifolds $\mu_0 = r_0 \delta_{(x_0,u_0)}$ and $\mu_1 = r_1 \delta_{(x_1,u_1)}$ where $u_0,u_1 \in \orientgrass \approx \Sp^{n-1}$ are two unit vectors giving the directions of the lines attached to the two Diracs and $r_0>0$, $r_1 \geq 0$. We have the following result which proof can be found in the appendix:  
\begin{proposition}
\label{prop:geod_1_var}
The LDDMM-FR distance between the Dirac 1-varifolds $r_0 \delta_{(x_0,u_0)}$ and $r_1 \delta_{(x_1,u_1)}$ is:
\begin{equation*}
 d_{V-\mathcal{FR}}(r_0 \delta_{(x_0,u_0)},r_1 \delta_{(x_1,u_1)}) = \frac{\|x_1-x_0\|^2}{2} + \frac{\tau}{2} \arccos(u_0^Tu_1)^2 + 2\tau \nu^2.
\end{equation*}
where $\tau= -\sigma^2/(2\rho'(0))>0$ and $\nu$ is given by:  
\begin{equation*}
 \nu = 
 \left\lbrace\begin{aligned}
&-\frac{1}{2} \ln\left(\frac{\chi-1}{\chi+1} \right), \ \ \text{if } r_1 = 0 \\
&\ln \left[-\sqrt{\frac{r_1}{r_0}} \  \frac{1- \sqrt{1+ \frac{r_0}{r_1} (\chi^2 -1)}}{\chi-1} \right], \ \ \text{if } 0<r_1\leq r_0  \\
&\ln \left[\sqrt{\frac{r_1}{r_0}} \ \frac{1+\sqrt{1+ \frac{r_0}{r_1} (\chi^2 -1)}}{\chi+1} \right], \ \ \text{if } 0<r_0\leq r_1
\end{aligned}\right.
\end{equation*}
with $\chi= \sqrt{1+\frac{\tau}{\gamma r_0}}$. Whenever $u_1 \neq -u_0$, the geodesic is unique and of the form $r(t) \delta_{(x(t),u(t))}$ where: 
\begin{equation*} 
 \left\lbrace\begin{aligned}
&x(t)=(1-t)x_0+t x_1 \\
&u(t) = \frac{\sin((1-t)\theta) u_0 + \sin(t\theta) u_1}{\sin(\theta)} \ \ \text{if } \theta = \arccos(u_0^Tu_1)\neq 0, \ \  u(t) = u_0 \ \text{otherwise} \\
&r(t) = \left(\frac{\sqrt{r_0} \sinh((1-t)\nu) + \sqrt{r_1}\sinh(\nu t)}{\sinh(\nu)}\right)^2
\end{aligned}\right.
\end{equation*}
\end{proposition}
The above geodesic equations still have a fairly natural interpretation. The Dirac is again being transported along the straight line joining $x_0$ to $x_1$, its direction vector $u(t)$ rotates along the great circle between $u_0$ and $u_1$ at constant angular speed. However the dynamics of the weight $r(t)$ becomes more complex than with $0$-varifolds. As noted earlier, the parameter $\gamma$ essentially controls how much of the weight change is performed by the diffeomorphism itself versus the weight rescaling function. In that regard, it is interesting to look at the two limit cases $\gamma \rightarrow 0$ and $\gamma \rightarrow +\infty$. In the former case, we see that $\chi \rightarrow +\infty$ and $\nu \rightarrow 0$ which leads to the distance $d_{V-\mathcal{FR}}(r_0 \delta_{(x_0,u_0)},r_1 \delta_{(x_1,u_1)}) = \frac{\vert x_1-x_0\vert^2}{2} + \frac{\tau}{2} \arccos(u_0^Tu_1)^2$. This is as expected the same distance as for two Dirac masses with the same weight $r_0=r_1$. On the other hand, when $\gamma \rightarrow +\infty$, the model reduces to the pure diffeomorphic matching of varifolds introduced in \cite{hsieh2021metrics}. In this case, we see that $\chi\rightarrow 1^{+}$ which gives the distance $d_{V-\mathcal{FR}}(r_0 \delta_{(x_0,u_0)},r_1 \delta_{(x_1,u_1)}) = \frac{\vert x_1-x_0\vert^2}{2} + \frac{\tau}{2} \arccos(u_0^Tu_1)^2 + \frac{\tau}{2} \ln\left(\frac{r_1}{r_0}\right)^2$. Unlike with LDDMM-FR, it becomes in that case impossible to reach the zero mass ($r_1=0$) in finite distance. We further illustrate the effect of $\gamma$ on the geodesic and optimal transformation with the numerical simulations of Figure \ref{fig:single dirac}. 

\section{Relaxed matching problem}
\label{sec:relaxed_problems}
As follows from their respective variational formulations, the estimation of $D_{\mathcal{D-L}^2}$ or $d_{\mathcal{D-FR}}$ and the associated optimal paths require solving optimal control problems with a prescribed terminal condition. This is often not directly tractable except for simple cases and not even necessarily desirable in practice. One of the reason may be that $\mu'$ does not belong to the orbit of $\mu$ under the joint action of deformations and reweighing functions or, even if it does, imperfections such as noise or segmentation inconsistencies in the data can lead to unnatural behaviour of the distances and of their geodesics if one enforces this terminal matching constraint exactly. Thus, as for many other problems in variational image and shape analysis, it is common to only enforce such a constraint through the addition of a fidelity term to the energy i.e. consider a relaxed (or inexact) matching problem. Yet, unlike images or landmarks, for objects such as measures and varifolds, deriving adequate fidelity metrics that can be nicely embedded within the type of variational problems considered here is not immediate. In fact, this issue has been the object of several different works in the past such as \cite{Glaunes2004,Charon2013,Roussillon2016,Feydy2017,kaltenmark2017general,feydy2019interpolating,charon2020fidelity}. In this paper, we shall rely on fidelity terms obtained from reproducing kernel metrics on the space of varifolds which have proved successful for that purpose in diffeomorphic registration problems. We thus give a succinct overview of their construction and properties in Section \ref{ssec:kernel_metrics} before focusing on the relaxed versions of the optimal control problems of Sections \ref{ssec:L2_model} and \ref{ssec:LDDMM_FR} for which we again study the existence of solutions. 

\subsection{Kernel metrics on varifolds}
\label{ssec:kernel_metrics}
Metrics on measure spaces derived from Reproducing Kernel Hilbert Spaces (RKHS), which are also referred to as maximum mean discrepancy in the field of statistics, provide a convenient class of fidelity terms that lead to fully explicit formulas for discrete measures. They essentially rely on the representation of measures as objects of a dual functional space. In a nutshell, to a (scalar) positive definite kernel on $\R^n \times \origrass$ with adequate regularity properties is associated, thanks to Aronszajn theorem \cite{Aronszajn1950}, a reproducing kernel Hilbert space of continuous functions on $\R^n \times \origrass$ which dual Hilbert metric will induce a (pseudo-)metric between varifolds. Specifically, we will consider tensor product kernels by relying on the following (Proposition 2 in \cite{hsieh2020diffeomorphic}): 
\begin{proposition}
\label{prop:ctru_kernel}
Let $k^{pos}$ and $k^{G}$ be continuous positive definite kernels on $\mathbb{R}^n$ and $\origrass$ respectively. Assume in addition that for any $x \in \mathbb{R}^n $, $k^{pos}(x,\cdot) \in C_0(\mathbb{R}^n)$. Then $k(x,U,x',U'):= k^{pos}(x,x') k^G(U,U')$ for all $x,x' \in \R^n$ and $U,U' \in \origrass$ defines a positive definite kernel on $\mathbb{R}^n \times \origrass$ and the RKHS $W$ associated to $k$ is continuously embedded into $C_0(\mathbb{R}^n \times \origrass)$.
\end{proposition}
Note that we use a small $k$ here as it refers to a scalar kernel in contrast with the matrix-valued kernel $K_V$ of the deformation field space introduced in Section \ref{ssec:L2_model}. Now, as any varifold can be viewed as a linear form on $C_0(\mathbb{R}^n \times \origrass)$ through \eqref{eq:var_mu_distribution}, the Hilbert norm $\|\cdot\|_W$ of $W$ induces the dual Hilbert metric:
\begin{equation}
    \label{eq:def_W_dist}
    d_{W^*}(\mu,\mu') = \| \mu'- \mu \|_{W^*} = \sup_{\|\omega\|_W \leq 1} (\mu'-\mu \vert \omega).
\end{equation}
The reason why \eqref{eq:def_W_dist} may only give a pseudo-distance rather than a distance between varifolds is because the RKHS $W$ may fail to be dense in $C_0(\mathbb{R}^n \times \origrass)$. Kernels for which this additional property does hold are called $C^0$\textit{-universal} and many examples, characterizations and constructions of $C^0$-universal kernels have been proposed, we shall refer the reader to the discussion in Section 3.2 of \cite{hsieh2020diffeomorphic} as well as \cite{Sriperumbudur10} for more details on this. Under this condition on $k$, we obtain a distance on $\spvard$ which can be shown to metrize the weak-* convergence of varifolds when restricted to specific subsets of $\spvard$, meaning that $d_{W^*}(\mu_n,\mu) \rightarrow 0$ as $n \rightarrow +\infty$ if and only if $\mu_n \stackrel{\ast}{\rightharpoonup} \mu$ i.e. for any continuous compactly supported function $\omega \in C_c(\R^n \times \origrass)$, $(\mu_n \vert \omega) \rightarrow (\mu \vert \omega)$. The precise result, which we shall need later, is the following (Corollary 1 in \cite{hsieh2021metrics}):
\begin{proposition}
\label{prop:dW_metrization}
For $M>0$ and $K$ a compact subset of $\R^n \times \origrass$, define $\mathcal{V}_{d,M,K}= \{\mu \in \spvard \ \vert \ \mu(\R^n \times \origrass) \leq M \text{ and } \text{supp}(\mu) \subset K\}$. If the kernel $k$ is $C^0$-universal then the dual RKHS distance $d_{W^*}$ metrizes the weak-* convergence of varifolds on $\mathcal{V}_{d,M,K}$. 
\end{proposition}

Yet the kernel metric $d_{W^*}$, as a dual Hilbert metric, remains essentially flat and, unlike the LDDMM-FR distance considered here, does not yield any relevant notion of geodesics. However, its key advantages is that it can be evaluated between any two varifolds $\mu$ and $\mu'$ in $\spvard$ and in closed form for given kernel $k$. Indeed, it follows from the Hilbert structure and reproducing kernel property that $d_{W^*}(\mu,\mu')^2 = \|\mu\|_{W^*}^2 - 2 \langle \mu , \mu' \rangle_{W^*} + \|\mu'\|_{W^*}^2$ with
\begin{equation*}
    \langle \mu , \mu' \rangle_{W^*} = \iint_{(\R^n \times \origrass)^2} k(x,U,x',U') d\mu(x,U) d\mu'(x',U'). 
\end{equation*}
In particular for discrete varifolds $\mu=\sum_{i=1}^N r_i \delta_{x_i,U_i}$ and $\mu=\sum_{j=1}{N'} r'_j \delta_{x'_j,U'_j}$, the above simply becomes:
\begin{equation}
\label{eq:var_kernel_inner_prod}
    \langle \mu , \mu' \rangle_{W^*} = \sum_{i=1}^{N} \sum_{j=1}^{N'} k(x_i,U_i,x'_j,U'_j) r_i r'_j. 
\end{equation}
Furthermore, by selecting kernels of higher order regularity, we see that the above expression depends smoothly on the positions, directions and weights of the respective varifolds. All these characteristics make kernel metrics $d_{W^*}$ well suited as relaxation terms for our variational problems.

We shall not discuss in many more details the different families of kernels $k^{pos}$ and $k^G$ that could be selected and the corresponding properties they induce on the distance: these questions have been examined quite thoroughly in previous publications notably \cite{kaltenmark2017general} and \cite{charon2020fidelity}. For some of the upcoming mathematical results and in all numerical applications, we shall specifically restrict $k^{pos}$ to a radial kernel $k^{pos}(x,x')=h^{pos}(\vert x-x'\vert)$ and $k^G$ to a \textit{zonal} kernel on $\origrass$ namely a function of the form $k^G(U,U') = h^G(\langle U, U'\rangle)$ for the inner product given by \eqref{eq:inner_product_Grass}.  

\subsection{Relaxed static $L^2$ metric problem}
\label{ssec:relaxed_L2_model}
%proof of existence of solutions, discrete model, write the explicit form of $\alpha$ in the discrete case.
We shall first focus on the static model of Section \ref{ssec:L2_model}. Given a varifold kernel $k$ and its RKHS $W$ as above, we will replace problem \eqref{eq:def_LDDMM_L2_distance} with the following relaxed version:
\begin{equation}
    \label{eq:relaxed_LDDMM_L2_distance}
    \inf \left\{\frac{1}{2}\int_0^1 \|v_t\|_V^2 dt + \frac{\gamma}{2} \int_{\R^n} (\alpha(x) - 1)^2 d\vert \mu \vert(x) + \frac{\lambda}{2} \|\mu_1 - \mu' \|_{W^*}^2 \right\}
\end{equation}
where the minimization is again over $v \in L^2([0,1],V)$ and $\alpha \in L^2(\vert \mu \vert)$, and with $\mu_1 = (\varphi_1^v)_{\#} (\alpha \mu_0)$ although we now only impose that $\mu_1$ approximately matches $\mu'$ as measured by their kernel distance $\|\mu_1 - \mu' \|_{W^*}$, the parameter $\lambda >0$ essentially controlling how small this distance should be. We have the following result of well-posedness for this problem which proof mainly follows the standard approach of calculus of variations and is given in Appendix \ref{appendix:relaxed_L2}:
\begin{theorem}
\label{thm:existence_relaxed_L2}
Assume that the kernel $k$ is $C_0$-universal, that $W$ is continuously embedded in $C_0^1(\mathbb{R}^n \times \origrass)$ and that $\text{supp}(\mu) \subset K$ for some compact subset $K$ of $\mathbb{R}^n \times \origrass$. Then there exist $(v,\alpha) \in L^2([0,1],V)\times L^2(\vert\mu\vert)$ achieving the infimum in \eqref{eq:relaxed_LDDMM_L2_distance}.
\end{theorem}
Note that the above embedding assumption of the RKHS $W$ into $C_0^1(\mathbb{R}^n \times \origrass)$ can be recovered quite simply by imposing some adequate regularity assumptions on the kernels $k$, as follows from e.g. the results of \cite{glaunes2014matrix}.  

In view of the implementation of this model, we shall now specify problem \eqref{eq:relaxed_LDDMM_L2_distance} to the particular situation of a discrete source varifold $\mu$ and derive the optimality conditions for the resulting finite-dimensional optimal control problem. Thus let us now assume that $\mu =\sum_{i=1}^{N} r_i \delta_{(x_i,U_i)}$ for some $N \in \mathbb{N}$, $x_i \in \R^n$, $U_i \in \origrass$ and $r_i>0$. Similarly, we shall assume that $\mu'$ is of the form $\mu' = \sum_{i=1}^{N'} r_i' \delta_{(x_i',U_i')}$. Following a similar approach as in \cite{hsieh2021metrics} and \cite{hsieh2021diffeomorphic}, we will represent each $U_i \in \origrass$ by an ordered frame of $d$ vectors $\{u_i^{(1)}, \ldots, u_i^{(d)}\}$ in $\R^n$ such that $U_i$ is the oriented $d$-dimensional space spanned by this frame and those vectors are furthermore chosen so that the $d$-volume of the corresponding parallelotope $\vert u_i^{(1)}\wedge \ldots \wedge u_i^{(d)} \vert = \sqrt{\text{det}(u_i^{(l)}\cdot u_i^{(l')})_{l,l'=1,\ldots,d}}$ is equal to $r_i$. Note that the choice of frame vectors satisfying those conditions is a priori not unique but this is not an issue here since each different term in the energy functional are independent of this choice. This allows to consider the state of the control system as the finite-dimensional variable $q = ((x_i,u_i^{(k)})_{1\leq i \leq N, 1 \leq k \leq d})$. In addition, the control $\alpha \in L^2(\vert \mu \vert)$ can be now represented more simply as the vector $\alpha = (\alpha_1,\ldots,\alpha_N) \in \R_+^N$ where each $\alpha_i$ stands as the reweighting factor for the $i$-th Dirac mass. The optimal control problem can be then rewritten as follows:
\begin{align*}
\min_{v,\alpha} C_{V-L^2}(v,\alpha) \doteq &\frac{1}{2}\int_0^1 \|v_t\|_V^2 dt + \frac{\gamma}{2} \sum_{i=1}^{N} r_i(\alpha_i-1)^2 \\
&+ \frac{\lambda}{2} \left \| \left(\sum_{i=1}^{N} \alpha_i r_i \vert u_i^{(1)}(1)\wedge \ldots \wedge u_i^{(d)}(1) \vert.  \delta_{x_i(1),U_i(1)} \right)- \mu' \right\|_{W^*}^2
\end{align*}
where $U_i(1)$ is the oriented $d$-plane spanned by $\{u_i^{(1)}(1), \ldots, u_i^{(d)}(1)\}$ and the evolution of the state $q(t) = ((x_i(t),u_i^{(k)}(t))_{1\leq i \leq N, 1 \leq k \leq d})$ at $t \in [0,1]$ is governed by the control system:
\begin{equation*}
 \left\lbrace\begin{aligned}
&\dot{x}_i(t)=v_t(x_i(t)) \\
&\dot{u}_i^{(k)}(t) = d_{x_i(t)} v_t(u_i^{(k)}(t))
\end{aligned}\right.
\end{equation*}
First, one can notice that for $v$ being fixed, the minimization with respect to $\alpha$ is a quadratic program with non-negativity constraints.  However, the problem is non-convex in the deformation field $v$. Still, necessary conditions satisfied by an optimal $v$ can be obtained from the Pontryagin maximum principle (PMP) of optimal control \cite{Pontryagin1962}. Their derivation is very similar to that of previous related models \cite{hsieh2021metrics,hsieh2021diffeomorphic} and we will skip some of the details for concision. Introducing the costate variables $p=(p_i^x,p_i^{u_k})$, the problem's Hamiltonian is:
\begin{align*}
H(p,q,v) \doteq \sum_{i=1}^N \langle p^x_i, v(x_i) \rangle + \sum_{i=1}^N \sum_{k=1}^d \langle p^{u_k}_i, d_{x_i}v(u_i^{(k)}) \rangle - \frac{1}{2} \|v\|_V^2.
\end{align*}
If $(v,\alpha)$ is a minimizer of \eqref{eq:relaxed_LDDMM_L2_distance}, the PMP leads to the existence of a costate function $t \mapsto p(t) \in H^1([0,1],\R^{nN(1+d)})$ such that $\partial_v H(p(t),q(t),v_t) =0$ i.e.:
\begin{align}\label{eq:opt_vf_L2}
v_t(\cdot) = \sum_{i=1}^N K_V(x_i(t),\cdot) p_i^{x}(t) + \sum_{k=1}^d \partial_1 K_V(x_i(t),\cdot)(u_i^{(k)}(t),p_i^{u_k}(t))
\end{align}
where $K_V$ is the kernel of the Hilbert space $V$ and $\partial_1 K_V$ denotes the differential of $K_V$ with respect to his first argument, and $p(t)$ is governed by the adjoint equations $\dot{p}(t) = - \partial_q H(p(t),q(t),v_t)$ which correspond to:
\begin{align}\label{eq:forward_L2}
\left\{ \begin{array}{l}
\dot{p}_i^{x}(t) = - d_{x_i(t)}v_t^T p_i^{x}(t)   - \sum_{k=1}^d d_{x_i(t)}^{(2)}v_t(\cdot,u_i^{(k)}(t))^T p_i^{u_k}(t) \\
\dot{p}_i^{u_k}(t) = - d_{x_i(t)} v_t^T p^{u_k}_i(t)
\end{array} \right. 
\end{align}
with the terminal condition that $p_i^x(1)$ and $p_i^{u_k}(1)$ are given by minus the derivative of the varifold term with respect to $x_i(1)$ and $u_i^{(k)}(1)$ respectively. To be more explicit, based on the expression of the kernel metric \eqref{eq:var_kernel_inner_prod}, we have $p_i^x(1) = -\partial_{x_i}g(q(1))$, $p_i^{u^{(k)}}(1) = -\partial_{u_i^{(k)}}g(q(1))$ where: 
\begin{align}
\label{eq:expr_g}
 g(q) = &\frac{\lambda}{2}\sum_{i=1}^N \sum_{j=1}^N \alpha_i \alpha_j r_i r_j \vert u_i^{(1)}\wedge \cdots \wedge u_i^{(d)}\vert.\vert u_j^{(1)}\wedge \cdots \wedge u_j^{(d)}\vert k(x_i,U_i,x_j,U_j) \nonumber\\
 &-\lambda\sum_{i=1}^N \sum_{j=1}^{N'}  \alpha_i r_i r_j' \vert u_i^{(1)}\wedge \cdots \wedge u_i^{(d)}\vert k(x_i,U_i,x_j',U_j') \ + \|\mu'\|_{W^*}^2.
\end{align}
Note that we did not expand explicitly the last term $\|\mu'\|_{W^*}^2$ since it is here only a constant with respect to $q$.
% To simplify the upcoming expressions, let us denote $\tilde{r}_i=r_i \vert u_i^{(1)}(1)\wedge \ldots \wedge u_i^{(d)}(1) \vert$. Then the first order optimality conditions on $\alpha$ can be derived, based on the explicit expression of the kernel metric \eqref{eq:var_kernel_inner_prod}, as follows:
% \begin{align*}
%     &0=\partial_{\alpha_i} C_{V-L^2}(v,\alpha) \\
%     &=\gamma r_i(\alpha_i-1) + \lambda \tilde{r}_i\left(\sum_{j=1}^{N} \alpha_j \tilde{r}_j k(x_i(1),U_i(1),x_j(1),U_j(1)) - \sum_{j=1}^{N'} r_j' k(x_i(1),U_i(1),x_j',U_j')  \right)
% \end{align*}
% which shows that the optimal $\alpha$ can be obtained from the terminal state $q(1)$ through the $N \times N$ linear system $M(q(1)) \alpha = b(q(1))$ where: 
% \begin{align*}
%     &M(q(1)) = \left[\tilde{r}_i \tilde{r}_j k(x_i(1),U_i(1),x_j(1),U_j(1)) \right]_{i,j=1,\ldots,N} + \frac{\gamma}{\lambda} \text{Diag}(r_1,\ldots,r_N) \\
%      &b(q(1)) = \left[\frac{\gamma}{\lambda} r_i + \sum_{j=1}^{N'} \tilde{r}_i r_j' k(x_i(1),U_i(1),x_j',U_j') \right]_{i=1,\ldots,N}.
% \end{align*}

\subsection{Relaxed LDDMM-FR problem}
\label{ssec:relaxed_FR_model}
Let us now introduce a corresponding inexact formulation for the estimation of the LDDMM-FR metric defined by \eqref{eq:def_LDDMM_FR_distance} and \eqref{eq:ener_LDDMM_FR2}. We can naturally use again the varifold kernel metric as data fidelity term and thus consider the minimization:
\begin{equation}
    \label{eq:relaxed_LDDMM_FR_distance}
    \inf \left\{\frac{1}{2}\int_0^1 \|v_t\|_V^2 dt +\frac{\gamma}{2} \int_0^1 \int_{\R^n \times \origrass} \eta_t^2(x) J_U \varphi_t^v(x) \ d\mu(x,U) dt + \frac{\lambda}{2} \|\mu_1 - \mu' \|_{W^*}^2 \right\}
\end{equation}
over $v \in L^2([0,1],V)$ and $\eta \in L^2([0,1],L^2(\vert\mu\vert))$, where we have $\mu_t = (\tilde{\alpha}_t^{\eta})^2\circ (\varphi_t^v)^{-1}) (\varphi_t^v)_{\#} \mu$. The study of existence of solutions for \eqref{eq:relaxed_LDDMM_FR_distance} is however more delicate than with the previous model since the Fisher-Rao penalty does not lead to the necessary weak compactness properties needed in the direct method of calculus of variations. By a slightly different argument, we are still able to prove existence although only in the particular case of a discrete varifold $\mu$ and under specific technical assumptions on the varifold kernel $k$.  
\begin{theorem}
\label{thm:existence_solutions_relaxed_LDDMM_FR}
Assume that the kernel $k$ satisfies $k(x,U,x',U') \geq 0$ and $k(x,U,x,U)>0$ for all $x,x' \in \R^n$ and $U,U' \in \origrass$ and that $W$ is continuously embedded in $C_0^1(\mathbb{R}^n \times \origrass)$. Let $\mu = \sum_{i=1}^{N} r_i \delta_{(x_i,U_i)}$ be a discrete varifold and $\mu^{tar} \in \spvard$. Then the relaxed problem \eqref{eq:relaxed_LDDMM_FR_distance} has a solution $(v,\eta) \in L^2([0,1],V) \times L^2([0,1],L^2(\vert\mu\vert))$.
\end{theorem}
The proof is detailed in Appendix \ref{appendix:relaxed_FR}. We point out that the above assumptions on the positivity of $k$ are satisfied by most of the kernels that we typically consider in numerical simulations (with the notable exception of the kernel norms related to the model of currents \cite{Glaunes2008}). Moreover, although in practice we are primarily interested in this relaxed approach for discrete varifolds as this lends itself to numerical implementation, we leave it to future investigations to extend the result of Theorem \ref{thm:existence_solutions_relaxed_LDDMM_FR} to more general varifolds $\mu$ and kernel families. 

Now, as in the previous section, we shall derive optimality conditions for minimizers of this relaxed problem. We again assume that $\mu=\sum_{i=1}^{N} r_i \delta_{(x_i,U_i)}$ and  $\mu' = \sum_{i=1}^{N'} r_i' \delta_{(x_i',U_i')}$. Similar to the above, we can represent $\eta \in L^2([0,1],L^2(\vert\mu\vert))$ more simply as a function $\eta \in L^2([0,1],\R^N)$ for each $t\in [0,1]$, we shall rewrite for simplicity $\tilde{\alpha}_t^{\eta} = (\tilde{\alpha}_i(t))_{i=1,\ldots,N} \in \R^N$, each component of the vectors being associated to one of the Dirac of $\mu$. We can then describe the state of the optimal control problem by the variable $q = ((x_i,u_i^{(k)},\tilde{\alpha}_i)_{1\leq i \leq N, 1 \leq k \leq d})$ and the control system is:
\begin{equation*}
 \left\lbrace\begin{aligned}
&\dot{x}_i(t)=v_t(x_i(t)) \\
&\dot{u}_i^{(k)}(t) = d_{x_i(t)} v_t(u_i^{(k)}(t)) \\
&\dot{\tilde{\alpha}}_i(t) = \frac{1}{2} \eta_{t,i}
\end{aligned}\right.
\end{equation*}
with the initial conditions $x_i(0)=x_i$, $(u_i^{(k)}(0))_{k=1,\ldots,d} = (u_i^{(k)})$ any frame spanning $U_i$ with $\vert u_i^{(1)}\wedge \ldots \wedge u_i^{(d)} \vert = r_i$ and $\tilde{\alpha}_i(0) = 1$. Furthermore, for the costate variable $p=(p_i^x,p_i^{u_k},p_i^{\tilde{\alpha}})$, the Hamiltonian is now given by:
\begin{align*}
H(q,p,v,\eta) &= \sum_{i=1}^{N} \langle p_i^x,v(x_i) \rangle + \sum_{i=1}^{N} \sum_{k=1}^{d} \langle p_i^{u_k}, d_{x_i}v(u_i^{(k)}) \rangle + \frac{1}{2} \sum_{i=1}^{N} p_i^{\tilde{\alpha}} \eta_i \\
&\phantom{aa}-\frac{1}{2} \|v\|_V^2 -\frac{\gamma}{2} \sum_{i=1}^{N}\vert u_i^{(1)}\wedge \cdots \wedge u_i^{(d)}\vert \eta_i^2.
\end{align*}
From this expression of the Hamiltonian and the PMP, we deduce that if $(v,\eta)$ is a minimizer then:
\begin{align*}
  &v_t(\cdot) = \sum_{i=1}^N K_V(x_i(t),\cdot) p_i^{x}(t) + \sum_{k=1}^d \partial_1 K_V(x_i(t),\cdot)(u_i^{(k)}(t),p_i^{u_k}(t)) \\
  &\eta_{t,i} =  \frac{p_i^{\tilde{\alpha}}}{2\gamma \vert u_i^{(1)}(t)\wedge \cdots \wedge u_i^{(d)}(t)\vert}
\end{align*}
with the costate function $p(t)$ satisfying the following adjoint equations:
\begin{align*}
\left\{\begin{array}{l}
\dot p_i^x(t) = -d_{x_i(t)} v_t^T p_i^x(t) - \sum_{k=1}^d d^2_{x_i(t)} v_t(\cdot,u_i^{(k)}(t)) p_i^{u_k}(t) \\
\dot p_i^{u_k}(t) = -d_{x_i(t)} v_t^T p_i^{u_k}(t) +\frac{\gamma \eta_{i,t}^2}{2 \vert u_i^{(1)}(t)\wedge \cdots \wedge u_i^{(d)}(t)\vert} \sum_{\ell=1}^{d} C_{i}^{\ell k}(t) u_i^{(k)}(t)  \\
\dot p_i^{\tilde{\alpha}}(t) = 0 \ \Rightarrow \ p_i^{\tilde{\alpha}}(t) = p_i^{\tilde{\alpha}} = \text{Cte}.
\end{array} \right.
\end{align*}
in which the $(C^{\ell k}_i(t))_{1\leq \ell,k \leq d}$ denote the coefficients of the cofactor matrix of the Gramian $(\langle u_i^{(\ell)}(t),u_i^{(k)}(t) \rangle)_{1\leq \ell,k \leq d}$. Moreover, one has again the terminal conditions on the costates $p_i^x(1) = -\partial_{x_i}g(q(1))$, $p_i^{u^{(k)}}(1) = -\partial_{u_i^{(k)}}g(q(1))$ and $p_i^{\tilde{\alpha}}(1) = -\partial_{\tilde{\alpha}_i}g(q(1))$ where $g$ is still given by the expression \eqref{eq:expr_g} with $\alpha_i$ replaced by $\tilde{\alpha}_i^2$. 

Finally, a last condition resulting from the PMP applied to this problem is the conservation of the Hamiltonian over time, namely for all $t\in[0,1]$, $H(q(t),p(t),v_t,\eta_t) = H(q(0),p(0),v_0,\eta_0)$. Furthermore, based on the above expressions of the optimal $v_t$ and $\eta_t$ and using the reproducing kernel property for the kernel $K_V$ and its derivatives (c.f. \cite{glaunes2014matrix}), one can show that:
\begin{align*}
    &H(q(t),p(t),v_t,\eta_t) = \frac{1}{2} \|v_t\|_V^2 + \frac{\gamma}{2} \sum_{i=1}^{N}\vert u_i^{(1)}(t)\wedge \cdots \wedge u_i^{(d)}(t)\vert \eta_{t,i}^2 \\
    &=\frac{1}{2}\sum_{i=1}^{N} \langle p_i^x(t),v_t(x_i(t)) \rangle + \frac{1}{2} \sum_{i=1}^{N} \sum_{k=1}^{d} \langle p_i^{u_k}(t), d_{x_i(t)}v_t(u_i^{(k)}(t)) \rangle + \frac{1}{4} \sum_{i=1}^{N} p_i^{\tilde{\alpha}} \eta_{t,i}
\end{align*}
from which we get that for the optimal $(v,\eta)$, the full transformation energy can be expressed as:
\begin{align}
\label{eq:reduced_Hamiltonian_FR}
    &\frac{1}{2} \int_0^1 \| v_t\|_V^2 dt + \frac{\gamma}{2} \sum_{i=1}^{N} \int_0^1 \vert u_i^{(1)}(t)\wedge \cdots \wedge u_i^{(d)}(t)\vert \eta_{t,i}^2 dt = H(q(0),p(0),v_0,\eta_0) \nonumber\\
    &=\frac{1}{2} \sum_{i,j=1}^N \Big\langle p_i^x(0), K_V(x_j,x_i) p_j^x(0) +\sum_{k=1}^d \partial_1 K_V(x_j,x_i)(u_j^{(k)},p_j^{u_k}(0)) \Big\rangle \nonumber \\
    &+\frac{1}{2} \sum_{i,j=1}^N \sum_{k=1}^{d} \Big\langle p_i^{u_k}(0), \partial_2 K_V(x_j,x_i)(u_i^{(k)},p_j^x(0)) \Big\rangle \nonumber\\ 
    &+\frac{1}{2}\sum_{i,j=1}^N \sum_{k,l=1}^d \Big \langle p_i^{u_k}(0), \partial_{1,2}^2 K_V(x_j,x_i)(u_i^{(k)},u_j^{(l)})p_j^{u_l}(0) \Big\rangle 
    +\frac{1}{8\gamma} \sum_{i=1}^{N} \frac{(p_i^{\tilde{\alpha}})^2}{r_i}
\end{align}

\section{Implementation and results}
\label{sec:implementation_results}
\subsection{Registration algorithms}
In order to numerically approximate solutions of the optimal control problems \eqref{eq:relaxed_LDDMM_L2_distance} and \eqref{eq:relaxed_LDDMM_FR_distance} in the case of discrete Dirac varifolds, we extend the method proposed in \cite{hsieh2021metrics}, namely we solve those problems using a shooting scheme based on the state and costate equations derived in Sections \ref{ssec:relaxed_L2_model} and \ref{ssec:relaxed_FR_model}. Our implementation in Python further leverages some recently developed libraries in order to efficiently evaluate the different functionals involved and their gradients. It is openly accessible on Github\footnote{\url{https://github.com/charoncode/Var_metamorph}}. In the next paragraphs, we detail a little more specifically the different components of our approach. 
\vskip1ex
\textbf{Optimization scheme.} In both models, we perform joint optimization over the deformation i.e., thanks to the above Hamiltonian equations, over the initial costate variables $p_i^x(0)$ and $p_i^{u^{(k)}}(0)$ and the weight changes, namely the $\alpha_i$'s in the static model and the $p_i^{\tilde{\alpha}}$'s in the LDDMM-FR case. This is done using the limited memory BFGS (L-BFGS) scheme of the SciPy library with the extra non-negativity constraints on the weight change variables. The actual computation of the cost functions and their gradients is explained below.
\vskip1ex
\textbf{Numerical integration.} Given the initial conditions $(x_i(0))$ and $(u_i^{(k)}(0))$ together with some values for the initial costates $(p_i^x(0))$ and $(p_i^{u^{(k)}}(0))$ (as well as $(p_i^{\tilde{\alpha}})$ for the LDDMM-FR model), we compute the approximate evolution of those variables for $t\in[0,1]$ by numerical integration of the coupled system of state and costate nonlinear ODEs. In our implementation, we fix a certain number $T$ of discrete time points (typically we take $T=15$ in our experiments) and use a Runge-Kutta scheme of order $4$ as numerical integrator. It is implemented through the PyTorch library in order to take advantage of the built-in back propagation pipelines and GPU computations. \vskip1ex
\textbf{Cost function evaluation and gradient computation.} The cost functions in \eqref{eq:relaxed_LDDMM_L2_distance} and \eqref{eq:relaxed_LDDMM_FR_distance} are made, on the one hand, of the transformation energies which are directly functions of the initial costates $(p_i^x(0))$, $(p_i^{u^{(k)}}(0))$ and $\alpha_i$ (or $(p_i^{\tilde{\alpha}})$) as shown in \eqref{eq:reduced_Hamiltonian_FR} and, on the other, of the fidelity term \eqref{eq:expr_g} which only depends indirectly on those costates via the terminal state $q(1)$. In both cases, our implementation simply relies on PyTorch and its CUDA bindings to compute those functions on GPU and automatically differentiate them. For the fidelity term, the gradient with respect to the final state $q(1)$ is then back propagated through the RK4 scheme as explained above to automatically recover its gradient with respect to the initial costates. We also point out that the recurrent and most numerically intensive operation in both the integration of the Hamiltonian equations as well as the computation of the varifold data fidelity term consists in the evaluation of kernel convolutions over all Diracs of the source shape. This is typically not handled efficiently in PyTorch itself as it operates by building large kernel matrices, which poses memory and time issues in the case of large sets of Diracs. We remedy this particular problem by taking advantage of the recently developed PyKeops library \cite{charlier2021kernel} that provides specialized implementations of kernel operations that remain compatible with PyTorch.
\vskip1ex
\begin{figure*}
%\fbox{\rule{0pt}{2in} \rule{.9\linewidth}{0pt}}
\begin{center}
\hspace*{-.5cm}
     \begin{tabular}{ccc}
    \includegraphics[width=3.7cm]{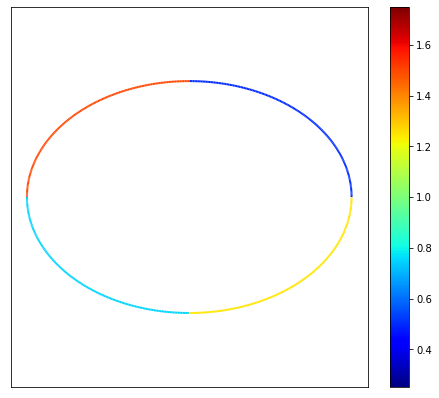}
    &\includegraphics[width=3.7cm]{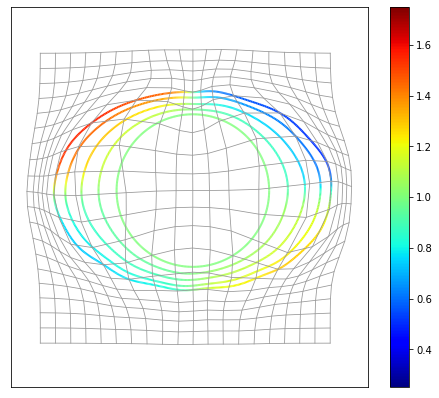}
    &\includegraphics[width=3.7cm]{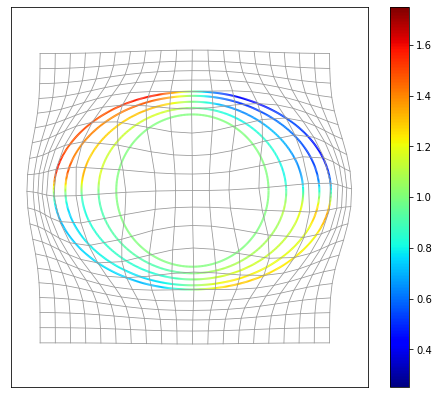} \\
    Target & LDDMM-L$^2$ & LDDMM-FR
    \end{tabular}
    \vskip3ex
    \hspace*{-.5cm}
    \begin{tabular}{ccc}
   \includegraphics[width=3.7cm]{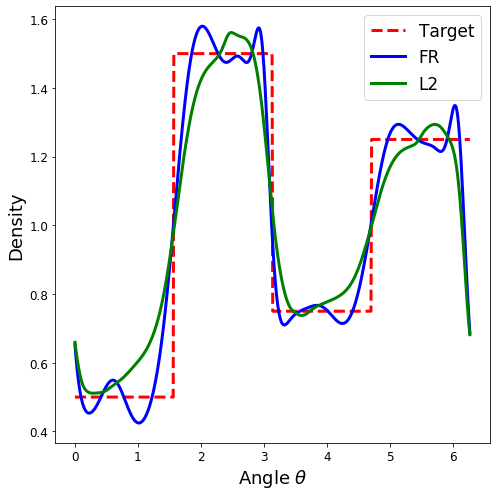} 
    &\includegraphics[width=3.7cm]{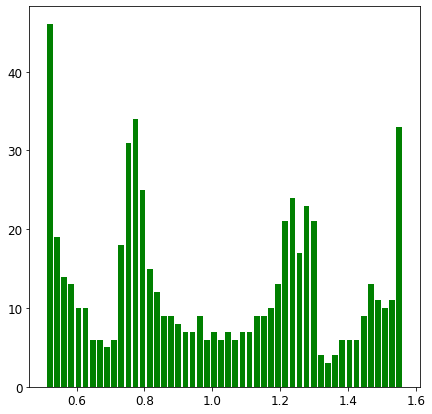} 
    &\includegraphics[width=3.7cm]{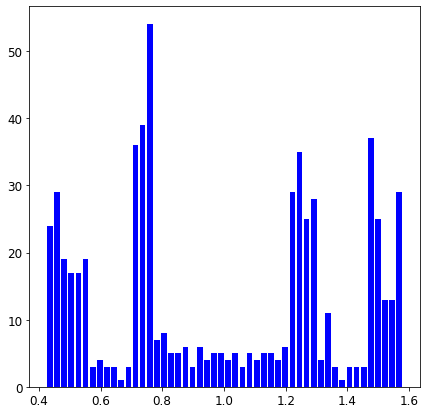} \\
    Densities &LDDMM-L$^2$ & LDDMM-FR
    \end{tabular}   
    \caption{Registration of a circle with uniform weight density $1$ to an ellipse with piecewise constant weights, $.5$, $.75$, $1.25$ and $1.75$. The first column shows the target ellipse with colors corresponding the weight at each location and the plot of those weights as functions of the angle coordinate along the curve. The second column displays the evolution of the transformed curve as well as the histogram of final weight values for the LDDMM-L$^2$ model. The corresponding plots for the LDDMM-FR model are on the third column.}
    \label{fig:circle to ellipse}
\end{center}
\end{figure*}
\textbf{Parameters.} The fundamental parameters in both registration models are the choice of deformation kernel $K_V$, of the varifold kernel $k$ defining the fidelity term and the coefficients $\gamma$ and $\lambda$ that weigh the relative importance of the different terms in the cost functional. In our implementation and the simulations of the next section, we use a Gaussian kernel $K_V(x,y) = \exp\left(-\vert x-y\vert^2/\sigma_V^2\right) \text{Id}_{n\times n}$ where $\sigma_V>0$ represents the deformation scale. As for $k$, we restrict to separable kernel of the form described at the end of Section \ref{ssec:kernel_metrics} and chosen among the different specific classes considered and discussed extensively in \cite{kaltenmark2017general,charon2020fidelity}. As in most those prior works, the scales of the different kernels are selected manually according to the typical size of the considered shapes and desired level of accuracy of the matching. Our implementation further allows for multistep and multiscale strategies to be used in particular as ways to improve the robustness of the minimization procedure. We point out that developing more data-driven and automatic selection methods for these various parameters is an important ongoing research topic. 

\subsection{Numerical results}
%colored circle ($L^2$ and F-R), amygdala (F-R, pure geometry), knee example (F-R,...)
We conclude this paper by presenting a few simulations based on the inexact registration algorithms presented above for $1$- and $2$- varifolds in $\R^n$ ($n=2$ or $3$). We stress that, in all the examples that involve curves of surfaces, these objects are first converted into corresponding $1$- or $2$- discrete varifolds to be processed by the algorithm. Specifically, we use the same discretization scheme for polygonal curves and triangulated surfaces that was introduced in earlier publications \cite{kaltenmark2017general,charon2020fidelity,hsieh2021metrics} in which each facet is approximated by a single Dirac $r_i \delta_{(x_i,U_i)}$ located at the center of mass $x_i$, with weight $r_i$ given by the area (or length) of the facet and the oriented space $U_i$ being represented by the $1$ or $2$ frame vectors spanned by the edges. The optimal costates estimated by the approach in turn allow to reconstruct the optimal deformation path $\varphi_t^v$ and transformed weights $\alpha_i r_i$. For better visualization, rather than showing the transformation of their associated varifolds, we typically plot the corresponding transformation induced on the original curve or surface with the weights being represented as changing colors over the respective facets. 

\begin{figure}
\begin{center}
    \begin{tabular}{cc}
        \includegraphics[trim= 25mm 20mm 15mm 10mm ,clip,width=3cm]{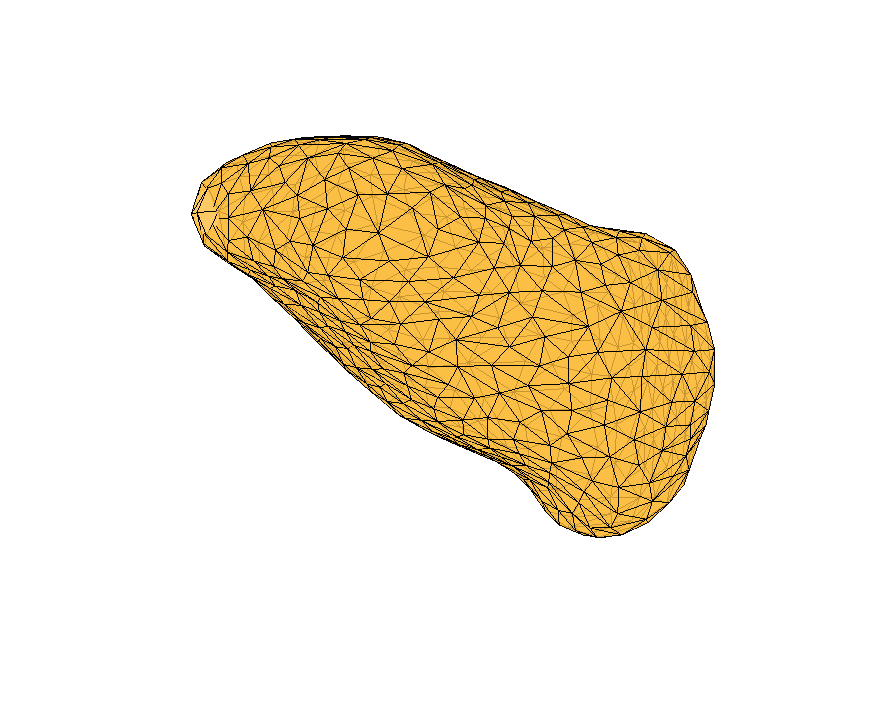}
         &\includegraphics[trim= 25mm 20mm 15mm 10mm ,clip,width=3cm]{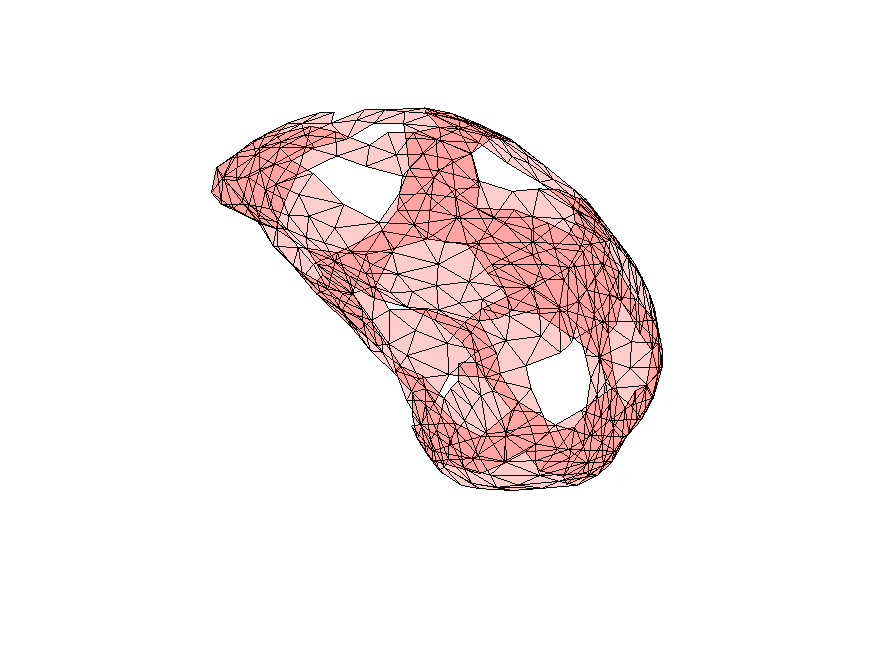} \\
         Source & Target
    \end{tabular}
    
     \begin{tabular}{cccccc}
     \rotatebox{90}{\phantom{LD}\small{LDDMM}}
    &\includegraphics[trim= 25mm 15mm 15mm 10mm ,clip, width=2.4cm]{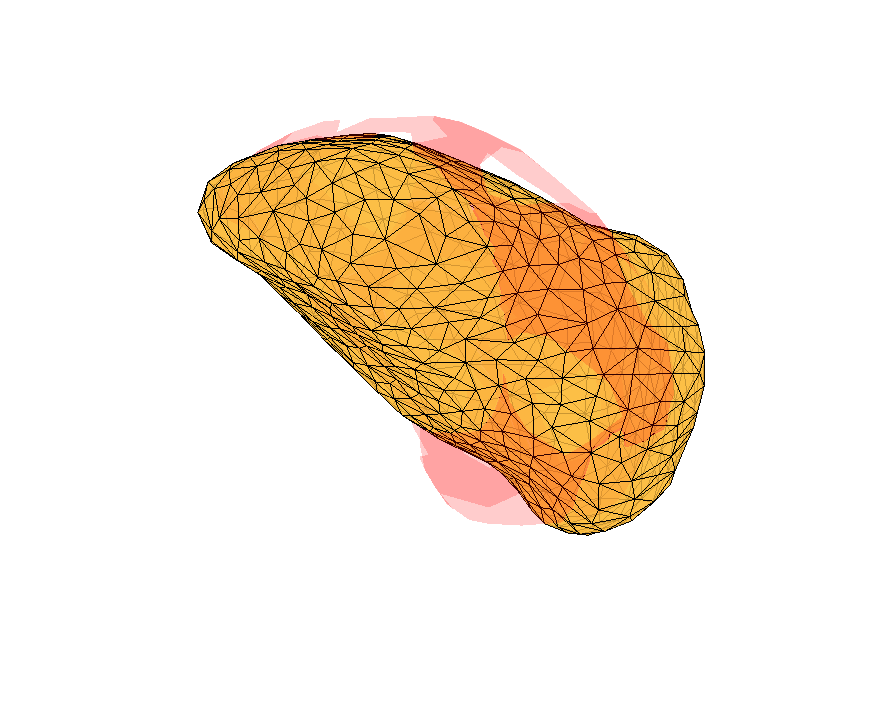}
    &\includegraphics[trim= 25mm 15mm 15mm 10mm ,clip, width=2.4cm]{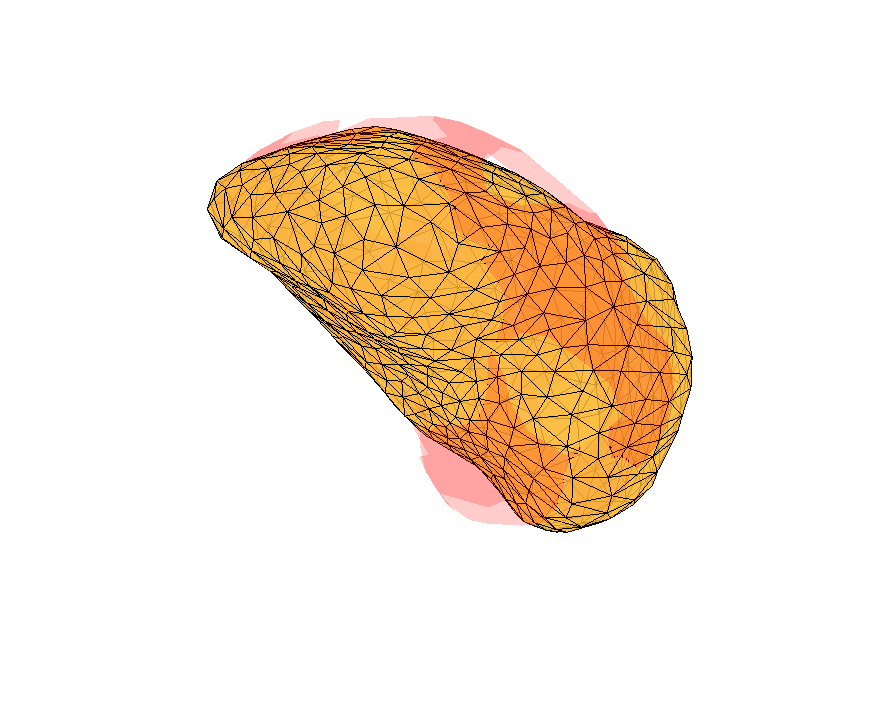}
    &\includegraphics[trim= 25mm 15mm 15mm 10mm ,clip, width=2.4cm]{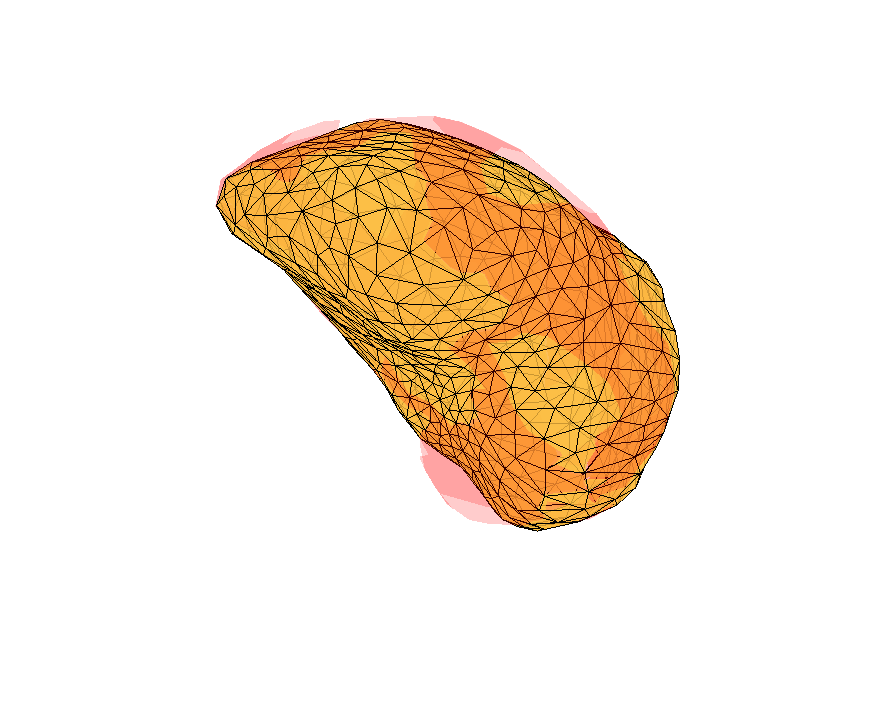}
    &\includegraphics[trim= 25mm 15mm 15mm 10mm ,clip, width=2.4cm]{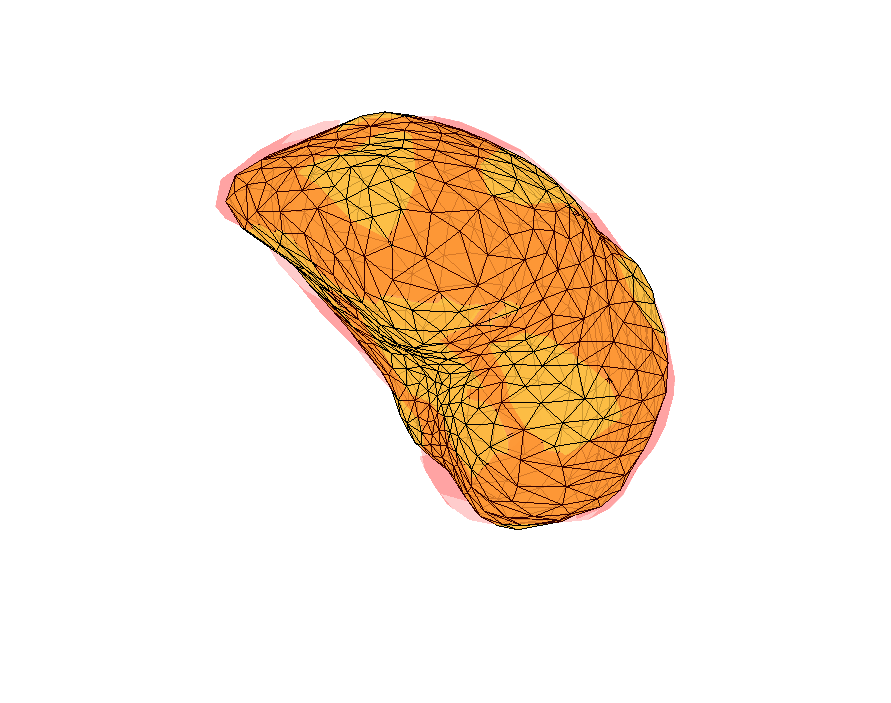}
    &\\
     
     \rotatebox{90}{\small{LDDMM-L$^2$}}
    &\includegraphics[trim= 25mm 15mm 15mm 10mm ,clip, width=2.4cm]{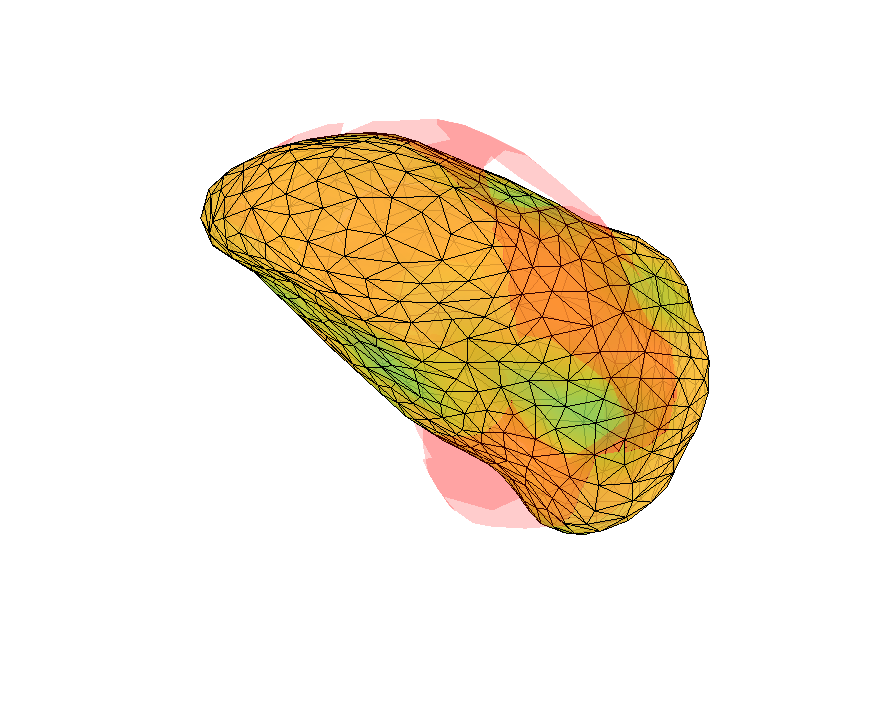}
    &\includegraphics[trim= 25mm 15mm 15mm 10mm ,clip, width=2.4cm]{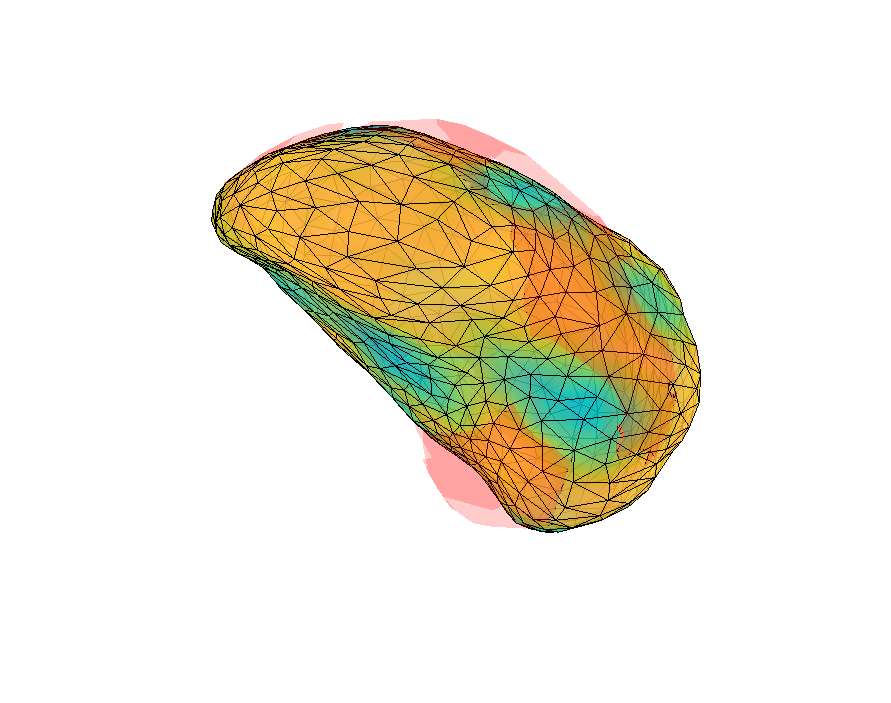}
    &\includegraphics[trim= 25mm 15mm 15mm 10mm ,clip, width=2.4cm]{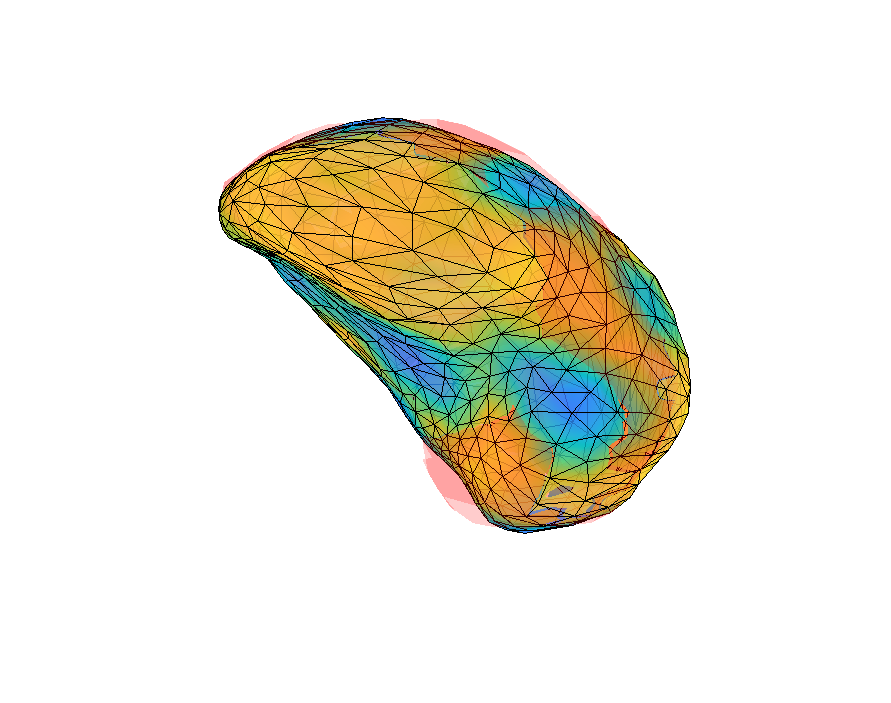}
    &\includegraphics[trim= 25mm 15mm 15mm 10mm ,clip, width=2.4cm]{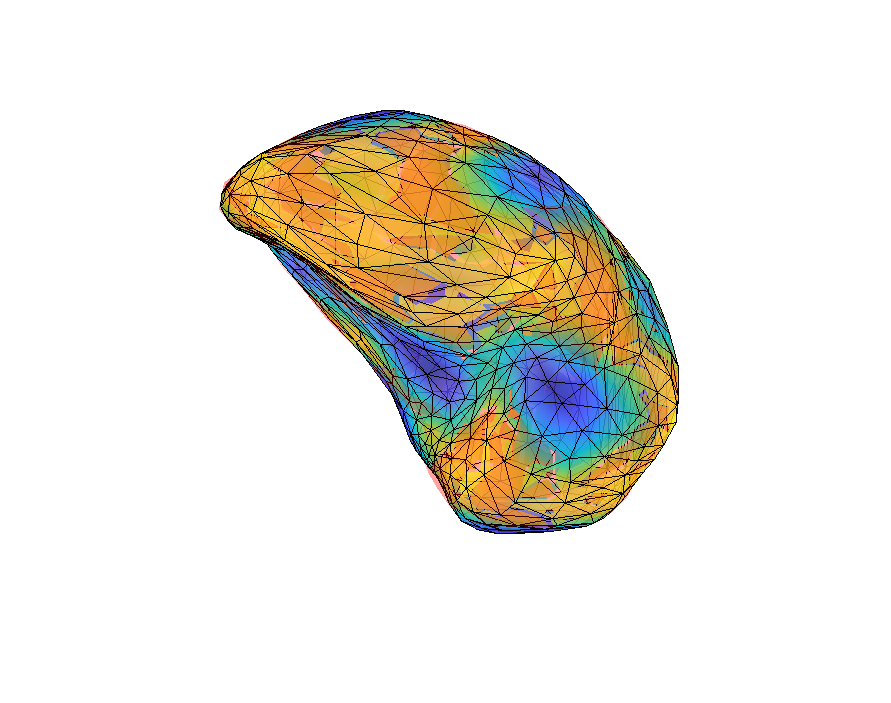}
    &\multirow{ 2}{*}[1.8cm]{\includegraphics[height=4cm]{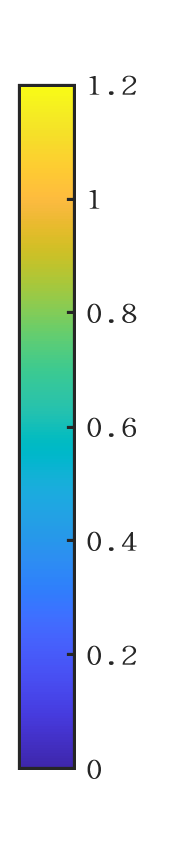}}\\
    
    \rotatebox{90}{\small{LDDMM-FR}}
    &\includegraphics[trim= 25mm 15mm 15mm 10mm ,clip, width=2.4cm]{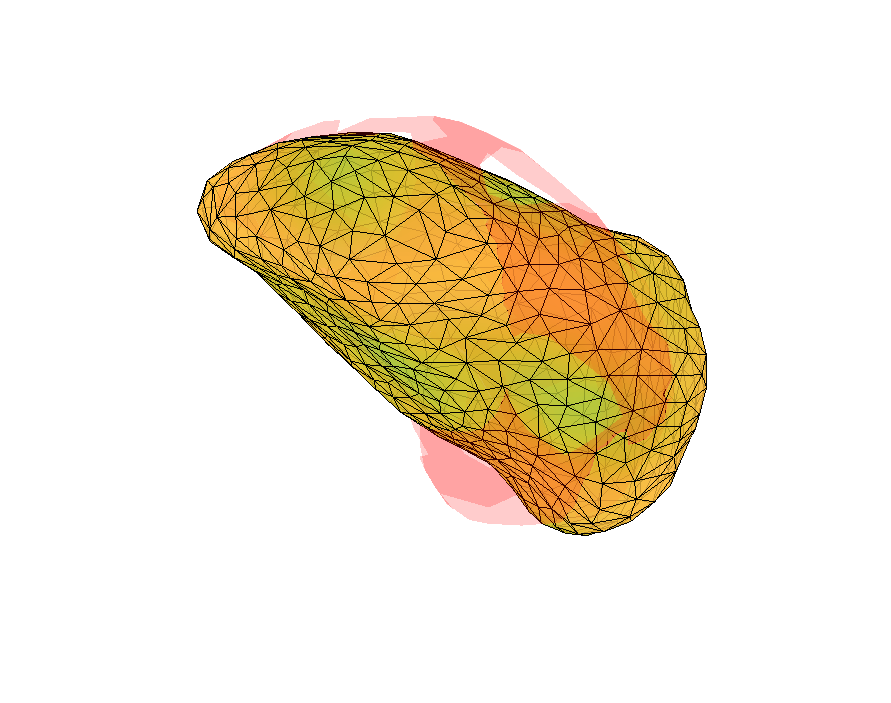}
    &\includegraphics[trim= 25mm 15mm 15mm 10mm ,clip, width=2.4cm]{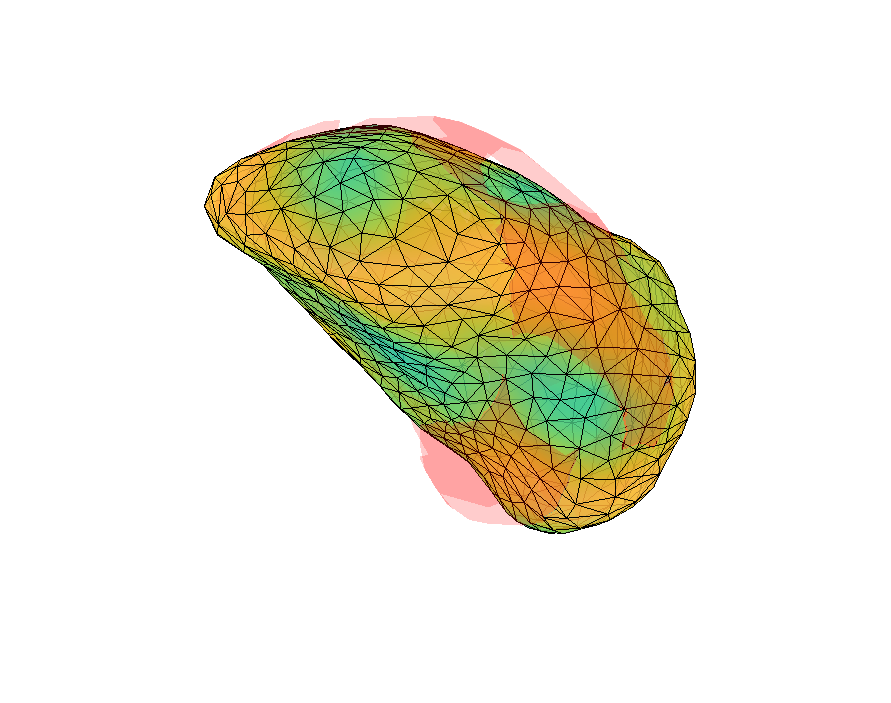}
    &\includegraphics[trim= 25mm 15mm 15mm 10mm ,clip, width=2.4cm]{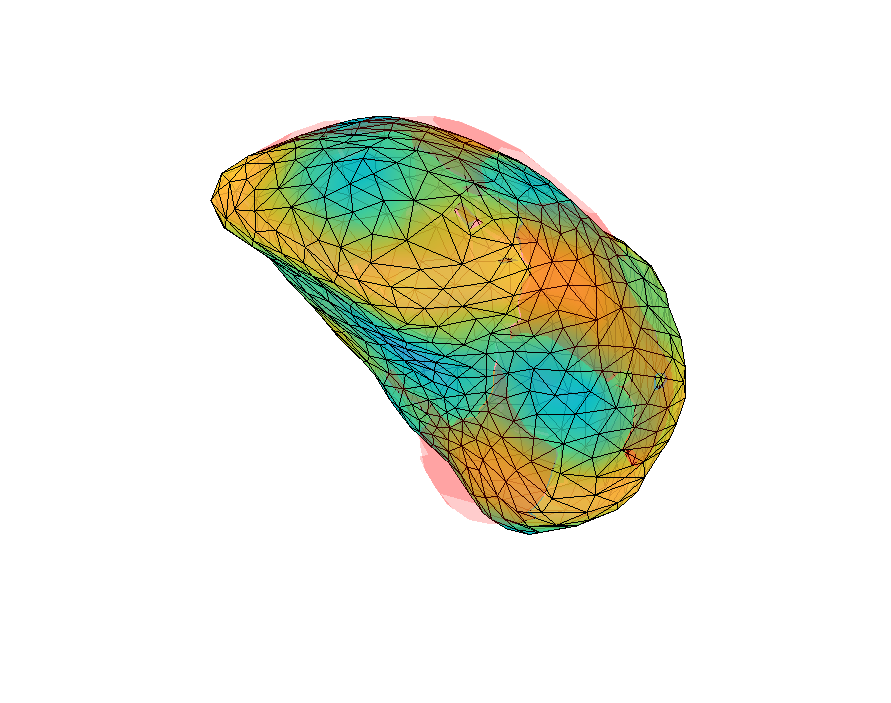}
    &\includegraphics[trim= 25mm 15mm 15mm 10mm ,clip, width=2.4cm]{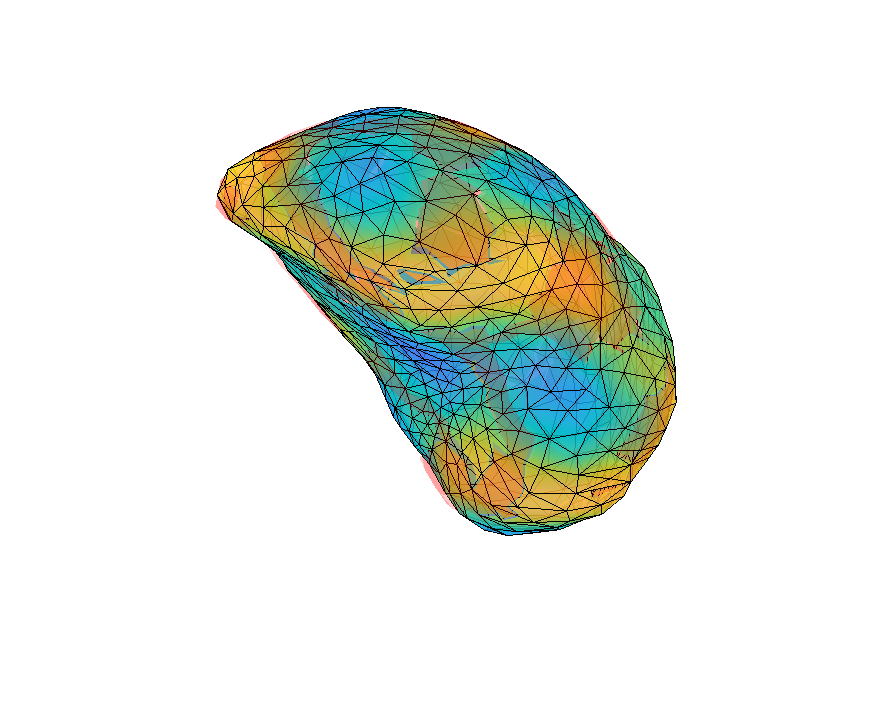}
    &\\
    
    &$t=1/3$ &  $t=7/15$  &  $t=4/5$ & $t=1 $& 
    \end{tabular}   
    \caption{Registration of two amygdala surfaces (data originating from the BIOCARD dataset \cite{Miller2015b}) with missing subregions computed with the LDDMM-L$^2$, LDDMM-FR, and standard LDDMM approach. Colors on the second and third rows correspond to the weight values at each location. The Chamfer distance to the ground truth is 0.7476 (standard LDDMM), 0.6252 (LDDMM-L$^2$), and $0.5618$ (LDDMM-FR.)} 
    \label{fig:amygdala}
\end{center}
\end{figure}

In Fig. \ref{fig:circle to ellipse}, we show the result for a simple simulated example in which the unit circle with uniform weight of $1$ is matched to an ellipse with piecewise constant weights (displayed by the different colors in the upper left image). On the first row is shown the obtained optimal deformations $\varphi_1^v$ via deformation grids as well as several intermediate time steps of the shape's evolution for both the LDDMM-L$^2$ and LDDMM-FR models (note that for the static $L^2$ model, the weights being displayed correspond to the linear interpolation $\left((1-t) + t \alpha_i\right)r_i$). Both models lead to a close geometric match with similar diffeomorphic deformations as well as a relatively good agreement with the weights of the target shape. We observe nevertheless a sharper concentration of weights around the ground truth values in the case of LDDMM-FR as evidenced by the histogram plots although this comes at the expense of small oscillations around those values.   

In Fig. \ref{fig:amygdala}, we consider the registration of two surfaces of amygdala each with uniform weights equal to $1$, the target surface being obtained by artificially removing portions of the ground truth target. In this experiment, we took $\lambda =10$ in all models and $\gamma = 0.1$ for LDDMM-L$^2$ and LDDMM-FR. Such missing parts typically induce, to different degrees, loss of precision in the registration when standard pure deformation models are used due to the mismatches between the mass of the source and target at the locations of those missing regions. In contrast, the joint estimation of a weight function in the LDDMM-L$^2$ and LDDMM-FR models, as evidenced by the plots in the last two rows of Fig. \ref{fig:amygdala}, allows to automatically identify the corresponding missing portions of the source shape and set their weights closer to $0$. This leads in turn to a more accurate match. To evaluate it more precisely, we also computed the Chamfer distance between the deformed source surface and the original (complete) ground truth, showing that the LDDMM-FR model achieves the lowest value in this example.

\begin{figure}[ht]
\begin{center}
    \begin{tabular}{cc}
        \includegraphics[trim= 25mm 20mm 15mm 10mm ,clip,width=3.5cm]{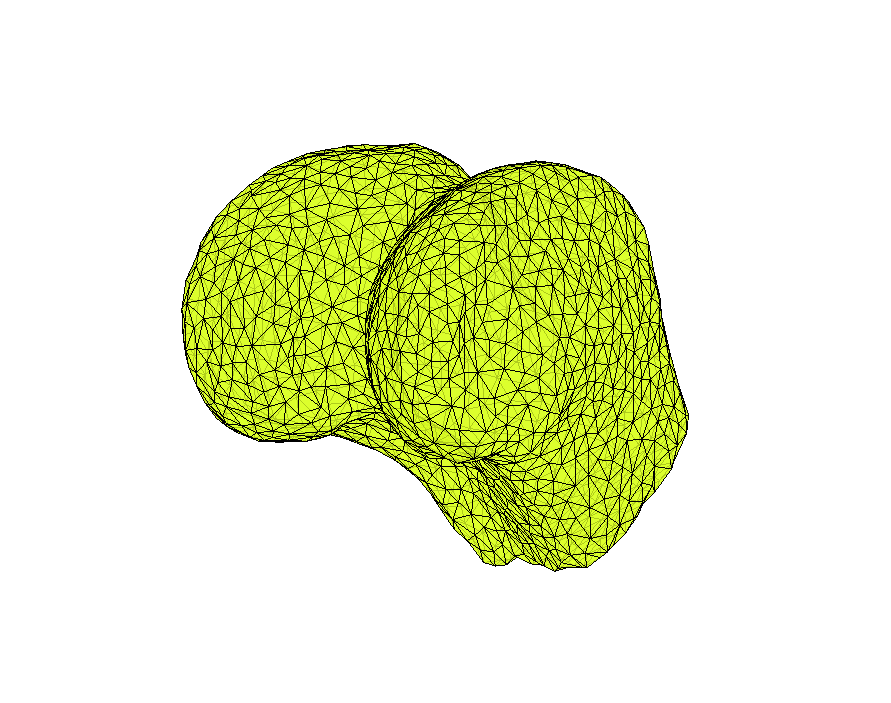}
         &\includegraphics[trim= 25mm 20mm 15mm 10mm ,clip,width=3.5cm]{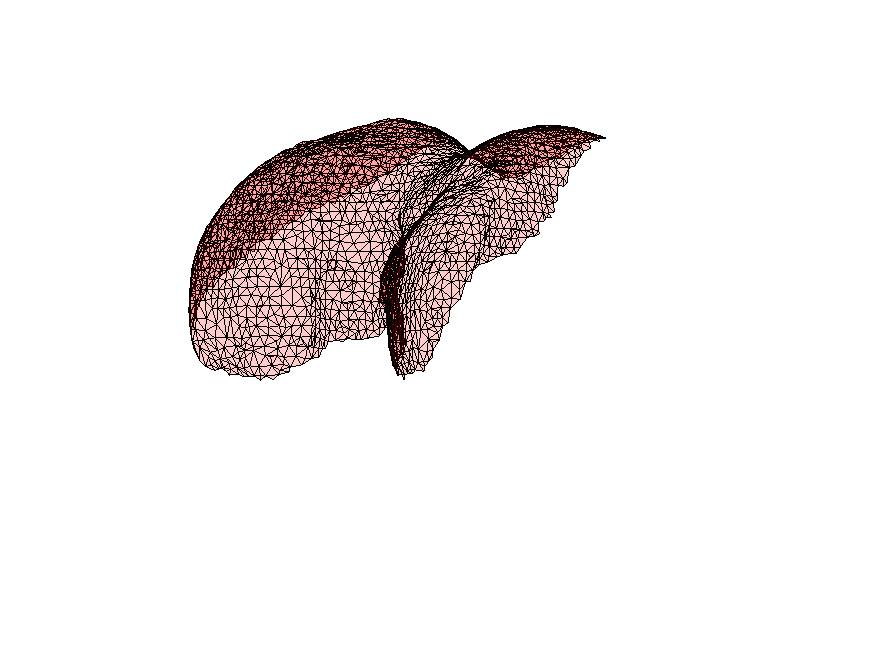} \\
         Source & Target
    \end{tabular}
    
     \begin{tabular}{cccccc}
     \rotatebox{90}{\phantom{LD}\small{LDDMM}}
    &\includegraphics[trim= 25mm 15mm 15mm 10mm ,clip, width=2.4cm]{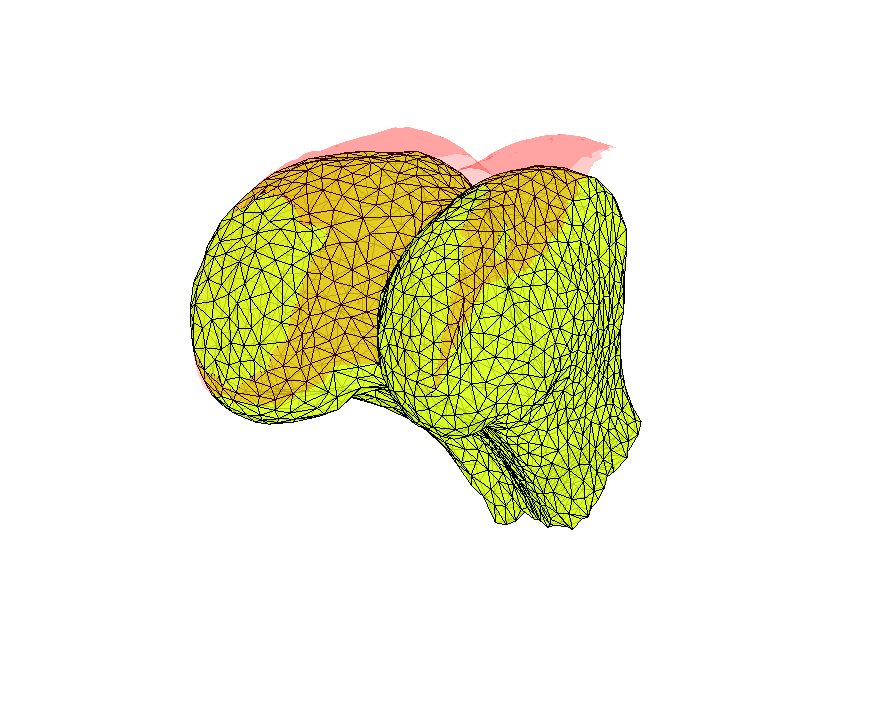}
    &\includegraphics[trim= 25mm 15mm 15mm 10mm ,clip, width=2.4cm]{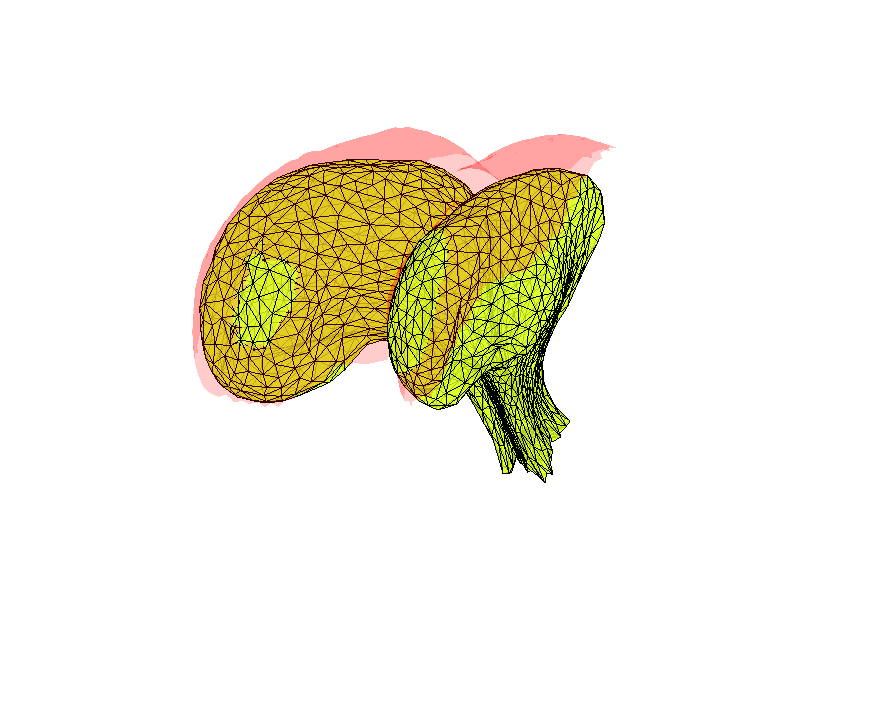}
    &\includegraphics[trim= 25mm 15mm 15mm 10mm ,clip, width=2.4cm]{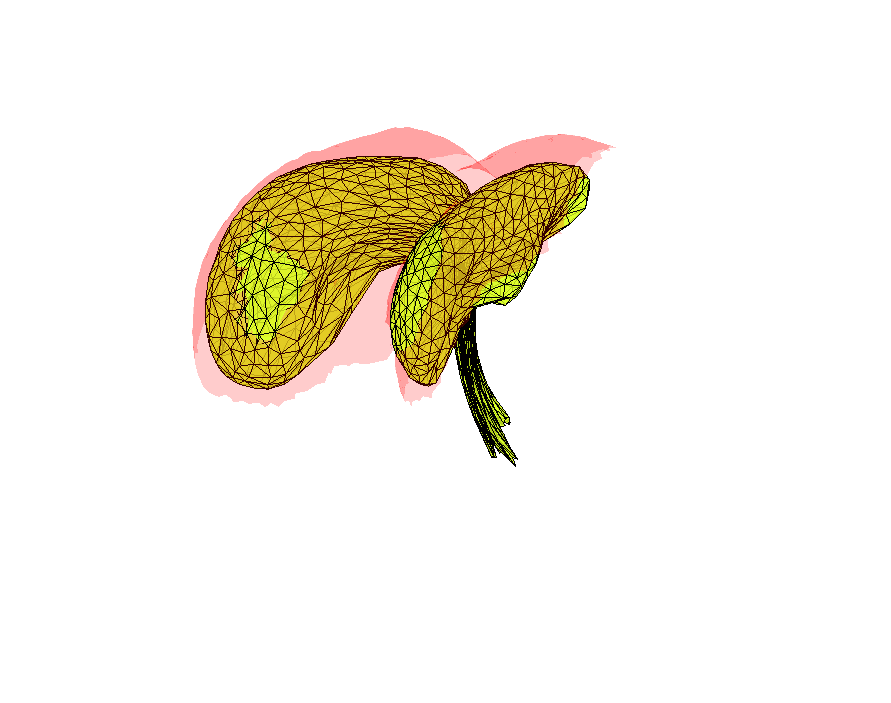}
    &\includegraphics[trim= 25mm 15mm 15mm 10mm ,clip, width=2.4cm]{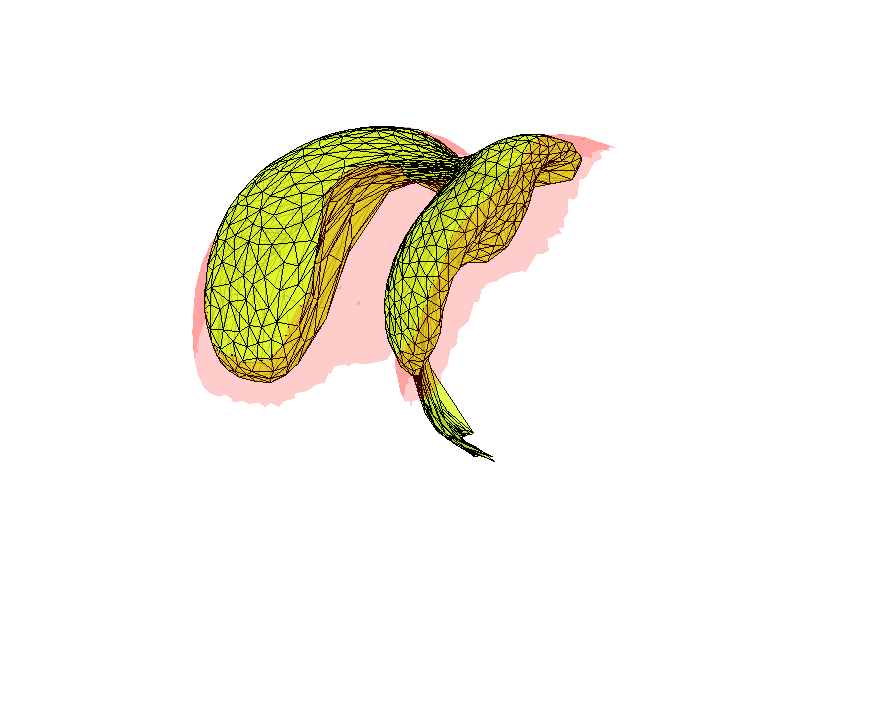}
    &\\
     
     \rotatebox{90}{\small{LDDMM-L$^2$}}
    &\includegraphics[trim= 25mm 15mm 15mm 10mm ,clip, width=2.4cm]{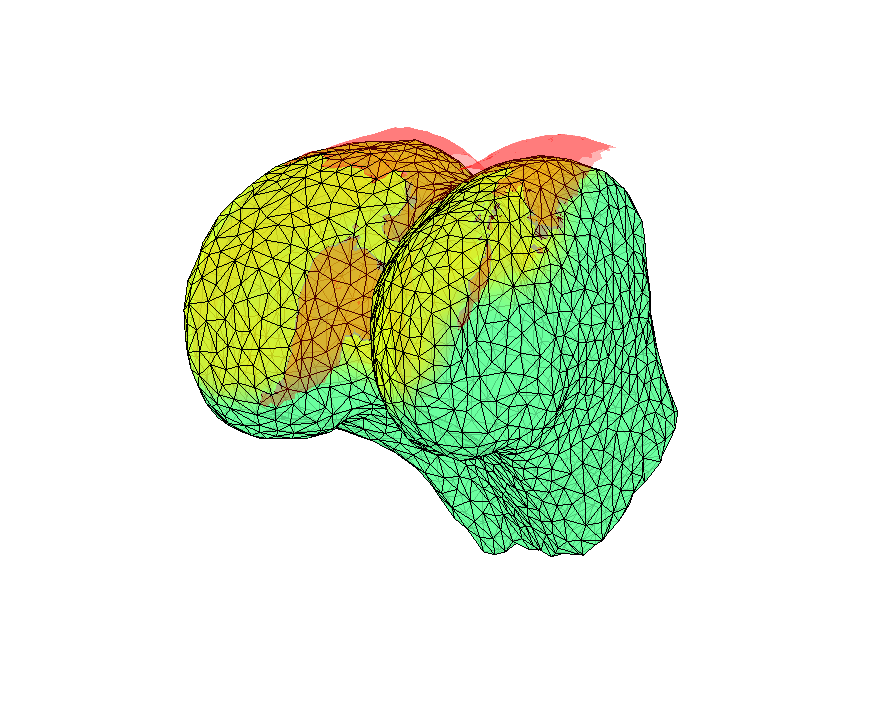}
    &\includegraphics[trim= 25mm 15mm 15mm 10mm ,clip, width=2.4cm]{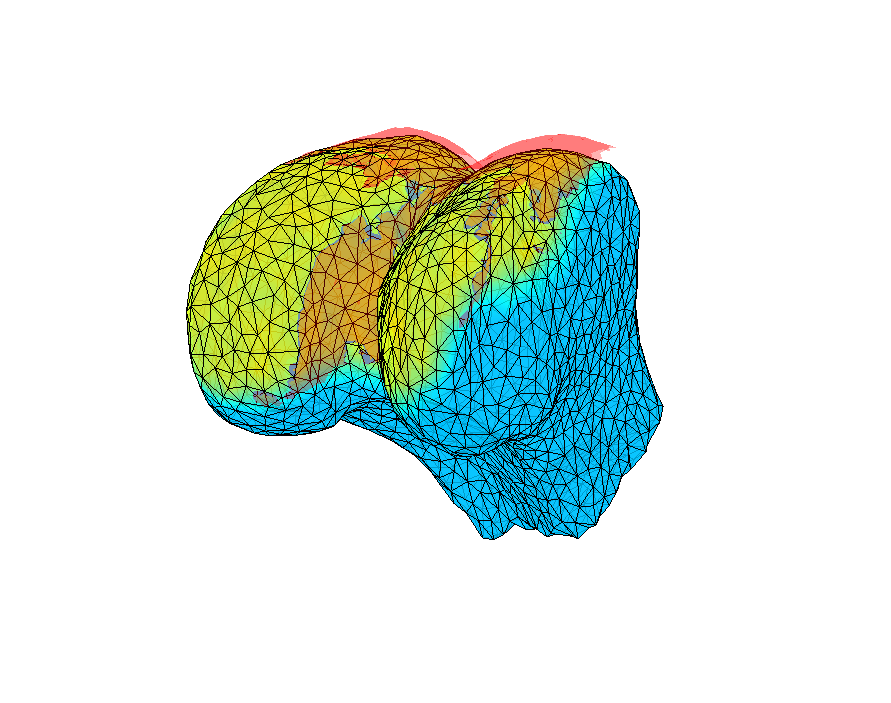}
    &\includegraphics[trim= 25mm 15mm 15mm 10mm ,clip, width=2.4cm]{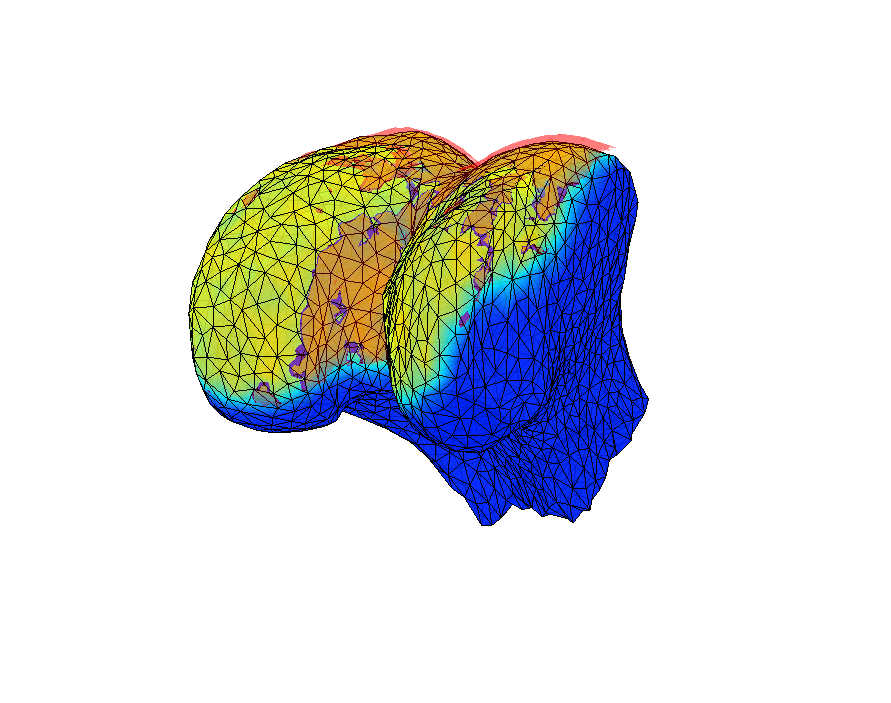}
    &\includegraphics[trim= 25mm 15mm 15mm 10mm ,clip, width=2.4cm]{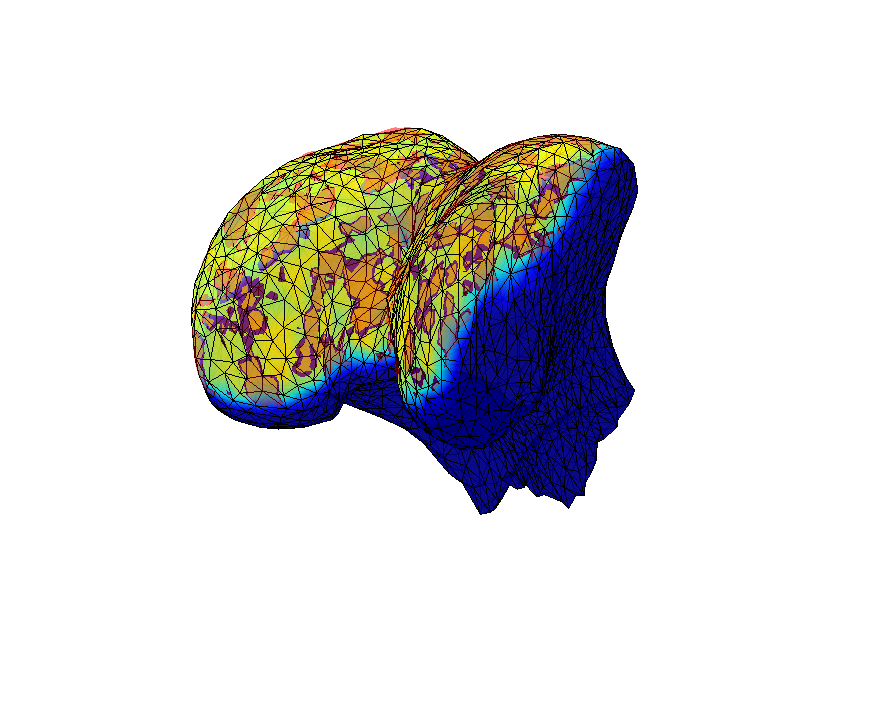}
    &\multirow{ 2}{*}[1.8cm]{\includegraphics[height=4cm]{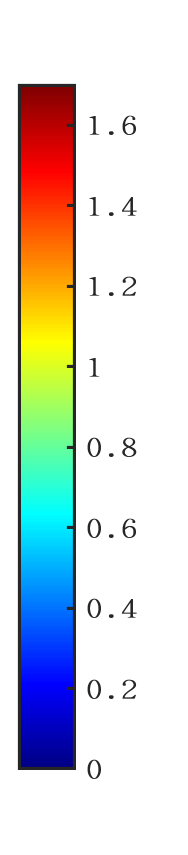}}\\
    
    \rotatebox{90}{\small{LDDMM-FR}}
    &\includegraphics[trim= 25mm 15mm 15mm 10mm ,clip, width=2.4cm]{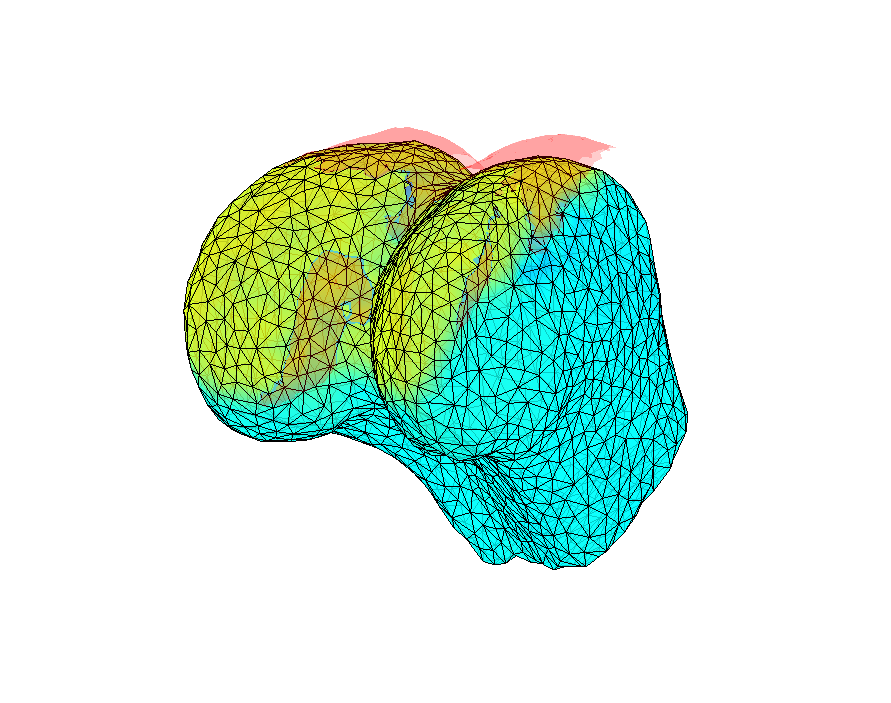}
    &\includegraphics[trim= 25mm 15mm 15mm 10mm ,clip, width=2.4cm]{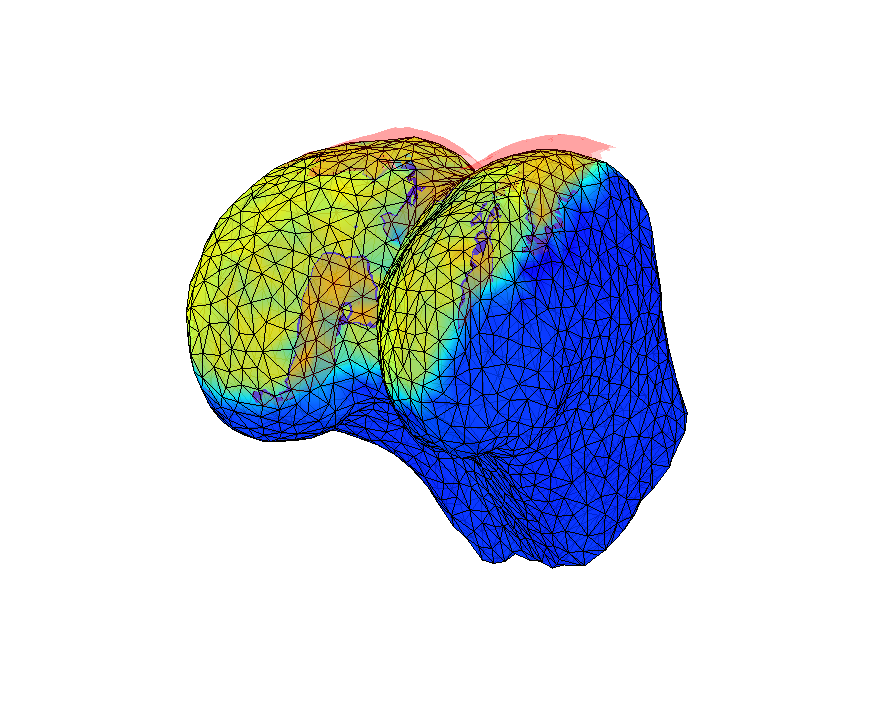}
    &\includegraphics[trim= 25mm 15mm 15mm 10mm ,clip, width=2.4cm]{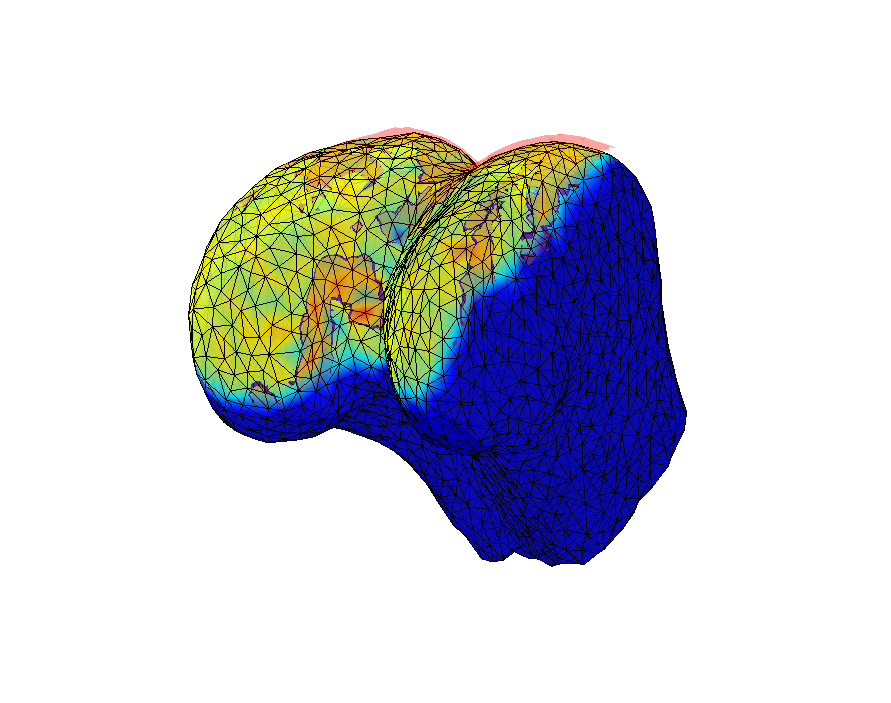}
    &\includegraphics[trim= 25mm 15mm 15mm 10mm ,clip, width=2.4cm]{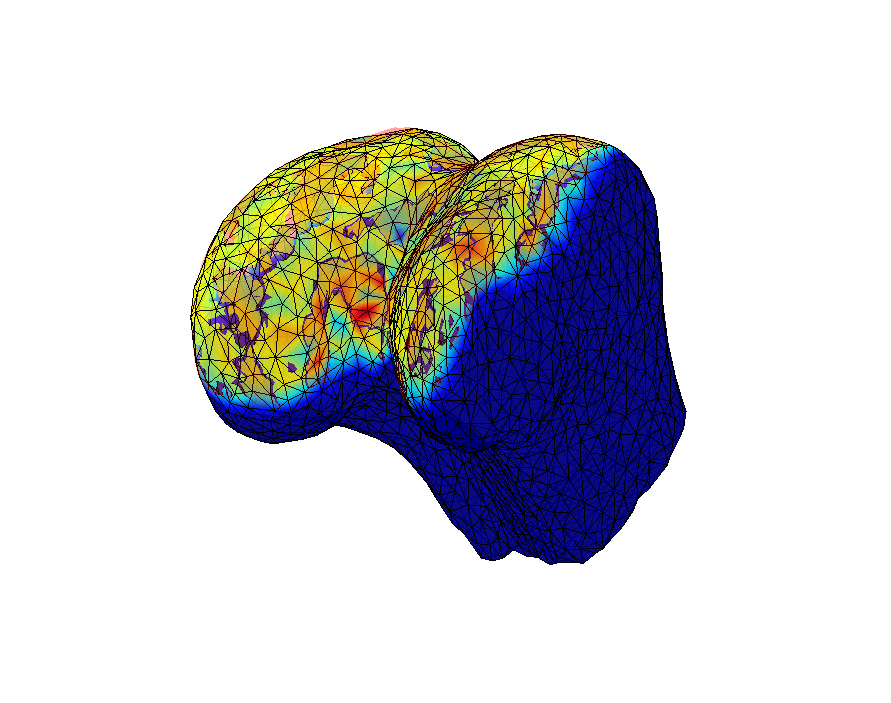}
    &\\
    
    &$t=1/3$ &  $t=7/15$  &  $t=4/5$ & $t=1 $& 
    \end{tabular}   
    \caption{Registration of two femur surfaces (data courtesy of W. Zbijewski from the Biomedical Imaging department at the Johns Hopkins University) with missing subregions obtained
with the standard LDDMM (corresponding to the limit case $\gamma \rightarrow +\infty$ of our models), the LDDMM-L$^2$ and the LDDMM-FR approach. The colors on the
second and third rows correspond to the values of the estimated weight density at each point.} 
    \label{fig:knee}
\end{center}
\end{figure}

\begin{figure}[ht]
    \centering
    \begin{tabular}{ccccc}
        \rotatebox{90}{\phantom{LD}\small{LDDMM}}
         &\includegraphics[width=2.4cm]{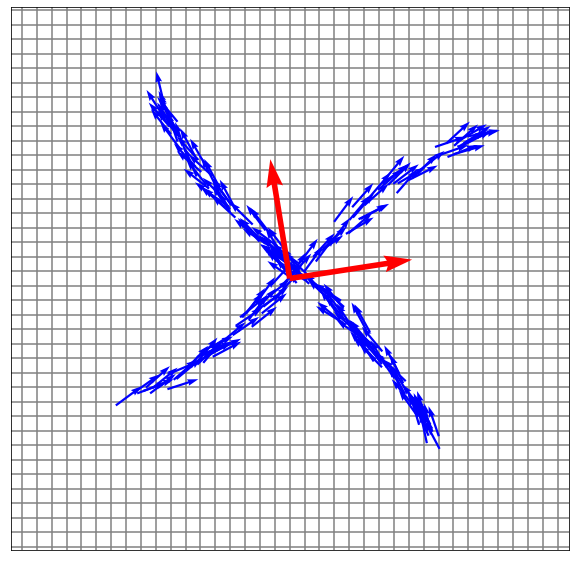}
         &\includegraphics[width=2.4cm]{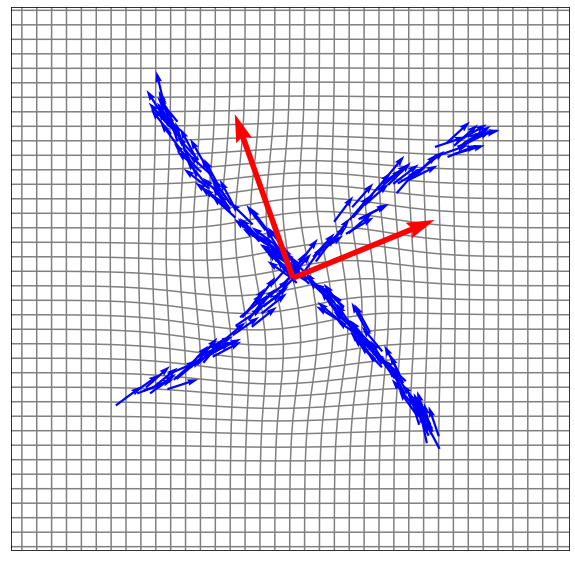}
         &\includegraphics[width=2.4cm]{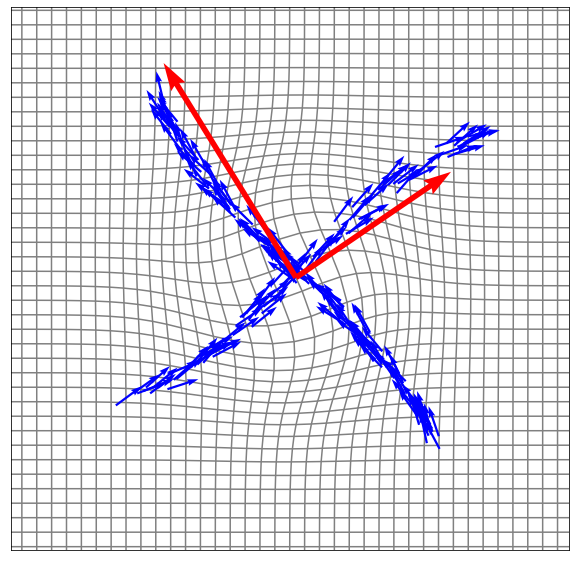}
         &\includegraphics[width=2.4cm]{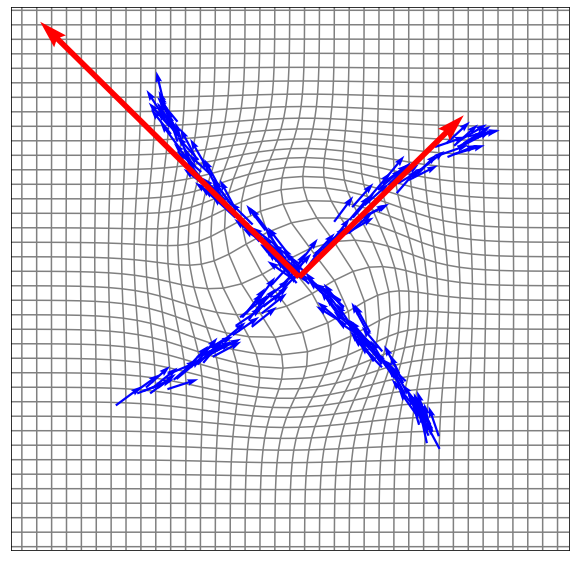}\\
         
         \rotatebox{90}{\phantom{L}\small{LDDMM-FR}}
         &\includegraphics[width=2.4cm]{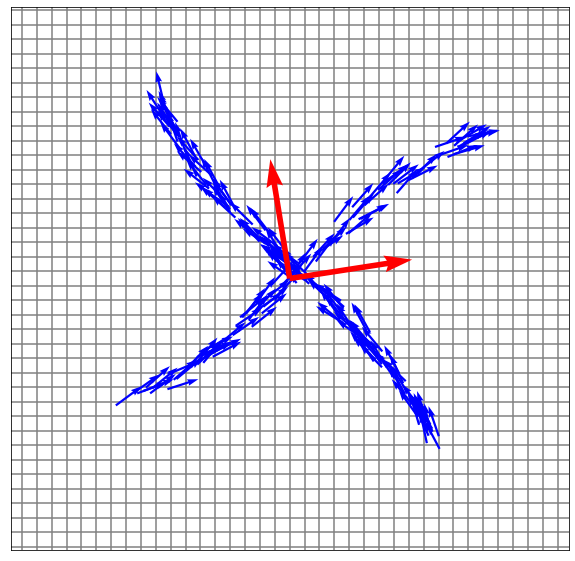}
         &\includegraphics[width=2.4cm]{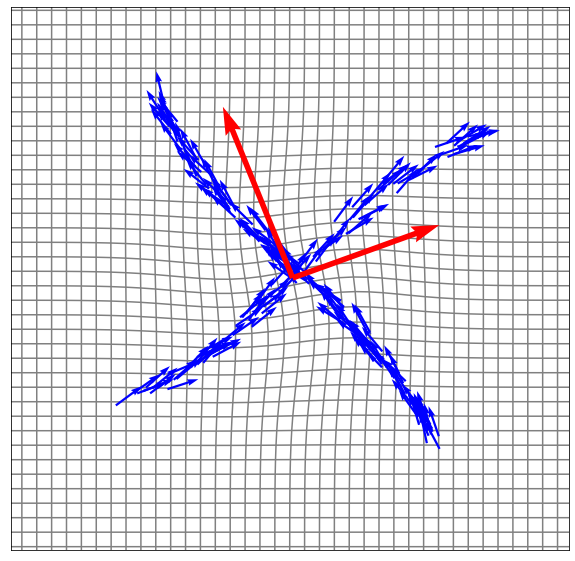}
         &\includegraphics[width=2.4cm]{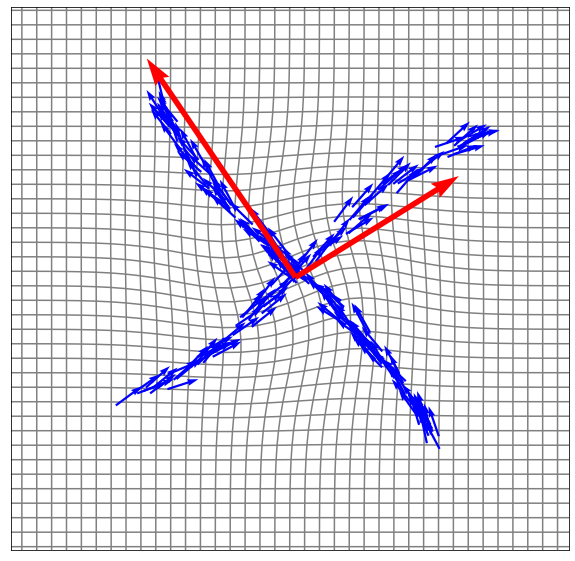}
         &\includegraphics[width=2.4cm]{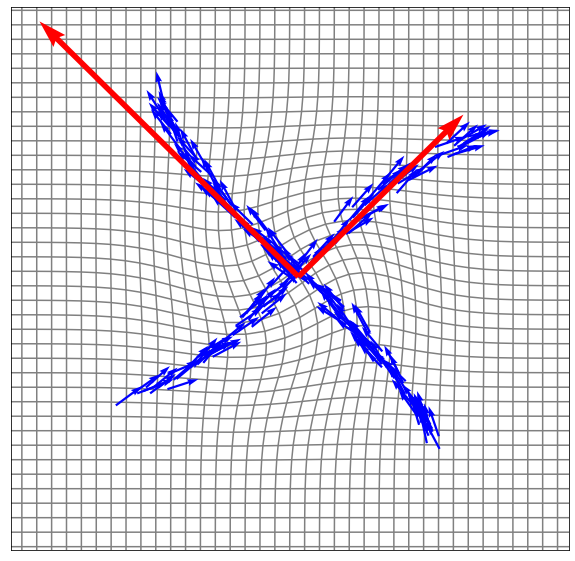}\\
         
         &$t=1/3$ &  $t=7/15$  &  $t=4/5$ & $t=1 $
    \end{tabular}
    
    \caption{Registration of a Dirac pair (red) onto crossing bundles of Diracs with distinct densities (120 and 80 Diracs respectively) obtained with the standard LDDMM and the LDDMM-FR approaches.}
    \label{fig:distribution_diracs}
\end{figure}

This is even more clearly exemplified by the result of Fig. \ref{fig:knee} that shows the registration of femur surfaces in which the target shape only corresponds to a relatively small part of the source. In such a case, standard diffeomorphic registration of the two varifolds (computed with the approach of \cite{hsieh2021metrics} or equivalently taking $\gamma$ very large in our models) shown on the second row leads to severe shrinking of the bottom section of the source femur as an attempt to geometrically eliminate this extra mass as well as important residual mismatch of the rest of the bone. This is mostly alleviated by the two models LDDMM-L$^2$ and LDDMM-FR which instead successfully erase the unmatching part by setting weights to $0$ at those locations while resulting in much more natural deformations.  

Lastly, with Fig. \ref{fig:distribution_diracs}, we show an example of registration involving non-manifold synthetic data in which we compute the registration of a pair of two $1$-varifold Diracs in $\R^2$ (shown as red arrows in the figure with length corresponding to their respective weight) onto the crossing distribution of Diracs shown in blue. For both the standard LDDMM and LDDMM-FR models (here with $\lambda =10$), we obtain a transformed Dirac pair in which each Dirac aligns with one of the two bundles while their respective weights (i.e. length of the arrows in the plot) match the relative densities of the corresponding bundles. In the case of the pure LDDMM model however, this is done via different local expansion of the diffeomorphic mapping along the two bundle directions. In contrast, with the LDDMM-FR model (for which we chose $\gamma =0.01$), the deformation mainly restricts to aligning the Dirac's directions whereas the change of mass is the result of the weight change function. Although only a very simple case, this example suggests the possible interest of these types of models to register or process, e.g., white matter fiber bundles with varying fiber densities as was already pointed out in the authors' previous work \cite{hsieh2021diffeomorphic}.

\section{Conclusion}
In this paper, we proposed extensions of the LDDMM registration model for geometric objects represented as varifolds in which the diffeomorphic deformation is combined with a transformation process of the source varifold weight function. We considered two different classes of cost functionals for these joint transformations, in particular the LDDMM-FR energy in which the weight change is penalized by the Fisher-Rao metric. We showed that the latter is associated to a metamorphosis model on the space of varifolds and induces a well-defined Riemannian metric. We further formulated and studied the corresponding inexact registration problems in which the terminal matching constraint is relaxed through the use of a kernel fidelity metric. Lastly, we derived numerical implementations of those approaches and showcased their potential interest, most notably when it comes to partial geometric data. 

Among the persistent shortcomings of the framework developed in this work, we should mention that it remains, to a certain degree, asymmetric in that it allows to modify the weight on the support of the source varifold (and even set it to zero) but cannot a priori "generate" new mass outside of the existing support. In applications to partial data registration, this implies that the model is well suited to deal with missing parts on the target shape but not on the source. Although it is a relatively common in shape analysis that the source shape is taken as a template and thus assumed to be complete, extending this approach to the situation of partial observations for both the source and target is an important and challenging open problem to address in future work. 

A second avenue for future research would be to replace the $L^2$ or Fisher-Rao penalties by different regularization metrics for the weight change function in either a static or metamorphosis setting, with the purpose of imposing spatially smoother weights. In the special case of rectifiable varifolds, one could for instance introduce higher-order Sobolev or total variation norms of the weight change function, by analogy with what has been considered in the context of functional shapes in \cite{Charon2018} or elastic shape analysis in \cite{sukurdeep2021new}.      

\section*{Acknowledgments}
This work was supported by the National Science Foundation (NSF) under the grant DMS-1945224. The authors would also like to thank Dr. Wojtek Zbijewski and the I-STAR lab for providing the femur data used in Fig. \ref{fig:knee} as well as Dr. Sylvain Arguill\`{e}re for some helpful discussions.

\begin{appendices}
\section{Proof of Theorem \ref{thm:distance_LDDMM_FR}}
We first prove symmetry by showing that for any $\mu, \mu' \in \Theta(\mu_0)$, if there exists an admissible path $\mu_t$ from $\mu$ to $\mu'$, then there exists a reversed admissible path $\mu'_t$ from $\mu'$ to $\mu$ which consumes the same energy as $\mu_t$. By definition,  there exists $(v,\eta) \in L^2([0,1],V \times L^2(\vert\mu \vert))$ such that $\mu_t \doteq ((\tilde{\alpha}_t^{\eta})^2 \circ (\varphi^t_v)^{-1})(\varphi_t^v)_{\#}  \mu$. If we define
\begin{align*}
(\bar{v}_t,\bar{\eta}_t) \doteq -(v_{1-t}, \left(\frac{\eta_{1-t}}{\tilde{\alpha}^{\eta}_1}\right) \circ (\varphi_1^v)^{-1}),
\end{align*}
then it follows from the results of \cite{Younes2019} (Chap. 7) that $\varphi_t^{\bar{v}} \circ \varphi_{1}^v = \varphi_{1-t}^v$. Furthermore 
\begin{align*}
&(\tilde{\alpha}_t^{\bar{\eta}} \circ (\varphi_t^{\bar{v}})^{-1})(\tilde{\alpha}_1^{\eta} \circ (\varphi_{1-t}^v)^{-1}) \\
&= \left(\left(1- \frac{1}{2\tilde{\alpha}^{\eta}_1(\cdot)}\int_0^t \eta_{1-s}(\cdot)ds  \right) \circ (\varphi^v_1)^{-1} \circ (\varphi^{\bar{v}}_t)^{-1} \right) \left( \tilde{\alpha}_1^{\eta} \circ (\varphi_{1-t}^v)^{-1} \right)\\
&=\left(\tilde{\alpha}^{\eta}_1 - \frac{1}{2}\int_{1-t}^1 \eta_s(\cdot) ds \right)\circ (\varphi_{1-t}^v)^{-1} \\
&=\tilde{\alpha}_{1-t}^{\eta} \circ (\varphi^v_{1-t})^{-1}.
\end{align*}
Moreover,
\begin{align*}
\mu'_t &\doteq ((\tilde{\alpha}_t^{\bar{\eta}})^2 \circ (\varphi_t^{\bar{v}})^{-1})(\varphi_t^{\bar{v}})_{\#} \mu' \\
&= ((\tilde{\alpha}_t^{\bar{\eta}})^2 \circ (\varphi_t^{\bar{v}})^{-1}) (\varphi_t^{\bar{v}})_{\#} [((\alpha^{\eta}_1)^2 \circ (\varphi^v_1)^{-1}) (\varphi_1^v)_{\#} \mu]  \\
&= (\tilde{\alpha}_{1-t}^{\eta})^2 \circ (\varphi^v_{1-t})^{-1} (\varphi_{1-t}^v)_{\#} \mu \\
&= \mu_{1-t}.
\end{align*}
From the equalities above, we can see that
\begin{align*}
\int_{\R^n} \bar{\eta}_t^2 \circ (\varphi_t^{\bar{v}})^{-1} d\vert (\varphi_t^{\bar{v}})_{\#} \mu'\vert
&=\int_{\R^n} \left(\frac{\eta_{1-t}}{\tilde{\alpha}^{\eta}_1}\right)^2 \circ (\varphi_1^v)^{-1} \circ (\varphi_t^{\bar{v}})^{-1} \frac{1}{(\tilde{\alpha}^{\bar{\eta}}_t)^2 \circ (\varphi_t^{\bar{v}})^{-1}} d\vert \mu'_t\vert\\
&= \int_{\R^n} \frac{\eta_{1-t}^2 \circ (\varphi_{1-t}^v)^{-1}}{(\tilde{\alpha}^{\eta}_1)^2 \circ (\varphi_{1-t}^v)^{-1}} 
\frac{(\tilde{\alpha}^{\eta}_1)^2 \circ (\varphi_{1-t}^v)^{-1}}{(\tilde{\alpha}_{1-t}^{\eta})^2 \circ (\varphi^v_{1-t})^{-1}} d\vert\mu_{1-t}\vert \\
& = \int_{\R^n} \eta_{1-t}^2 \circ (\varphi_{1-t}^v)^{-1} d \vert (\varphi_{1-t}^v)_{\#}\mu\vert.
\end{align*}
Therefore, the cost for the reverse trajectory is the same as the original one:
\begin{align*}
C_{\mu'}(\bar{v},\bar{\eta}) &= \frac{1}{2} \int_0^1 \|v_{1-t}\|_V^2 dt + \frac{\gamma}{2} \int_0^1\int_{\R^n}\eta_{1-t}^2 \circ (\varphi_{1-t}^v)^{-1} d \vert (\varphi_{1-t}^v)_{\#}\mu\vert dt \\
&= \frac{1}{2} \int_0^1 \|v_t\|_V^2 dt + \frac{\gamma}{2} \int_0^1\int_{\R^n}\eta_t^2 \circ (\varphi_t^v)^{-1} d \vert (\varphi_t^v)_{\#}\mu\vert dt \\
&=C_{\mu}(v,\eta) < \infty.
\end{align*}
The above also implies that $(\bar{v},\bar{\eta}) \in L^2([0,1],V \times L^2(\vert \mu'\vert)$. Indeed, as $v \in L^2([0,1],V)$ and $V \hookrightarrow C^1_0(\R^n,\R^n)$, one has from Theorem 7.10 in \cite{Younes2019}, that $\sup_{t\in[0,1]} \{ \|\varphi_t^v - \text{id}\|_{1,\infty},\|(\varphi_t^v)^{-1} - \text{id}\|_{1,\infty}\} <+\infty$ which implies in particular the existence of $M>0$ such that $1/M\leq J_T \varphi^{\bar{v}}_t(x)\leq M$ for all $t \in [0,1]$, $x \in \R^n$ and $T \in \origrass$. Then 
\begin{align*}
\int_0^1 \int_{\R^n} \bar{\eta}^2_t(x) d\vert \mu'\vert(x) dt &\leq M \int_0^1 \int_{\R^n \times \origrass} \bar{\eta}^2_t(x) J_T \varphi^{\bar{v}}_t(x) d\mu'(x,T) dt \\
&\leq M C_{\mu'}(\bar{v},\bar{\eta}) < \infty.
\end{align*}
where we used \eqref{eq:ener_LDDMM_FR2} for the first inequality. This shows that $\mu_t'$ is an admissible path from $\mu'$ to $\mu$ and that $d_{\mathcal{D-FR}}$ is symmetric.

%This leads to the symmetry of $d_{\mathcal{D-FR}}$.

Second, we show that if $d_{\mathcal{D-FR}}(\mu,\mu')=0$, then $\mu=\mu'$.
Now, from the definition of $d_{\mathcal{D-FR}}$, there exists a minimizing sequence $(v^j,\eta^j)$ such that $\mu' = (\varphi^{v^j}_1)_{\#} (\tilde{\alpha}_1^{\eta^j}\mu)$ and $\lim_{j \rightarrow \infty} C_{\mu}(v^j,\eta^j)=0$. This implies in particular that $\int_0^1 \|v^j\|_V^2 dt$ is uniformly bounded over $j$ and from the same argument as above we can find $M>0$ such that $1/M\leq J_T \varphi^{v^j}_t(x)\leq M$ for all $j\in \N$, $t \in [0,1]$, $x \in \R^n$ and $T \in \origrass$. We deduce that
 \begin{align*}
 \int_{\R^n} (\tilde{\alpha}^{\eta^j}_1(x)-1)^2 d \vert\mu\vert(x) 
 &= \int_{\R^n} \left(\int_0^1 \eta_t^j(x) dt \right)^2 d\vert\mu\vert(x) \\
 &\leq \int_{\R^n} \int_0^1 (\eta_t^j(x))^2 dt \ d \vert\mu\vert(x) \\
 &\leq M \int_0^1 \int_{\R^n \times \tilde{R}^n_d} (\eta_t^j(x))^2 J_T \varphi_t^{v^j}(x) d \mu(x,T) \rightarrow 0, 
 \end{align*}
as $j \rightarrow \infty$, the first bound following from Cauchy-Schwarz inequality and the last limit from the fact that the integral is upper bounded by $C_{\mu}(v^j,\eta^j) \rightarrow 0$. This also implies that $\int_{\R^n} (\tilde{\alpha}_1^{\eta^j}(x))^2 d\vert\mu\vert(x) \rightarrow \vert\mu\vert(\R^n)<\infty$ and $(\tilde{\alpha}_1^{\eta^j}(x))^2 \rightarrow 1$ $\vert\mu\vert$-$a.e.$ up to a subsequence. From Theorem 4.6.2 in \cite{durrett2019probability}, $(\alpha^{\eta^j}_1)^2$ converges to $1$ in $L^1(\vert\mu\vert)$. For any bounded continuous function $\omega$ and $j\in \N$
\begin{align*}
(\mu-\mu'\vert\omega) &= \underbrace{\int_{\mathbb{R}^n \times \tilde{G}^n_d} \omega(x,T) d\mu(x,T) - \int_{\mathbb{R}^n \times \tilde{G}^n_d} \omega(\varphi_1^{v^j}(x), d_x \varphi_1^{v^j} \cdot T)J_T\varphi_1^{v^j}(x) d\mu(x,T)}_{\MakeUppercase{\romannumeral 1}} \\
&+ \underbrace{ \int_{\mathbb{R}^n \times \tilde{G}^n_d} \omega(\varphi_1^{v^j}(x), d_x \varphi_1^{v^j} \cdot T)(1-(\tilde{\alpha}^{\eta^j}_1(x))^2) J_T\varphi_1^{v^j}(x) d\mu(x,T)}_{\MakeUppercase{\romannumeral 2}}.
\end{align*}
It is clear that $\MakeUppercase{\romannumeral 1} \rightarrow 0$ from bounded convergence theorem. Also, from the fact that $(\tilde{\alpha}_1^{\eta^j})^2$ converges to $1$ in $L^2(\vert\mu\vert)$ and the uniform control of Jacobian above, we obtain that
\begin{align*}
\vert\MakeUppercase{\romannumeral 2}\vert
\leq  C \int_{\mathbb{R}^n} \vert 1-(\tilde{\alpha}^{\eta^j}(x))^2\vert  \, d\vert\mu\vert(x) \rightarrow 0.
\end{align*}
This shows that $\mu = \mu'$. 

Next, we prove the triangular inequality. Take $\mu$, $\mu'$ and $\mu''$ in $\Theta(\mu_0)$ such that $0<d_{\Theta}(\mu,\mu'), \ d_{\Theta}(\mu',\mu'')<\infty$. Let $(v,\eta) \in L^2([0,1],V \times L^2(\mathbb{R}^n,\vert\mu \vert)$ and $(v',\eta') \in L^2([0,1],V \times L^2(\vert\mu' \vert))$ be such that
\begin{align*}
\mu_1 = \mu' \textrm{ and } \mu'_1 = \mu'',
\end{align*}
where
\begin{align*}
\mu_t \doteq (\varphi^v_t)_{\#} (\tilde{\alpha}^{\eta}_t \mu) \textrm{ and }
\mu'_t \doteq (\varphi^{v'}_t)_{\#} (\tilde{\alpha}^{\eta'}_t \mu').
\end{align*}
Now, let $a,b>1$ such that $1/a+1/b=1$ and denote $u_t \doteq \min{\{at,1\}}$ and $u'_t \doteq \max{\{b(t-1/a),0\}}$. Defining
\begin{align*}
(\bar{v}_t,\bar{\eta}_t) \doteq  a (v_{u_t},\eta_{u_t}) \mathbf{1}_{[0,1/a)}(t)
 + b(v'_{u'_t},\tilde{\alpha}^{\eta}_1\eta'_{u'_t} \circ \varphi_1^v) \mathbf{1}_{[1/a,1]}(t),
\end{align*}
then we can easily check that
\begin{align*}
&\varphi^{\bar{v}}_t = \varphi^{v'}_{u'_t} \circ \varphi^v_{u_t}, 
\end{align*}
and
\begin{align*}
\tilde{\alpha}^{\bar{\eta}}_t &= 1+ \frac{a}{2} \int_0^{\min \{t,1/a\}} \eta_{u_s} ds + \frac{\tilde{\alpha}^{\eta}_1 b}{2} \int_{1/a}^{\max\{t,1/a\}} \eta'_{u_s'} \circ \varphi^v_1 ds \\
&= 1 + \frac{1}{2}\int_0^{u_t} \eta_s ds + \frac{\tilde{\alpha}^{\eta}_1}{2} \int_0^{u_t'} \eta'_s \circ \varphi^v_1 ds \\
&=\tilde{\alpha}^{\eta}_{u_t} + \tilde{\alpha}^{\eta}_1 \int_0^{u_t'} \eta_s' \circ \varphi^v_1 ds \\
&= \left(1+ \int_0^{u_t'} \eta_s' \circ \varphi^v_1 ds  \right) \tilde{\alpha}^{\eta}_{u_t} = (\tilde{\alpha}^{\eta'}_{u'_t} \circ \varphi^v_1) \tilde{\alpha}^{\eta}_{u_t}.
\end{align*}
From the above, we can further obtain
\begin{align*}
&C_{\mu}(\bar{v},\bar{\eta}) = \frac{a}{2} \left( \int_0^{1/a} \left( \|v_{u_t}\|_V^2 + \gamma \int_{\R^n} \eta^2_{u_t} \circ (\varphi_{u_t}^v)^{-1}(x) d\vert(\varphi_{u_t}^v)_{\#}\mu\vert(x) \right) \frac{du_t}{dt} dt \right) \\
&+ \frac{b}{2} \left(\int_{1/a}^1 \left(\|v'_{u_t'}\|^2_V  + \gamma \int_{\R^n} \eta^2_{u'_t} \circ (\varphi^{v'}_{u_t'})^{-1}(x) (\alpha^{\eta}_1)^2 \circ (\varphi_t^{\bar{v}})^{-1} d\vert(\varphi_{u_t'}^{\bar{v}})_{\#}\mu\vert(x)\right) \frac{du'_t}{dt} dt \right) \\
&=aC_{\mu}(v,\eta) + \frac{b}{2} \left(\int_{1/a}^1 \left(\|v'_{u_t'}\|^2_V  + \gamma \int_{\R^n} \eta^2_{u'_t} \circ (\varphi^{v'}_{u_t'})^{-1}(x)  d\vert(\varphi_{u_t'}^{v'})_{\#}\mu'\vert(x)\right) \frac{du'_t}{dt} dt \right) \\
& = aC_{\mu}(v,\eta) + bC_{\mu'}(v',\eta')< \infty,
\end{align*}
where the second equality came from the fact that 
\begin{align*}
((\tilde{\alpha}^{\eta}_1)^2 \circ (\varphi^{\bar{v}}_t)^{-1}) (\varphi^{\bar{v}}_t)_{\#} \mu &= (\varphi_t^{\bar{v}})_{\#} (\tilde{\alpha}_1^{\eta})^2 \mu = (\varphi^{v'}_{u_t'})_{\#} (\varphi^{v}_{1})_{\#} ((\alpha^{\eta}_1)^2 \mu) = (\varphi^{v'}_{u_t'})_{\#} \mu'.
\end{align*}
Again, similar to the argument above, this implies that $(\bar{v},\bar{\eta}) \in L^2([0,1],V\times L^2(\vert\mu\vert))$. Then, taking $a^* = 1+\frac{C_{\mu'}(v',\eta')^{1/2}}{C_{\mu}(v,\eta)^{1/2}}$, we obtain that
\begin{align*}
C_{\mu}(\bar{v},\bar{\eta})^{1/2} = C_{\mu}(v,\eta)^{1/2} + C_{\mu'}(v',\eta')^{1/2}.
\end{align*}
Now, for any $\varepsilon >0$, we can find $(v,\eta)$ and $(v',\eta')$ satisfying $C_{\mu}(v,\eta)^{1/2} \leq d_{\mathcal{D-FR}}(\mu,\mu') + \varepsilon$ and $C_{\mu}(v',\eta')^{1/2} \leq d_{\mathcal{D-FR}}(\mu',\mu'') + \varepsilon$. Then, taking $(\bar{v},\bar{\eta})$ as constructed above, we get:
\begin{equation*}
   d_{\mathcal{D-FR}}(\mu,\mu'') \leq C_{\mu}(\bar{v},\bar{\eta})^{1/2} \leq d_{\mathcal{D-FR}}(\mu,\mu') + d_{\mathcal{D-FR}}(\mu',\mu'') +2\varepsilon.
\end{equation*}
for any $\varepsilon >0$ leading to the triangle inequality.

Finally, we can easily prove that for any $\mu,\mu' \in \Theta(\mu_0)$, $d_{\mathcal{D-FR}}(\mu,\mu')<\infty$. Indeed, from the definition of $\Theta(\mu_0)$, there are admissible paths from $\mu_0$ to $\mu$ and $\mu'$ separately. We assume that $\mu$ and $\mu'$ are distinct from $\mu_0$ since this case is trivial. From the arguments in the proof of symmetry, there exists an admissible path from $\mu$ to $\mu_0$. Moreover, we can obtain an admissible path from $\mu$ to $\mu'$ by concatenating the admissible paths from $\mu$ to $\mu_0$ and from $\mu_0$ to $\mu'$ with the same argument used above for proving the triangular inequality. Since admissible path have finite energy, we have $d_{\mathcal{D-FR}}(\mu,\mu')<\infty$.

\section{Proof of Theorem \ref{thm:geodesics_LDDMM_FR}}
\label{appendix:proof_geod_FR}
\textbf{Proof of the Lemma.} First, we point out that since $\mu' = (\varphi_1^{v},(\tilde{\alpha}_1^{\eta_0})^2) \cdot \mu = (\tilde{\alpha}_1^{\eta_0}\circ (\varphi_1^{v})^{-1})^2 (\varphi_1^{v})_{\sharp}\mu$ and since we have $(\tilde{\alpha}_1^{\eta_0}\circ (\varphi_1^{v})^{-1})^2 \in L^1(\vert(\varphi_1^{v})_{\sharp}\mu\vert)$, $\mu'$ is absolutely continuous with respect to $\vert(\varphi_1^{v})_{\sharp}\mu\vert$ and the Radon-Nykodym derivative of $\vert\mu'\vert$ with respect to the transported measure $\vert(\varphi_1^{v})_{\sharp}\mu\vert$ is precisely $(\tilde{\alpha}_1^{\eta_0}\circ (\varphi_1^{v})^{-1}(x))^2$ for all $x \in \text{supp}(\vert(\varphi_1^{v})_{\sharp}\mu\vert) = \text{supp}(\vert\mu'\vert)$. Now, denoting $h_t(x) = \int_{\origrass} J_U\varphi_t^{v}(x) d\nu_x(U)$, we see that the variational problem consists in minimizing:
 \begin{align*}
 &\int_0^1 \int_{\R^n \times \origrass} \eta^2_t(x) J_U\varphi_t^v(x) \, d\mu(x,U) dt \\
 &= \int_0^1 \int_{\R^n}  \eta^2_t(x) \left(\int_{\origrass} J_U\varphi_t^v(x) \, d\nu_x(U) \right) d\vert\mu\vert(x) dt \\
 &= \int_0^1 \int_{\R^n}  \eta^2_t(x) h_t(x) d\vert\mu\vert(x) dt
\end{align*}
over $\eta \in L^2([0,1],L^2(\vert\mu\vert))$ subject to $\tilde{\alpha}_t^{\eta}(x) = 1 + \frac{1}{2} \int_0^t \eta_s(x) ds$ with $\tilde{\alpha}_t^{\eta}(x)\geq 0$ and the boundary constraint $\tilde{\alpha}_1^{\eta}(x) = \tilde{\alpha}_1^{\eta_0}(x)$ for $\vert\mu\vert$-a.e $x \in \R^n$. Since this is a linear (although infinite-dimensional) control system with a quadratic cost function, ignoring for now the non-negativity constraint on $\tilde{\alpha}_t^{\eta}$, it has an essentially unique solution which is given by:
\begin{equation*}
 \bar{\eta}_t(x) = 2 \frac{\tilde{\alpha}_1^{\eta_0}(x)-1}{h_t(x)\int_0^1 1/h_s(x)ds}
\end{equation*}
and we see that for all $t \in [0,1]$ and for $\vert\mu\vert$-a.e $x \in \R^n$:
\begin{equation*}
 \alpha_t^{\bar{\eta}}(x) = 1 + \frac{1}{2} \int_0^t \bar{\eta}_s(x) ds = 1 + (\tilde{\alpha}_1^{\eta_0}(x)-1)\frac{\int_0^t 1/h_s(x) ds}{\int_0^1 1/h_s(x)ds} \geq 0
\end{equation*}
where the last inequality follows from the fact that $h_t(x)>0$ for all $x \in \R^n$ and $t\in [0,1]$ giving $\int_0^t 1/h_s(x) ds \leq \int_0^1 1/h_s(x)ds$ as well as $\tilde{\alpha}_1^{\eta_0}(x) \geq 0$ for $\vert\mu\vert$-a.e $x \in \R^n$. Therefore $\bar{\eta}$ is indeed the unique solution of the problem with fixed deformation field $v$.
\vskip1ex
\noindent \textbf{Proof of the Theorem.} Let $(v^m,\tilde{\eta}^m)$ be a minimizing sequence in $L^2([0,1],V \times L^2(\vert\mu\vert))$. We have by definition that for all $m \in \N$, $(\varphi_1^{v^m},(\tilde{\alpha}_1^{\eta^m})^2) \cdot \mu = \mu'$. In addition, since the sequence $C_{\mu}(v^m,\eta^m)$ is bounded, we obtain in particular that $(v^m)$ is a bounded sequence in $L^2([0,1],V)$ from which we deduce that, up to extracting a subsequence, there exists $v^* \in L^2([0,1],V)$ such that $v^m$ converges to $v^*$ weakly in $L^2([0,1],V)$. By the result of Theorem 7.13 in \cite{Younes2019}, this implies that the mapping $\varphi_t^{v^m}$ converges to the diffeomorphism $\varphi_t^{v^*}$ and $(\varphi_t^{v^m})^{-1}$ converges to $(\varphi_t^{v^*})^{-1}$ in $\|\cdot\|_{1,\infty}$ on any compact subset of $\R^n$ (and thus on $\text{supp}(\vert\mu\vert))$, and the convergence is also uniform over $t \in [0,1]$. From this it follows in addition that there exist $A>1$ such that $1/A \leq J_T\varphi_t^{v^m}(x) \leq A$ and $1/A \leq J_T(\varphi_t^{v^m})^{-1}(x) \leq A$ for all $t \in [0,1]$, $T \in \origrass$, $x \in \text{supp}(\vert\mu\vert))$ and $m \in \N$, with the same bounds also holding for $J_T\varphi_t^{v^*}(x)$.

Let us write $\mu = \vert \mu\vert \otimes \nu_x$ with $\nu_x$ a probability measure on $\origrass$ for all $x \in \R^n$ the disintegration of the varifold $\mu$. For each $m \in \N$, we have $\mu'= (\varphi_1^{v^m})_{\#}((\tilde{\alpha}_1^{\eta^m})^2 \mu) = (\tilde{\alpha}_1^{\eta^m})^2 \circ (\varphi_1^{v^m})^{-1}.(\varphi_1^{v^m})_{\#} \mu$ with $\tilde{\alpha}_1^{\eta^m}(x) \geq 0$ for $\vert\mu\vert-$ a.e. $x \in \R^n$, which we can equivalently rewrite as $(\varphi_1^{v^m})^{-1}_{\#}\mu' = (\alpha_1^{\tilde{\eta}^m})^2 \mu$. This shows that the disintegration of $(\varphi_1^{v^m})^{-1}_{\#}\mu'$ must take the form $(\varphi_1^{v^m})^{-1}_{\#}\mu' = \vert(\varphi_1^{v^m})^{-1}_{\#}\mu'\vert \otimes \nu_x$. Also, for any measurable subset $B\subseteq \R^n \times \origrass$, we have: 
\begin{equation*}
 (\varphi_1^{v^m})^{-1}_{\#}\mu'(B) = \int_{\R^n \times \origrass} \mathds{1}\vert_{B}((\varphi_1^{v^m})^{-1}(x),d_x(\varphi_1^{v^m})^{-1}\cdot T) \, J_T (\varphi_1^{v^m})^{-1}(x) d\mu'(x,T)
\end{equation*}
and by uniform convergence of $(\varphi_1^{v^m})^{-1}$ to $(\varphi_1^{v^*})^{-1}$ and $d_x(\varphi_1^{v^m})^{-1}$ to $d_x(\varphi_1^{v^*})^{-1}$ on the compact $\text{supp}(\vert\mu'\vert)$, we obtain by applying Lebesgue's dominated convergence theorem:
\begin{align*}
 (\varphi_1^{v^m})^{-1}_{\#}\mu'(B) \xrightarrow[m \rightarrow \infty]{} &\int_{\R^n \times \origrass} \mathds{1}\vert_{B}((\varphi_1^{v^*})^{-1}(x),d_x(\varphi_1^{v^*})^{-1}\cdot T) \, J_T (\varphi_1^{v^*})^{-1}(x) d\mu'(x,T) \\
 &= (\varphi_1^{v^*})^{-1}_{\#}\mu'(B)
\end{align*}
On the other hand, let $g(x)\doteq \int_{\origrass} \mathds{1}\vert_{B}(x,T) d\nu_x(T)$. We obtain from the bounded convergence theorem that 
\begin{align*}
(\varphi_1^{v^m})^{-1}_{\#}\mu'(B) &= \vert(\varphi_1^{v^m})^{-1}_{\#}\mu'\vert \otimes \nu_x(B) = \int_{\R^n \times \origrass} g(x) d(\varphi_1^{v^m})^{-1}_{\#}\mu'(x,T) \\
&\xrightarrow[m \rightarrow \infty]{} \int_{\R^n \times \origrass} g((\varphi_1^{v^*})^{-1}(x)) \, J_T (\varphi_1^{v^*})^{-1}(x) d\mu'(x,T) \\
&= \int_{\R^n} g(x) d\vert(\varphi_1^{v^*})^{-1}_{\#}\mu'|(x) = |(\varphi_1^{v^*})^{-1}_{\#}\mu'\vert \otimes \nu_x(B).
\end{align*}
It results in particular that $(\varphi_1^{v^*})^{-1}_{\#}\mu' = \vert(\varphi_1^{v^*})^{-1}_{\#}\mu'\vert \otimes \nu_x$ and that for any measurable $E \subseteq \R^n$, $\vert(\varphi_1^{v^m})^{-1}_{\#}\mu'\vert(E) \rightarrow \vert(\varphi_1^{v^*})^{-1}_{\#}\mu'\vert(E)$. Now assume that $E\subseteq \R^n$ is a measurable subset such that $\vert\mu\vert(E)=0$. Then, for all $m\in \N$, using again the equality $(\varphi_1^{v^m})^{-1}_{\#}\mu' = (\alpha_1^{\tilde{\eta}^m})^2 \mu$ 
\begin{equation*}
 \vert(\varphi_1^{v^m})^{-1}_{\#}\mu'\vert(E) = \int_{E} (\alpha_1^{\tilde{\eta}^m})^2(x) d\vert\mu\vert(x) = 0 .
\end{equation*}
Therefore we get $\vert(\varphi_1^{v^*})^{-1}_{\#}\mu'\vert(E)=0$ from which we deduce that $\vert(\varphi_1^{v^*})^{-1}_{\#}\mu'\vert$ is absolutely continuous with respect to $\vert\mu\vert$. By Radon-Nykodym's theorem, there exists $\check{\alpha}^* \in L^1(\vert\mu\vert)$ such that $\vert(\varphi_1^{v^*})^{-1}_{\#}\mu'\vert = \check{\alpha}^* \vert\mu\vert$ and thus $(\varphi_1^{v^*})^{-1}_{\#}\mu' = \check{\alpha}^* \vert\mu\vert \otimes \nu_x = \check{\alpha}^* \mu$ which leads to $\mu' = (\varphi_1^{v^*})_{\#} (\check{\alpha}^* \mu)$. Setting for all $t\in [0,1]$, $\eta_t(x) = 2(\sqrt{\check{\alpha}^*} -1)$, we get $\eta \in L^2([0,1],L^2(\vert\mu\vert))$, $\alpha_t^{\eta}(x) \geq 0$ for $\vert\mu\vert$-a.e. $x \in \R^n$ and the previous equality is equivalent to $\mu' =(\varphi_1^{v^*})_{\#}((\tilde{\alpha}^{\eta}_1)^2\mu)=(\varphi_1^{v^*},(\alpha_1^{\eta})^2) \cdot \mu$.

It follows that we can apply Lemma \ref{lemma:LDDMM_FR_fixed_deform} and thus find $\bar{\eta}^* \in L^2([0,1],L^2(\vert\mu\vert))$ such that $\bar{\eta}^*$ minimizes: 
\begin{equation*}
 \int_0^1 \int_{\R^n \times \origrass}\eta^2_t(x) J_T\varphi_t^{v^*}(x) \, d\mu(x,T) dt
\end{equation*}
among all $\eta \in  L^2(\vert\mu\vert)$ such that $\tilde{\alpha}_t^{\eta}(x)\geq 0$ and $(\varphi_1^{v^*},(\tilde{\alpha}_1^{\eta})^2) \cdot \mu = \mu'$. Therefore we have on the one hand $(\varphi_1^{v^*},(\tilde{\alpha}_1^{\bar{\eta}^*})^2) \cdot \mu = \mu'$ and by construction, for all $m\in \N$, 
\begin{equation}
\label{eq:proof_existence_geodesics_LDDMM_FR1}
 \int_0^1 \int_{\R^n \times \origrass}\bar{\eta}^*_t(x)^2 J_T\varphi_t^{v^*}(x) \, d\mu(x,T) dt \leq \int_0^1 \int_{\R^n \times \origrass}\eta^m_t(x)^2 J_T\varphi_t^{v^*}(x) \, d\mu(x,T) dt
\end{equation}
Moreover since $\varphi_t^{v^m}$ converges to $\varphi_t^{v^*}$ in $\|\cdot\|_{1,\infty}$ uniformly over $t \in [0,1]$ and on a compact set that contains $\text{supp}(\vert\mu\vert)$, for any $\varepsilon>0$, there exists $p \in \mathbb{N}$ such that for all $m\geq p$, $t\in[0,1]$, $T \in \origrass$ and $x \in \text{supp}(\vert\mu\vert)$, we have $\vert J_{T} \varphi_t^{v^m}(x) - J_{T} \varphi_t^{v^*}(x)\vert\leq \varepsilon$ which also leads to $\vert J_{T} \varphi_t^{v^m}(x) - J_{T} \varphi_t^{v^*}(x)\vert \leq A \varepsilon J_{T} \varphi_t^{v^m}(x)$. Going back to \eqref{eq:proof_existence_geodesics_LDDMM_FR1}, we obtain:
\begin{align*}
 &\int_0^1 \int_{\R^n \times \origrass}\bar{\eta}^*_t(x)^2 J_T\varphi_t^{v^*}(x) \, d\mu(x,T) dt \\
 &\leq \int_0^1 \int_{\R^n \times \origrass}\eta^m_t(x)^2 J_T\varphi_t^{v^m}(x) \, d\mu(x,T) dt \\
 &\phantom{aa}+ \int_0^1 \int_{\R^n \times \origrass}\eta^m_t(x)^2 \left[J_T\varphi_t^{v^*}(x)- J_T\varphi_t^{v^m}(x)\right] \, d\mu(x,T) dt \\
 &\leq(1+A\varepsilon) \int_0^1 \int_{\R^n \times \origrass}\tilde{\eta}^m_t(x)^2 J_T\varphi_t^{v^m}(x) \, d\mu(x,T) dt 
\end{align*}
As this holds for all $\varepsilon>0$, we obtain that:
\begin{align*}
 &\int_0^1 \int_{\R^n \times \origrass}\bar{\eta}^*_t(x)^2 J_T\varphi_t^{v^*}(x) \, d\mu(x,T) dt \\
 &\leq \lim \inf_{m\rightarrow +\infty} \int_0^1 \int_{\R^n \times \origrass}\eta^m_t(x)^2 J_T\varphi_t^{v^m}(x) \, d\mu(x,T) dt
\end{align*}
which, combined with the weak lower semicontinuity of the squared Hilbert norm $v \mapsto \int_0^1 \|v_t\|_V^2 dt$, leads to:
\begin{equation*}
 C_{\mu}(v^*,\bar{\eta}^*) \leq \lim \inf_{m\rightarrow +\infty} C_{\mu}(v^m,\eta^m) = d_{\mathcal{D-FR}}(\mu,\mu').
\end{equation*}
In conclusion, $(v^*,\bar{\eta}^*)$ is a minimizer of the energy defining the distance between $\mu$ and $\mu'$.

\section{Proof of Propositions \ref{prop:geod_0_var} and \ref{prop:geod_1_var}} 
As we know from Section \ref{ssec:LDDMM_FR} that solutions to the geodesic boundary value problem between the two Diracs exist, we will derive necessary conditions satisfied by the solutions summoning the Pontryagin maximum principle (PMP) of optimal control \cite{Pontryagin1962} (more specifically the version derived in \cite{arguillere14:_shape} for the type of infinite-dimensional problems considered here). 

\vskip1ex 

\textbf{Proof of Proposition \ref{prop:geod_0_var}.} We start with the $0$-varifold case. The optimal control problem here simplifies to the minimization of:
\begin{equation*}
 C(v,\eta) = \frac{1}{2} \int_0^1 \|v_t\|_V^2 dt + \frac{\gamma}{2} \int_0^1 \eta_t^2 r_0 dt 
\end{equation*}
over $(v,\eta) \in L^2([0,1],V\times \R)$ subject to the state equations:
\begin{equation*}
 \left\lbrace\begin{aligned}
\dot{x}(t) &= v_t(x(t)) \\
\dot{\tilde{\alpha}}(t) &= \frac{1}{2} \eta_t
\end{aligned}\right.
\end{equation*}
and the boundary conditions $x(0)=x_0, \tilde{\alpha}(0)=1$ and $x(1)=x_1, \tilde{\alpha}(1)=\sqrt{r_1/r_0}$. Let $(v,\eta)$ be a solution. The Pontryagin maximum principle then states the existence of the costate functions $t \mapsto p_x(t) \in \R^n$ and $t\mapsto p_{\tilde{\alpha}}(t) \in \R$ that satisfy the adjoint equations: 
\begin{equation*}
 \left\lbrace\begin{aligned}
\dot{p}_x(t) &= -(d_{x(t)} v_t)^T p_x(t) \\
\dot{p}_{\tilde{\alpha}}(t) &= 0
\end{aligned}\right.
\end{equation*}
from which we immediately deduce that $p_{\tilde{\alpha}}$ is constant. Furthermore, the Hamiltonian of the system is here given by
\begin{equation*}
 H(x,\tilde{\alpha},p_x,p_{\tilde{\alpha}},v,\eta) = p_x^T v(x) + \frac{1}{2} p_{\tilde{\alpha}} \eta - \frac{1}{2} \|v\|_V^2 - \frac{\gamma}{2} r_0 \eta^2
\end{equation*}
and so the optimality conditions of the PMP that write $(v_t,\eta_t) = \text{argmin}_{(v',\eta')\in V \times \R} \  H(x(t),\tilde{\alpha}(t),p_x(t),p_{\tilde{\alpha}}(t),v',\eta')$ become here $v_t(\cdot) = K_V(x(t),\cdot) p_x(t)$ (c.f. \cite{arguillere14:_shape} for details on this derivation) and $\eta_t = \frac{p_{\tilde{\alpha}}}{2\gamma r_0}$. From this and the reproducing kernel formula, we first deduce that $\|v_t\|_V^2 = p_x(t)^T K_V(x(t),x(t)) p_x(t)$. Moreover, we get $\dot{x}(t) = K_V(x(t),x(t)) p_x(t)$ and so $p_x(t) = K_V(x(t),x(t))^{-1} \dot{x}(t)$. It follows that the first term in the cost function is equal to
\begin{equation*}
  \int_0^1 \dot{x}(t)^T G_{x(t)} \dot{x}(t) dt 
\end{equation*}
with $G_{x(t)} = K_V(x(t),x(t))^{-1}$. Therefore, as the second term of the cost is independent of $v$ and thus $\dot{x}(t)$, the path $t\mapsto x(t)$ minimizes the above energy subject to the boundary constraints $x(0)=x_0$ and $x(1)=x_1$ which means precisely that $x(t)$ follows the geodesic from $x_0$ to $x_1$ for the Riemannian metric on $\R^n$ given by the field of positive definite matrices $G_x = K_V(x,x)^{-1}$. On the other hand, we have:
\begin{equation*}
  \dot{\tilde{\alpha}}(t) = \frac{1}{2} \eta_t = \frac{p_{\tilde{\alpha}}}{4\gamma r_0} 
\end{equation*}
which, together with the boundary conditions $\tilde{\alpha}(0)=1$ and $\tilde{\alpha}(1)=\sqrt{r_1/r_0}$, leads to $\tilde{\alpha}(t) = (1-t) + t \sqrt{r_1/r_0}$ and so the varifold weight $r(t)$ is:
\begin{equation*}
    r(t) = \tilde{\alpha}(t)^2 r_0 = ((1-t)\sqrt{r_0} + t \sqrt{r_1})^2.
\end{equation*}

\vskip1ex

\textbf{Proof of Proposition \ref{prop:geod_1_var}.} Let us now move on to the 1-varifold case. For concision, we will write again $\tilde{\alpha}(t) = \tilde{\alpha}_t^\eta$ and define the auxiliary variable $\tilde{u}(t) = d_{x(t)}\varphi_t^v(r_0 u_0)$ so that $u(t)= \frac{\tilde{u}(t)}{\vert \tilde{u}(t)\vert}$ and $r(t) = \tilde{\alpha}(t)^2 \vert \tilde{u}(t) \vert$. The cost to minimize is then:
\begin{equation*}
 C(v,\eta) = \frac{1}{2} \int_0^1 \|v_t\|_V^2 dt + \frac{\gamma}{2} \int_0^1 \eta_t^2 \vert \tilde{u}(t) \vert dt
\end{equation*}
with the state equations:
\begin{equation*}
 \left\lbrace\begin{aligned}
\dot{x}(t) &= v_t(x(t)) \\
\dot{\tilde{u}}(t) &= d_{x(t)}v_t(\tilde{u}(t)) \\
\dot{\tilde{\alpha}}(t) &= \frac{1}{2} \eta_t
\end{aligned}\right.
\end{equation*}
and the boundary conditions $x(0)=x_0$, $x(1)=x_1$, $\tilde{u}(0)= r_0 u_0$, $\tilde{\alpha}(0)=1$, $\frac{\tilde{u}(1)}{\vert \tilde{u}(1) \vert}= u(1)$ and $\tilde{\alpha}(1)^2\vert \tilde{u}(1)\vert=r_1$. The Hamiltonian for this optimal control problem is now:
\begin{equation*}
 H(x,\tilde{u},\tilde{\alpha},p_x,p_{\tilde{u}},p_{\tilde{\alpha}},v,\eta) = p_x^T v(x) + p_{\tilde{u}}^T d_{x}v(\tilde{u}) + \frac{1}{2} p_{\tilde{\alpha}} \eta - \frac{1}{2} \|v\|_V^2 - \frac{\gamma}{2} \eta^2 \vert \tilde{u} \vert
\end{equation*}
from which we get the following adjoint equations:
\begin{equation}
\label{eq:adjoint_eq_1var}
 \left\lbrace\begin{aligned}
\dot{p}_x(t) &= -(d_{x(t)} v_t)^T p_x(t) - d^2_{x(t)} v_t(\tilde{u}(t),\cdot)^T p_{\tilde{u}}(t)\\
\dot{p}_{\tilde{u}}(t) &= -(d_{x(t)} v_t)^T p_{\tilde{u}}(t) + \frac{\gamma}{2} \eta_t^2 \frac{\tilde{u}(t)}{\vert\tilde{u}(t) \vert} \\
\dot{p}_{\tilde{\alpha}}(t) &= 0
\end{aligned}\right.
\end{equation}
where $d^2_x v(u,\cdot)$ denotes the matrix of the linear mapping $h \in \R^n \mapsto d^2_x v(u,h) \in \R^n$. Since $p_{\tilde{\alpha}}(t)$ is constant from the last equation above, we will simply write $p_{\tilde{\alpha}}$ in what follows. The optimality conditions of the PMP lead to the following expressions for the optimal controls:
\begin{equation}
\label{eq:expr_opt_control1}
\left\lbrace\begin{aligned}
 &v_t(\cdot) = K_V(x(t),\cdot) p_x(t) + \partial_1 K_V(x(t),\cdot)(\tilde{u}(t),p_{\tilde{u}}(t)) \\
 &\eta_t = \frac{p_{\tilde{\alpha}}}{2 \gamma \vert \tilde{u}(t) \vert}
 \end{aligned}\right.
\end{equation}
Now, using the expression of the kernel $K_V(x,y) = \rho\left(\frac{\vert x-y \vert^2}{\sigma^2} \right)$, the first equation in \eqref{eq:expr_opt_control1} gives for all $x \in \R^n$:
\begin{equation}
\label{eq:expr_opt_v_2}
   v_t(x) = \rho\left(\frac{\vert x-x(t)\vert^2}{\sigma^2}\right) p_x(t) - \frac{2}{\sigma^2}\rho'\left(\frac{\vert x-x(t)\vert^2}{\sigma^2}\right) [(x-x(t))^T \tilde{u}(t)] p_{\tilde{u}}(t).
\end{equation}
In particular, we get that $\dot{x}(t) = v_t(x(t)) = p_x(t)$. Furthermore, differentiating the above equation at $x =x(t)$, we find that for all $h \in \R^n$:
\begin{align*}
 &d_{x(t)}v_t(h) = \frac{1}{\tau} \left[\tilde{u}(t)^T h \right] p_{\tilde{u}}(t) \\ 
 &d_{x(t)}v_t^T(h) = \frac{1}{\tau} \left[p_{\tilde{u}}(t)^T h \right] \tilde{u}(t) \\
 &d^2_{x(t)} v_t(\cdot,\bar{u}(t))^T h = -\frac{1}{\tau} (p_x(t)^T h) \tilde{u}(t) 
\end{align*}
where we remind the reader that $\tau = -\frac{\sigma^2}{2\rho'(0)}>0$ and thus 
\begin{align*}
  \dot{p}_x(t) = -\frac{1}{\tau} \left[p_{\tilde{u}}(t)^T p_x(t) \right] \tilde{u}(t)  + \frac{1}{\tau} (p_x(t)^T p_{\tilde{u}}(t)) \tilde{u}(t)  = 0.
\end{align*}
Therefore $p_x(t)$ is constant and $\ddot{x}(t) =0$ which, with the two boundary conditions, leads to $p_x = x_1 - x_0$ and $x(t) = (1-t)x_0 + t x_1$. \\
We next analyze the behaviour of the direction $u(t)$. Rewriting the state and adjoint equations based on \eqref{eq:adjoint_eq_1var} and \eqref{eq:expr_opt_control1}, we have:
\begin{equation*}
\left\lbrace\begin{aligned}
 \dot{\tilde{u}}(t) &= \frac{1}{\tau} \vert\tilde{u}(t)\vert^2 p_{\tilde{u}}(t)  \\
 \dot{p}_{\tilde{u}}(t) &= -\frac{1}{\tau} \vert p_{\tilde{u}}(t) \vert^2 + \frac{p_{\tilde{\alpha}}^2}{8\gamma} \frac{\tilde{u}(t)}{\vert\tilde{u}(t)\vert^3}
 \end{aligned}\right.
\end{equation*}
Going back to $u(t)= \frac{\tilde{u}(t)}{\vert \tilde{u}(t) \vert}$, one has:
\begin{align*}
    \dot{u}(t) &= \frac{\dot{\tilde{u}}(t)}{\vert \tilde{u}(t) \vert} - \frac{1}{\vert\tilde{u}(t)\vert^3} \left[\dot{\tilde{u}}(t)^T \tilde{u}(t) \right] \tilde{u}(t) \\
    &=\frac{1}{\tau} \vert \tilde{u}(t) \vert \left(p_{\tilde{u}}(t) - \left[u(t)^T p_{\tilde{u}}(t)\right] u(t)\right) =\frac{1}{\tau} \vert \tilde{u}(t) \vert \text{Proj}_{u(t)^\bot}(p_{\tilde{u}}(t))
\end{align*}
with $\text{Proj}_{u(t)^\bot}$ denoting the orthogonal projector onto the hyperplane normal to $u(t)$. Letting $p_u(t) \doteq \vert \tilde{u}(t) \vert \text{Proj}_{u(t)^\bot}(p_{\tilde{u}}(t))$, we find after calculations that $\dot{p}_u(t) = -\frac{1}{\tau} \vert p_u(t)\vert^2 u(t)$. Therefore, we have the following coupled system of ODEs on $(u(t),p_u(t))$:
\begin{equation}
\label{eq:ODE_u_pu}
\left\lbrace\begin{aligned}
 \dot{u}(t) &= \frac{1}{\tau} p_u(t)  \\
 \dot{p}_{u}(t) &= -\frac{1}{\tau} \vert p_u(t)\vert^2 u(t)
 \end{aligned}\right.
\end{equation}
As $p_u(t)$ is orthogonal to $u(t)$ by definition, the above equations imply in particular that $\vert p_u(t)\vert^2= \kappa^2$ is constant. It also follows immediately that for all $t$, $(u(t),p_u(t))$ stay in the plane spanned by $\{u(0),p_u(0)\}$. Let us therefore identify the plane $\text{Span}\{u(0),p_u(0)\}$ with $\C$, choosing without loss of generality $u(0)=1$. We can then write $u(t) = e^{i\theta(t)}$ with $\theta(0) = 0, \theta(1) = \theta_{0,1}$ the angle between $u_0$ and $u_1$ and $p_{u}(t) = \pm i \kappa e^{i\theta(t)}$. With this identification, \eqref{eq:ODE_u_pu} leads to $\dot{\theta}(t) = \pm \frac{\kappa}{\tau}$, in other words the direction $\bar{u}(t)$ rotates with constant angular velocity. With the boundary conditions, it shows that $\kappa=\tau \theta_{0,1}$ and thus $\theta(t) = t\theta_{0,1}$. Eventually, assuming $u_1 \neq u_0$, this leads to the expression of $u(t)$:
\begin{equation}
 \label{eq:opt_path_u}
 u(t) = 
  \left\lbrace\begin{aligned}
 &\frac{1}{\sin(\theta_{0,1})} (\sin((1-t)\theta_{0,1}) u_0 + \sin(t\theta_{0,1}) u_1) \ \ \text{if } \theta_{0,1} \neq 0 \\
 &u_0 (= u_1) \ \ \text{if } \theta_{0,1} = 0
 \end{aligned}\right.
\end{equation}
with $\theta_{0,1} = \arccos(u_0^Tu_1)$. Note that when $u_1 = - u_0$, the geodesic is not unique as one can rotate from $u_0$ to $u_1$ in infinitely many ways. 

We are left with determining $r(t) = \tilde{\alpha}(t)^2 \vert \tilde{u}(t)\vert$. Let us start by introducing the auxiliary function $m(t) = \sqrt{r(t)} = \tilde{\alpha}(t) \sqrt{\vert \tilde{u}(t)\vert}$. Using the state and costate equations, we first see that $\dot{m}(t) = \frac{p_{\tilde{\alpha}}}{4\gamma\sqrt{\vert \tilde{u}(t)\vert}} + \frac{1}{2\tau} m(t) (\tilde{u}(t)^T p_{\tilde{u}}(t))$. Differentiating a second time, we get:
\begin{align*}
    \ddot{m}(t) &= -\frac{p_{\tilde{\alpha}}}{8\tau \gamma \sqrt{\vert \tilde{u}(t)\vert}} (\tilde{u}(t)^T p_{\tilde{u}}(t)) + \frac{p_{\tilde{\alpha}}}{8\tau \gamma \sqrt{\vert \tilde{u}(t)\vert}} (\tilde{u}(t)^T p_{\tilde{u}}(t)) \\
    &\phantom{aa}+\frac{1}{4 \tau^2}  (\tilde{u}(t)^T p_{\tilde{u}}(t))^2 m(t) + \frac{1}{2 \tau^2} \vert \tilde{u}(t)\vert^2 \vert p_{\tilde{u}}(t)\vert^2 m(t) \\
    &\phantom{aa}- \frac{1}{2 \tau^2} \vert \tilde{u}(t)\vert^2 \vert p_{\tilde{u}}(t)\vert^2 m(t) + \frac{p_{\tilde{\alpha}}^2}{16\tau \gamma} \frac{m(t)}{\vert \tilde{u}(t)\vert} \\
    &=\frac{1}{4\tau} \left[\frac{1}{\tau} (\tilde{u}(t)^T p_{\tilde{u}}(t))^2 + \frac{p_{\tilde{\alpha}}^2}{4 \gamma \vert \tilde{u}(t) \vert} \right] m(t).
\end{align*}
Moreover, it turns out that the term inside brackets is constant in time. Indeed, on the one hand we have $\frac{d}{dt} \tilde{u}(t)^T p_{\tilde{u}}(t) = \frac{p_{\tilde{\alpha}}^2}{8\gamma\vert \tilde{u}(t) \vert}$ so that:
\begin{equation*}
    \frac{d}{dt} (\tilde{u}(t)^T p_{\tilde{u}}(t))^2 = \frac{p_{\tilde{\alpha}}^2}{4\gamma \vert \tilde{u}(t)\vert} (\tilde{u}(t)^T p_{\tilde{u}}(t))
\end{equation*}
and on the other:
\begin{equation*}
    \frac{d}{dt} \frac{p_{\tilde{\alpha}}^2}{4 \gamma \vert \tilde{u}(t) \vert} =  -\frac{p_{\tilde{\alpha}}^2}{4\gamma\vert \tilde{u}(t)\vert^3} (\tilde{u}(t)^T \dot{\tilde{u}}(t)) = \frac{p_{\tilde{\alpha}}^2}{4\gamma \vert \tilde{u}(t)\vert} (\tilde{u}(t)^T p_{\tilde{u}}(t)). 
\end{equation*}
Let us therefore set $\nu^2 \doteq \frac{1}{4\tau} \left[\frac{1}{\tau} (\tilde{u}(0)^T p_{\tilde{u}}(0))^2 + \frac{p_{\tilde{\alpha}}^2}{4 \gamma \vert \tilde{u}(0) \vert} \right] >0$ so that $\ddot{m}(t) = \nu^2 m(t)$ and therefore $r(t) = m(t)^2$ takes the form $r(t) = C_0 \sinh(C_1 + \nu t)^2$ for some constants $C_0,C_1 \in \R$. Based on the two boundary conditions for $r(t)$ and hyperbolic trigonometry identities, we find eventually that:
\begin{equation}
\label{eq:expr_r_1var}
    r(t) = \left(\frac{\sqrt{r_0} \sinh((1-t)\nu) + \sqrt{r_1}\sinh(\nu t)}{\sinh(\nu)}\right)^2.
\end{equation}

To recover the explicit expression of $\nu$, let us first express $\dot{r}(0)$ based on \eqref{eq:expr_r_1var}. We obtain 
\begin{equation}
\label{eq:expr_dotr_1}
    \dot{r}(0) = 2\nu r_0 \left(\sqrt{\frac{r_1}{r_0}} \frac{1}{\sinh(\nu)} -\coth(\nu) \right).
\end{equation}
On the other hand, $\dot{r}(0) = 2 m(0) \dot{m}(0)$ and using the expression of $\dot{m}(t)$ obtained earlier, we get
\begin{equation}
\label{eq:expr_dotr_2}
    \dot{r}(0) = 2\sqrt{r_0} \left(\frac{p_{\tilde{\alpha}}}{4\gamma \sqrt{r_0}} + \frac{1}{2\tau} \sqrt{r_0} (\tilde{u}(0)^T p_{\tilde{u}}(0)) \right)
\end{equation}
Furthermore, the term $(\tilde{u}(0)^T p_{\tilde{u}}(0))$ may be expressed with respect to $p_{\tilde{\alpha}}$. Indeed, consider the function $h(t) = \frac{\tilde{\alpha}(t)p_{\tilde{\alpha}}}{2} - (\tilde{u}(t)^T p_{\tilde{u}}(t))$. We see that:
\begin{equation*}
    \dot{h}(t) = \frac{p_{\tilde{\alpha}}^2}{8\gamma\vert \tilde{u}(t) \vert} - \frac{p_{\tilde{\alpha}}^2}{8\gamma\vert \tilde{u}(t) \vert} =0
\end{equation*}
and thus $h(t)=h(0)=h(1)$. Now one of the terminal constraint for the optimal control problem is $\tilde{\alpha}(1)^2\vert \tilde{u}(1)\vert=r_1$ and therefore the transversality condition of the PMP yields that the vector $(p_{\tilde{u}}(1),p_{\tilde{\alpha}})^T$ is parallel to $\left(\tilde{\alpha}(1)^2\frac{\tilde{u}(1)}{\vert \tilde{u}(1) \vert}, 2 \tilde{\alpha}(1)\vert \tilde{u}(1) \vert\right)^T$ from which we deduce that $\tilde{u}(1)^T p_{\tilde{u}}(1) = \frac{\tilde{\alpha}(1)p_{\tilde{\alpha}}}{2}$, in other words $h(1) = h(0) =0$. Since $\tilde{\alpha}(0)=1$, this implies that $\tilde{u}(0)^T p_{\tilde{u}}(0) = \frac{p_{\tilde{\alpha}}}{2}$ from which we can rewrite $\nu^2$ as:
\begin{equation}
\label{eq:expr_nu_2}
    \nu^2 = \left(\frac{p_{\tilde{\alpha}}}{4\tau} \right)^2 \left[1+\frac{\tau}{\gamma r_0} \right] \implies \frac{p_{\tilde{\alpha}}}{\nu} = \pm \frac{4 \tau}{\sqrt{1+\frac{\tau}{\gamma r_0}}}.
\end{equation}
Moreover, \eqref{eq:expr_dotr_2} becomes 
\begin{equation*}
    \dot{r}(0) = \frac{p_{\tilde{\alpha}}}{2\gamma} \left(1+\frac{\gamma}{\tau} r_0 \right)
\end{equation*}
which combined with \eqref{eq:expr_dotr_1} and \eqref{eq:expr_nu_2} leads to the implicit equation on $\nu$:
\begin{equation}
\label{eq:implicit_nu}
\sqrt{\frac{r_1}{r_0}} \frac{1}{\sinh(\nu)} -\coth(\nu) = \pm \underbrace{\sqrt{1+\frac{\tau}{\gamma r_0}}}_{\doteq \chi}.
\end{equation}
A quick analysis of the function on the left hand side shows that this equation only has a solution with $+\chi$ on the right hand side when $r_1>r_0$ and with $-\chi$ when $r_1\leq r_0$. Using Mathematica, we find specifically the explicit expressions for $\nu$ given in Proposition \ref{prop:geod_1_var}. 

Finally, we can express the distance between $r_0 \delta_{x_0,u_0}$ and $r_1 \delta_{x_1,u_1}$. First, the kernel norm $\|v_t\|_V^2$ can be calculated based on the expressions \eqref{eq:expr_opt_control1} and \eqref{eq:expr_opt_v_2} using the reproducing kernel property for kernel derivatives (\cite{glaunes2014matrix} Theorem 2.11). Skipping some of the details for brevity, we obtain:
\begin{align*}
    \|v_t\|_V^2 &= \vert p_x \vert^2 +\frac{1}{\tau} \vert p_{\tilde{u}}(t) \vert^2  \vert \tilde{u}(t) \vert^2 \\
    &= \vert x_1 - x_0 \vert^2 +  \frac{1}{\tau} \left(\vert \text{Proj}_{u(t)^\bot}(p_{\tilde{u}}(t)) \vert^2 + (u(t)^T p_{\tilde{u}}(t))^2  \right) \vert \tilde{u}(t) \vert^2  \\
    &= \vert x_1 - x_0 \vert^2 + \frac{1}{\tau} \vert p_{u}(t) \vert^2 + \frac{1}{\tau} (\tilde{u}(t)^T p_{\tilde{u}}(t))^2 
\end{align*}
where we have used the fact that $p_x = x_1 - x_0$, $\tilde{u}(t) = u(t)/\vert u(t) \vert$ and $p_{u}(t) = \vert \tilde{u}(t) \vert \text{Proj}_{u(t)^\bot}(p_{\tilde{u}}(t)$. Furthermore, by the results above, we know that $\vert p_{u}(t) \vert^2 = \kappa^2 = (\tau \arccos(u_0^Tu_1))^2$. Therefore:
\begin{align*}
    \|v_t\|_V^2 &= \vert x_1 - x_0 \vert^2 + \tau \arccos(u_0^Tu_1)^2 + \frac{1}{\tau} (\tilde{u}(t)^T p_{\tilde{u}}(t))^2.
\end{align*}
Now for the second term in the energy, we have using \eqref{eq:expr_opt_control1}:
\begin{align*}
    \frac{\gamma}{2} \eta_t^2 \vert \tilde{u}(t) \vert = \frac{p_{\tilde{\alpha}}^2}{8 \gamma \vert \tilde{u}(t) \vert}
\end{align*}
leading to:
\begin{align*}
    \frac{1}{2}\|v_t\|_V^2 + \frac{\gamma}{2} \eta_t^2 \vert \tilde{u}(t) \vert  &= \frac{\vert x_1 - x_0 \vert^2}{2} + \frac{\tau}{2} \arccos(u_0^Tu_1)^2 \\
    &\phantom{aa}+ \frac{1}{2} \left(\frac{1}{\tau} (\tilde{u}(t)^T p_{\tilde{u}}(t))^2 + \frac{p_{\tilde{\alpha}}^2}{4 \gamma \vert \tilde{u}(t) \vert} \right)
\end{align*}
Recalling that the last term is constant in time and equal to $2 \tau \nu^2$ thanks to our earlier derivations, we finally get:
\begin{align*}
    d_{V-\mathcal{FR}}(r_0 \delta_{(x_0,u_0)},r_1 \delta_{(x_1,u_1)})^2 &=\int_0^1  \left(\frac{1}{2}\|v_t\|_V^2 + \frac{\gamma}{2} \eta_t^2 \vert \tilde{u}(t) \vert \right) dt\\
    &=\frac{\vert x_1 - x_0 \vert^2}{2} + \frac{\tau}{2} \arccos(u_0^Tu_1)^2 + 2 \tau \nu^2.
\end{align*}
% As we have an autonomous control system with time-independent Lagrangian, the Hamiltonian is conserved over time namely for all $t \in [0,1]$:
% \begin{align*}
%   &H(x(t),\tilde{u}(t),\tilde{\alpha}(t),p_x(t),p_{\tilde{u}}(t),p_{\tilde{\alpha}}(t),v_t,\eta_t) \\
%   &= H(x(0),\tilde{u}(0),\tilde{\alpha}(0),p_x(0),p_{\tilde{u}}(0),p_{\tilde{\alpha}}(0),v_0,\eta_0) \\
%   &= p_x(0)^T v(x_0) + p_{\tilde{u}}(0)^T d_{x_0}v_0(\tilde{u}(0)) + \frac{1}{2} p_{\tilde{\alpha}} \eta_0 -\frac{1}{2} \|v_0\|_V^2 - \frac{\gamma}{2} \eta_0^2 \vert \tilde{u}(0)\vert \\
%   &=
% \end{align*}

\section{Proof of Theorem \ref{thm:existence_relaxed_L2}}
\label{appendix:relaxed_L2}
Let us denote by $C(v,\alpha)$ the functional to minimize in \eqref{eq:relaxed_LDDMM_L2_distance} and let $(v^j,\alpha^j)$ be a minimizing sequence. Using a similar argument as in the previous existence proofs, since $(v^j)$ is bounded in $L^2([0,1],V)$, we can assume (by extracting a subsequence if necessary) that $(v^j)$ converges weakly to some $\bar{v}$ in $L^2([0,1],V)$ giving that $\|(\varphi_1^{v_j} - \varphi_1^{\bar{v}})\vert_K \|_{1,\infty} \rightarrow 0$ as $j\rightarrow +\infty$. Also, since $\gamma>0$, $(\alpha^j)$ is a bounded sequence in $L^2(\vert\mu\vert)$ and thus up to extraction of another subsequence, we may assume that we have weak convergence $\alpha^j \rightharpoonup \bar{\alpha}$ to some $\bar{\alpha}\in L^2(\vert\mu\vert)$.

From the weak lower semicontinuity of the first two terms in the energy, we deduce: 
\begin{align}
    &\int_0^1 \|\bar{v}\|_V^2 dt \leq \liminf_{j\rightarrow +\infty} \int_0^1 \|v^j\|_V^2 dt \label{eq:proof_prop_LDDMM_L2_1} \\
    &\int_{\R^n}(\bar{\alpha}(x)-1)^2 d\vert\mu\vert(x) \leq \liminf_{j\rightarrow +\infty} \int_{\R^n}(\alpha^j(x)-1)^2 d\vert\mu \vert(x) \label{eq:proof_prop_LDDMM_L2_2}
\end{align}
In addition, since $\|\varphi_1^{v_j}\|_{1,\infty}$ is bounded and $\mu$ is supported in the compact set $K$, there exists another compact subset $K\subset K' \subset \R^n \times \origrass$  such that for all $j\in \N$, $\text{supp}( (\varphi_1^{v^j})_{\#} (\alpha^j\mu)) \subset K'$. Moreover, using the disintegration theorem on the measure $\mu$ (c.f. Proposition \ref{prop:disintegration}), we see that:
\begin{align*}
\vert (\varphi_1^{v^j})_{\#}(\alpha^j\mu)\vert(\R^n) &= \int_{\R^n \times \origrass} \alpha^j(x) \vert J_U \varphi_1^{v^j}(x)\vert d\mu(x,U) \\
&=\int_{\R^n} \alpha^j(x) \left(\int_{\origrass} \vert J_U \varphi_1^{v^j}(x)\vert d\nu_x(U) \right) d\vert \mu\vert(x) 
\end{align*}
For the same reason as above, we have that $\vert J_U \varphi_1^{v^j}(x)\vert$ is bounded uniformly over $x \in \text{supp}(\vert\mu\vert)\subset K$, $j\in \N$ and $U \in \origrass$ from which we get, applying the Cauchy-Schwarz inequality:
\begin{align*}
\vert (\varphi_1^{v^j})_{\#} (\alpha^j\mu)\vert(\R^n) &\leq C \sqrt{\int_{\R^n} \alpha^j(x)^2 d\vert \mu\vert(x)} = C \|\alpha^j\|_{L^2(\vert\mu \vert)}
\end{align*}
for some constant $C>0$. Now since $\|\alpha^j\|_{L^2(\vert\mu \vert)}$ is bounded, we deduce that there exists $M>0$ such that $\vert\varphi_1^{v^j})_{\#} (\alpha^j\mu)\vert(\R^n)\leq M$ for all $j\in \N$. In other words, we have obtained that for all $j$, $\varphi_1^{v^j})_{\#} (\alpha^j\mu)$ belong to the space $\mathcal{V}_{d,M,K'}$ defined in Proposition \ref{prop:dW_metrization}. To show the convergence of $(\varphi_1^{v^j})_{\#} (\alpha^j\mu)$ for the $W^*$ metric, we are thus left to show that it converges for the weak-$^*$ topology. 

Thus, let $\omega \in C_c(\R^n \times \origrass)$. We have: 
\begin{align*}
\left( (\varphi_1^{v^j})_{\#}(\alpha^j \mu) \vert \omega \right) = \left((\varphi_1^{\bar{v}})_{\#}(\alpha^j \mu)\vert \omega \right) + \left((\varphi_1^{v^j})_{\#}(\alpha^j \mu) - (\varphi_1^{\bar{v}})_{\#} (\alpha^j \mu)\vert \omega \right)
\end{align*}
Looking at the first term, we see that:
\begin{align*}
  &\left((\varphi_1^{\bar{v}})_{\#}(\alpha^j \mu) \vert \omega \right) = \int_{\R^n \times \origrass} \alpha^j(x) J_U\varphi_1^{\bar{v}}(x) \omega(\varphi_1^{\bar{v}}(x),d_x \varphi_1^{\bar{v}} \cdot U) d\mu(x,U) \\
  &=\int_{\R^n} \alpha^j(x) \left(\int_{\origrass} J_U\varphi_1^{\bar{v}}(x) \omega(\varphi_1^{\bar{v}}(x),d_x \varphi_1^{\bar{v}} \cdot U) d\nu_x(U) \right) d\vert \mu \vert(x) \\
  &\xrightarrow[j\rightarrow \infty]{} \int_{\R^n} \bar{\alpha}(x) \left(\int_{\origrass} J_U\varphi_1^{\bar{v}}(x) \omega(\varphi_1^{\bar{v}}(x),d_x \varphi_1^{\bar{v}} \cdot U) d\nu_x(U) \right) d\vert \mu \vert(x) \\
  &=\left((\varphi_1^{\bar{v}})_{\#} (\bar{\alpha} \mu) \vert \omega \right)
\end{align*}
where the convergence in the third row follows from the weak convergence of $\alpha^j$ to $\bar{\alpha}$ in $L^2(\vert\mu\vert)$ and the fact the function between parentheses is measurable and bounded on $\R^n$ and thus in $L^2(\vert \mu \vert)$. As for the second term, we can expand it and see that:
\begin{align*}
&\left\vert \left( (\varphi_1^{v^j})_{\#} (\alpha^j\mu_0) - (\varphi_1^{\bar{v}})_{\#} (\alpha^j\mu))\vert \omega \right) \right\vert \\
&\leq \int_{K} \vert \alpha^j(x)\vert \left\vert J_U \varphi_1^{v^j}(x) \omega(\varphi_1^{v^j}(x),d_x \varphi_1^{v^j} \cdot U) -  J_U\varphi_1^{\bar{v}}(x) \omega(\varphi_1^{\bar{v}}(x),d_x \varphi_1^{\bar{v}} \cdot U) \right\vert \\ 
&\phantom{\leq \int_{K} \vert \alpha^j(x)\vert \left\vert J_U \varphi_1^{v^j}(x) \omega(\varphi_1^{v^j}(x),d_x \varphi_1^{v^j} \cdot U) -  J_U\varphi_1^{\bar{v}}(x) \omega(\varphi_1^{\bar{v}}(x),d_x\right.}d\mu(x,U) \\
&\leq \|\alpha^j\|_{L^2(\vert\mu\vert)}. \\
&\phantom{a} \left(\int_{K} \left\vert J_U \varphi_1^{v^j}(x) \omega(\varphi_1^{v^j}(x),d_x \varphi_1^{v^j} \cdot U) -  J_U\varphi_1^{\bar{v}}(x) \omega(\varphi_1^{\bar{v}}(x),d_x \varphi_1^{\bar{v}} \cdot S) \right\vert^2 d\mu(x,U) \right)^{1/2}
\end{align*}
Now, $\|\alpha^j\|_{L^2(\vert\mu\vert)}$ is bounded and by convergence of $\varphi_1^{v^j}$ to $\varphi_1^{\bar{v}}$ in $\|\cdot\|_{1,\infty}$ and uniform continuity of $\omega$, we deduce that 
\begin{equation*}
 \left( (\varphi_1^{v^j})_{\#} (\alpha^j\mu_0) - (\varphi_1^{\bar{v}})_{\#} (\alpha^j\mu))\vert \omega \right) \rightarrow 0.
\end{equation*}
Therefore, we have shown that $(\varphi_1^{v^j})_{\#} (\alpha^j\mu) \overset{\ast}{\rightharpoonup} (\varphi_1^{\bar{v}})_{\#}(\bar{\alpha}\mu)$. By Proposition  \ref{prop:dW_metrization}, this implies that $(\varphi_1^{v^j})_{\#}(\alpha^j\mu) \xrightarrow[]{d_{W^*}} (\varphi_1^{\bar{v}})_{\#}(\bar{\alpha}\mu)$. As a result, 
\begin{equation}
\label{eq:proof_prop_LDDMM_L2_3}
\lim_{j\rightarrow \infty} \| (\varphi_1^{v^j})_{\#}(\alpha^j\mu) - \mu' \|_{W^*}^2 = \| (\varphi_1^{\bar{v}})_{\#} (\bar{\alpha}\mu) - \mu' \|_{W^*}^2.
\end{equation}

Combining \eqref{eq:proof_prop_LDDMM_L2_1}, \eqref{eq:proof_prop_LDDMM_L2_2} and \eqref{eq:proof_prop_LDDMM_L2_3}, we obtain:
\begin{equation*}
    C(\bar{v},\bar{\alpha})\leq \lim \inf_{j\rightarrow \infty} C(v^j,\alpha^j)
\end{equation*}
and thus $(\bar{v},\bar{\alpha})$ is a minimizer of $C$.

\section{Proof of Theorem \ref{thm:existence_solutions_relaxed_LDDMM_FR}}
\label{appendix:relaxed_FR}
Let us again denote by $C(v,\eta)$ the functional to minimize in \eqref{eq:relaxed_LDDMM_FR_distance} and let $(v^j,\eta^j)$ be a minimizing sequence. Since $\mu = \sum_{i=1}^{N} r_i \delta_{(x_i,U_i)}$ is here assumed to be a discrete varifold, we can equivalently view the functions $\eta^j$ as functions in $L^2([0,1],\R^N)$ in  which each component of the vector $\eta^j(t)$ is associated to the corresponding Dirac in $\mu$. Similarly, $\tilde{\alpha}_1^{\eta^j}$ can be viewed as a vector in $\R_{+}^N$ and to simplify the following derivations, we shall write with a slight abuse of notations $\tilde{\alpha}_1^{\eta^j} = \alpha^j=(\alpha_i^j)_{i=1,\ldots,N}$. Then for each $j\in \N$, define $\mu^j = \sum_{i=1}^{N} r^j_i \delta_{(x^j_i,U^j_i)}$ to be the transformed varifold at $t=1$ for $(v^j,\eta^j)$ i.e. $x_i^j = \varphi_1^{v^j}(x_i)$, $U_i^j = d_{x_i}\varphi_1^{v^j} \cdot U_i$ and $r_i^j = (\alpha_i^j)^2 J_{U_i}\varphi_1^{v^j} r_i$. 

We will first show that, modulo extraction of subsequences, for each $i=1,\ldots,N$ the sequences $(x_i^j)$, $(U_i^j)$ and $(r_i^j)$ are converging in $\R^n$, $\origrass$ and $\R_+$ respectively. Let $K$ be a compact subset of $\R^n$ that contains $x_1,\ldots,x_N$. Since $C(v^j,\eta^j)$ and thus $\int_0^1 \|v^j_t\|_V^2 dt$ is bounded uniformly in $j$, using similar arguments as in the previous proofs, we have, up to extraction of a subsequence, that $v^j$ converges weakly to some $v^* \in L^2([0,1],V)$ and thus $\varphi_1^{v^j}$ converges to $\varphi_1^{v*} \in \textrm{Diff}(\R^n)$ in $\|\cdot\|_{1,\infty}$ on every compact subset of $\R^n$ in particular on $K$ leading to $x_i^j \xrightarrow[j \rightarrow \infty]{} \varphi_1^{v*}(x_i)\doteq x_i^*$ and $U_i^j \xrightarrow[j \rightarrow \infty]{} d_{x_i}\varphi_1^{v*} \cdot U_i \doteq U_i^*$ for all $i=1,\ldots,N$. Now we also have that $\|\mu^j - \mu^{tar}\|_{W^*}$ is bounded uniformly in $j$ and furthermore:
\begin{align*}
    \|\mu^j -\mu^{tar} \|_{W^*} 
    &\geq \| \mu^j  \|_{W^*} -\| \mu^{tar}  \|_{W^*} \\ 
    &=\sqrt{\sum_{i,i'=1}^N r_i^j r_{i'}^j k(x_i^j,U_i^j,x_{i'}^j,U_{i'}^j)} - \| \mu^{tar} \|_{W^*} \\
    &\geq \sqrt{\sum_{i=1}^N (r_i^j)^2 k(x_i^j,U_i^j,x_i^j,U_i^j)} - \| \mu^{tar} \|_{W^*}
    %&= \sqrt{\sum_{i,j=1}^N r_i^m r_j^m \rho(|x_i^m - x_j^m|^2) \gamma(\langle T_i^m,T_j^m \rangle)} - \| \mu^{tar} \|_{W^*} \\
    %&\geq \sqrt{\sum_{i=1}^N (r_i^m)^2 \rho(0) \gamma(1)} - \| \mu^{tar} \|_{W^*}
\end{align*}
where the last equality follows from the fact that $k(x_i^j,U_i^j,x_{i'}^j,U_{i'}^j) \geq 0$ by assumption on the kernel of $W$. Moreover, $(x,U) \mapsto k(x,U,x,U)$ is strictly positive, continuous and we can find a compact subset $K'\subset \R^n$ with $(x_i^j,U_i^j)$ belonging to the compact $K' \times \origrass$ for all $i=1,\ldots,N$ and $j \in \N$. Consequently, taking $\delta = \min \{k(x,U,x,U) \ \vert \ (x,U) \in K' \times \origrass\} >0$ we have $k(x_i^j,U_i^j,x_i^j,U_i^j)\geq \delta$ for all $i$ and $j$ and therefore we must have that each sequence $(r_i^j)_{j \in \mathbb{N}}$ is bounded since otherwise we could find a subsequence making the right hand side in the above inequality go to infinity. Thus, up to extracting once again a subsequence, we have $r_i^j \xrightarrow[j \rightarrow \infty]{} r_i^*$ with $r_i^* \geq 0$ for all $i=1,\ldots,N$. 

Let us now set $\mu^* = \sum_{i=1}^{N} r_i^* \delta_{(x_i^*,T_i^*)}$ and define $\eta^*\in L^2([0,1],\R^N)$ as given by Lemma \ref{lemma:LDDMM_FR_fixed_deform} which means here specifically:
\begin{equation*}
 \eta^*_{i}(t) = 2 \frac{\sqrt{\frac{r_i^*}{r_i h_i(1)}}-1}{h_i(t)\int_0^1 1/h_i(s)ds}
\end{equation*} 
where $h_i(t) = J_{U_i}\varphi_t^{v^*}(x_i)$. Then it is easy to check that $\mu^* = (\varphi_1^{v^*},(\alpha_1^{\eta^*})^2) \cdot \mu$ and thus $\mu^* \in \bar{\Theta}(\mu)$. Then by the same argument as in the proof of Theorem \ref{thm:geodesics_LDDMM_FR} (c.f. Appendix \ref{appendix:proof_geod_FR}), we get that:
\begin{align*}
    &\frac{1}{2} \int_0^1 \|v_t^*\|_V^2 dt + \frac{\gamma}{2} \sum_{i=1}^{N}\int_0^1 \eta^*_i(t)^2 J_{U_i} \varphi_t^{v^*}(x_i) r_i dt \\
    &\leq \lim \inf_{j\rightarrow \infty} \left( \frac{1}{2} \int_0^1 \|v_t^j\|_V^2 dt + \frac{\gamma}{2} \sum_{i=1}^{N}\int_0^1 \eta^j_i(t)^2 J_{U_i} \varphi_t^{v^j}(x_i) r_i dt \right).
\end{align*}

Thus it only remains to examine the convergence of the kernel fidelity term for which it is enough to show that the sequence $(\mu^j)$ converges to $\mu^*$ for $\|\cdot\|_{W^*}$. Let $\omega \in W$ such that $\|\omega\|_{W}\leq 1$. Given the continuous embedding assumption, we also have $\|\omega\|_{1,\infty} \leq c_W \|\omega\|_W = c_W$ for some constant $c_W >0$. Then 
\begin{align*}
 \left\vert \left(\mu^j-\mu^* \vert \omega\right)\right\vert &= \left\vert\sum_{i=1}^{N} \left(\omega(x_i^j,U_i^j) r_i^j - \omega(x_i^*,U_i^*) r_i^* \right)\right\vert \\
 &\leq \sum_{i=1}^{N} \left(\left\vert\omega(x_i^j,U_i^j)(r_i^j-r_i^*) \right\vert + \left\vert(\omega(x_i^j,U_i^j)- \omega(x_i^*,U_i^*)) r_i^*\right\vert \right)\\
 &\leq \sum_{i=1}^{N} \left(\|\omega\|_{\infty} \vert r_i^j-r_i^*\vert +  r_i^* \|\omega\|_{1,\infty} \max\{\vert x_i^j-x_i^*\vert,\vert U_i^j-U_i^*\vert\}\right) \\
 &\leq 2N \|\omega\|_{1,\infty} \max_{i=1,\ldots,N} \{\vert x_i^j-x_i^*\vert,\vert U_i^j-U_i^*\vert,\vert r_i^j-r_i^*\vert\} \\ 
 &\leq 2N c_W \max_{i=1,\ldots,N} \{\vert x_i^m-x_i^*\vert,\vert U_i^j-U_i^*\vert,\vert r_i^j-r_i^*\vert\}
\end{align*}
Therefore, we see that $\|\mu^j-\mu^*\|_{W^*} = \sup_{\|\omega\|_W \leq 1} \left\vert\left(\mu^j-\mu^*\vert \omega\right)\right\vert \xrightarrow[j\rightarrow \infty]{} 0$ since $x_i^j \rightarrow x_i^*$, $U_i^j \rightarrow U_i^*$ and $r_i^j \rightarrow r_i^*$.

Finally, combining the previous two estimates, we obtain:
\begin{equation*}
 C(v^*,\eta^*) \leq \lim \inf_{j\rightarrow +\infty} C(v^j,\eta^j)
\end{equation*}
which shows that $(v^*,\eta^*)$ is indeed a minimizer of \eqref{eq:relaxed_LDDMM_FR_distance}.
\end{appendices}

\bibliography{biblio}


\begin{thebibliography}{51}
\ifx \bisbn   \undefined \def \bisbn  #1{ISBN #1}\fi
\ifx \binits  \undefined \def \binits#1{#1}\fi
\ifx \bauthor  \undefined \def \bauthor#1{#1}\fi
\ifx \batitle  \undefined \def \batitle#1{#1}\fi
\ifx \bjtitle  \undefined \def \bjtitle#1{#1}\fi
\ifx \bvolume  \undefined \def \bvolume#1{\textbf{#1}}\fi
\ifx \byear  \undefined \def \byear#1{#1}\fi
\ifx \bissue  \undefined \def \bissue#1{#1}\fi
\ifx \bfpage  \undefined \def \bfpage#1{#1}\fi
\ifx \blpage  \undefined \def \blpage #1{#1}\fi
\ifx \burl  \undefined \def \burl#1{\textsf{#1}}\fi
\ifx \doiurl  \undefined \def \doiurl#1{\url{https://doi.org/#1}}\fi
\ifx \betal  \undefined \def \betal{\textit{et al.}}\fi
\ifx \binstitute  \undefined \def \binstitute#1{#1}\fi
\ifx \binstitutionaled  \undefined \def \binstitutionaled#1{#1}\fi
\ifx \bctitle  \undefined \def \bctitle#1{#1}\fi
\ifx \beditor  \undefined \def \beditor#1{#1}\fi
\ifx \bpublisher  \undefined \def \bpublisher#1{#1}\fi
\ifx \bbtitle  \undefined \def \bbtitle#1{#1}\fi
\ifx \bedition  \undefined \def \bedition#1{#1}\fi
\ifx \bseriesno  \undefined \def \bseriesno#1{#1}\fi
\ifx \blocation  \undefined \def \blocation#1{#1}\fi
\ifx \bsertitle  \undefined \def \bsertitle#1{#1}\fi
\ifx \bsnm \undefined \def \bsnm#1{#1}\fi
\ifx \bsuffix \undefined \def \bsuffix#1{#1}\fi
\ifx \bparticle \undefined \def \bparticle#1{#1}\fi
\ifx \barticle \undefined \def \barticle#1{#1}\fi
\bibcommenthead
\ifx \bconfdate \undefined \def \bconfdate #1{#1}\fi
\ifx \botherref \undefined \def \botherref #1{#1}\fi
\ifx \url \undefined \def \url#1{\textsf{#1}}\fi
\ifx \bchapter \undefined \def \bchapter#1{#1}\fi
\ifx \bbook \undefined \def \bbook#1{#1}\fi
\ifx \bcomment \undefined \def \bcomment#1{#1}\fi
\ifx \oauthor \undefined \def \oauthor#1{#1}\fi
\ifx \citeauthoryear \undefined \def \citeauthoryear#1{#1}\fi
\ifx \endbibitem  \undefined \def \endbibitem {}\fi
\ifx \bconflocation  \undefined \def \bconflocation#1{#1}\fi
\ifx \arxivurl  \undefined \def \arxivurl#1{\textsf{#1}}\fi
\csname PreBibitemsHook\endcsname

%%% 1
\bibitem{Grenander1993}
\begin{bbook}
\bauthor{\bsnm{Grenander}, \binits{U.}}:
\bbtitle{General Pattern Theory: A Mathematical Study of Regular Structures}.
\bpublisher{Clarendon Press Oxford},
(\byear{1993})
\end{bbook}
\endbibitem

%%% 2
\bibitem{christensen1996deformable}
\begin{barticle}
\bauthor{\bsnm{Christensen}, \binits{G.E.}},
\bauthor{\bsnm{Rabbitt}, \binits{R.D.}},
\bauthor{\bsnm{Miller}, \binits{M.I.}}:
\batitle{Deformable templates using large deformation kinematics}.
\bjtitle{IEEE transactions on image processing}
\bvolume{5}(\bissue{10}),
\bfpage{1435}--\blpage{1447}
(\byear{1996})
\end{barticle}
\endbibitem

%%% 3
\bibitem{Michor2007}
\begin{barticle}
\bauthor{\bsnm{Michor}, \binits{P.}},
\bauthor{\bsnm{Mumford}, \binits{D.}}:
\batitle{{An overview of the Riemannian metrics on spaces of curves using the
  Hamiltonian approach}}.
\bjtitle{Applied and Computational Harmonic Analysis}
\bvolume{23}(\bissue{1}),
\bfpage{74}--\blpage{113}
(\byear{2007})
\end{barticle}
\endbibitem

%%% 4
\bibitem{Bauer2011b}
\begin{barticle}
\bauthor{\bsnm{Bauer}, \binits{M.}},
\bauthor{\bsnm{Harms}, \binits{P.}},
\bauthor{\bsnm{Michor}, \binits{P.}}:
\batitle{{Sobolev metrics on shape space of surfaces}}.
\bjtitle{J. Geom. Mech.}
\bvolume{3}(\bissue{4}),
\bfpage{389}--\blpage{438}
(\byear{2011})
\end{barticle}
\endbibitem

%%% 5
\bibitem{Beg2005}
\begin{botherref}
\oauthor{\bsnm{Beg}, \binits{M.F.}},
\oauthor{\bsnm{Miller}, \binits{M.I.}},
\oauthor{\bsnm{Trouv\'{e}}, \binits{A.}},
\oauthor{\bsnm{Younes}, \binits{L.}}:
Computing large deformation metric mappings via geodesic flows of
  diffeomorphisms.
International journal of computer vision
\textbf{61}(139-157)
(2005)
\end{botherref}
\endbibitem

%%% 6
\bibitem{Joshi2000}
\begin{barticle}
\bauthor{\bsnm{Joshi}, \binits{S.C.}},
\bauthor{\bsnm{Miller}, \binits{M.I.}}:
\batitle{Landmark matching via large deformation diffeomorphisms}.
\bjtitle{Image Processing, IEEE Transactions on}
\bvolume{9}(\bissue{8}),
\bfpage{1357}--\blpage{1370}
(\byear{2000})
\end{barticle}
\endbibitem

%%% 7
\bibitem{Glaunes2}
\begin{botherref}
\oauthor{\bsnm{Glaunès}, \binits{J.}},
\oauthor{\bsnm{Vaillant}, \binits{M.}}:
Surface matching via currents.
Proceedings of Information Processing in Medical Imaging (IPMI), Lecture Notes
  in Computer Science
\textbf{3565}(381-392)
(2006)
\end{botherref}
\endbibitem

%%% 8
\bibitem{durrleman2009statistical}
\begin{barticle}
\bauthor{\bsnm{Durrleman}, \binits{S.}},
\bauthor{\bsnm{Pennec}, \binits{X.}},
\bauthor{\bsnm{Trouv{\'e}}, \binits{A.}},
\bauthor{\bsnm{Ayache}, \binits{N.}}:
\batitle{Statistical models of sets of curves and surfaces based on currents}.
\bjtitle{Medical image analysis}
\bvolume{13}(\bissue{5}),
\bfpage{793}--\blpage{808}
(\byear{2009})
\end{barticle}
\endbibitem

%%% 9
\bibitem{Charon2013}
\begin{barticle}
\bauthor{\bsnm{Charon}, \binits{N.}},
\bauthor{\bsnm{Trouvé}, \binits{A.}}:
\batitle{The varifold representation of non-oriented shapes for diffeomorphic
  registration}.
\bjtitle{SIAM journal of Imaging Science}
\bvolume{6}(\bissue{4}),
\bfpage{2547}--\blpage{2580}
(\byear{2013})
\end{barticle}
\endbibitem

%%% 10
\bibitem{kaltenmark2017general}
\begin{bchapter}
\bauthor{\bsnm{Kaltenmark}, \binits{I.}},
\bauthor{\bsnm{Charlier}, \binits{B.}},
\bauthor{\bsnm{Charon}, \binits{N.}}:
\bctitle{A general framework for curve and surface comparison and registration
  with oriented varifolds}.
In: \bbtitle{Proceedings of the IEEE Conference on Computer Vision and Pattern
  Recognition},
pp. \bfpage{3346}--\blpage{3355}
(\byear{2017})
\end{bchapter}
\endbibitem

%%% 11
\bibitem{hsieh2020diffeomorphic}
\begin{bchapter}
\bauthor{\bsnm{Hsieh}, \binits{H.-W.}},
\bauthor{\bsnm{Charon}, \binits{N.}}:
\bctitle{Diffeomorphic registration of discrete geometric distributions}.
In: \bbtitle{Mathematics Of Shapes And Applications},
pp. \bfpage{45}--\blpage{74}.
\bpublisher{World Scientific},
(\byear{2020})
\end{bchapter}
\endbibitem

%%% 12
\bibitem{hsieh2021metrics}
\begin{barticle}
\bauthor{\bsnm{Hsieh}, \binits{H.-W.}},
\bauthor{\bsnm{Charon}, \binits{N.}}:
\batitle{Metrics, quantization and registration in varifold spaces}.
\bjtitle{Foundations of Computational Mathematics}
\bvolume{21},
\bfpage{1317}--\blpage{1361}
(\byear{2021})
\end{barticle}
\endbibitem

%%% 13
\bibitem{Gori2016}
\begin{botherref}
\oauthor{\bsnm{Gori}, \binits{P.}},
\oauthor{\bsnm{Colliot}, \binits{O.}},
\oauthor{\bsnm{Marrakchi-Kacem}, \binits{L.}},
\oauthor{\bsnm{Worbe}, \binits{Y.}},
\oauthor{\bsnm{Fallani}, \binits{F.D.V.}},
\oauthor{\bsnm{Chavez}, \binits{M.}},
\oauthor{\bsnm{Poupon}, \binits{C.}},
\oauthor{\bsnm{Hartmann}, \binits{A.}},
\oauthor{\bsnm{Ayache}, \binits{N.}},
\oauthor{\bsnm{Durrleman}, \binits{S.}}:
{Parsimonious Approximation of Streamline Trajectories in White Matter Fiber
  Bundles}.
IEEE Transactions on Medical Imaging
\textbf{PP}(99)
(2016)
\end{botherref}
\endbibitem

%%% 14
\bibitem{Trouve1}
\begin{barticle}
\bauthor{\bsnm{Trouv\'{e}}, \binits{A.}},
\bauthor{\bsnm{Younes}, \binits{L.}}:
\batitle{Metamorphoses through lie group action}.
\bjtitle{Foundation of computational mathematics}
\bvolume{5},
\bfpage{173}--\blpage{198}
(\byear{2005})
\end{barticle}
\endbibitem

%%% 15
\bibitem{Holm2009}
\begin{barticle}
\bauthor{\bsnm{Holm}, \binits{D.}},
\bauthor{\bsnm{Trouv\'e}, \binits{A.}},
\bauthor{\bsnm{Younes}, \binits{L.}}:
\batitle{{The Euler-Poincar\'e theory of Metamorphosis}}.
\bjtitle{Quart. Appl. Math}
\bvolume{67}(\bissue{2}),
\bfpage{661}--\blpage{685}
(\byear{2009})
\end{barticle}
\endbibitem

%%% 16
\bibitem{Richardson2013}
\begin{barticle}
\bauthor{\bsnm{Richardson}, \binits{C.}},
\bauthor{\bsnm{Younes}, \binits{L.}}:
\batitle{Computing metamorphoses between discrete measures}.
\bjtitle{Journal of Geometric Mechanics}
\bvolume{5}(\bissue{1}),
\bfpage{131}--\blpage{150}
(\byear{2013})
\end{barticle}
\endbibitem

%%% 17
\bibitem{Richardson2015}
\begin{botherref}
\oauthor{\bsnm{Richardson}, \binits{C.}},
\oauthor{\bsnm{Younes}, \binits{L.}}:
{Metamorphosis of images in reproducing kernel Hilbert spaces}.
Advances in Computational Mathematics,
1--31
(2015)
\end{botherref}
\endbibitem

%%% 18
\bibitem{berkels2015time}
\begin{barticle}
\bauthor{\bsnm{Berkels}, \binits{B.}},
\bauthor{\bsnm{Effland}, \binits{A.}},
\bauthor{\bsnm{Rumpf}, \binits{M.}}:
\batitle{Time discrete geodesic paths in the space of images}.
\bjtitle{SIAM Journal on Imaging Sciences}
\bvolume{8}(\bissue{3}),
\bfpage{1457}--\blpage{1488}
(\byear{2015})
\end{barticle}
\endbibitem

%%% 19
\bibitem{Charon2018}
\begin{barticle}
\bauthor{\bsnm{Charon}, \binits{N.}},
\bauthor{\bsnm{Charlier}, \binits{B.}},
\bauthor{\bsnm{Trouv\'{e}}, \binits{A.}}:
\batitle{{Metamorphoses of functional shapes in Sobolev spaces}}.
\bjtitle{J. Foundations of Comput. Math}
\bvolume{18}(\bissue{6}),
\bfpage{1535}--\blpage{1596}
(\byear{2018})
\end{barticle}
\endbibitem

%%% 20
\bibitem{Liero2016}
\begin{barticle}
\bauthor{\bsnm{Liero}, \binits{M.}},
\bauthor{\bsnm{Mielke}, \binits{A.}},
\bauthor{\bsnm{Savar\'{e}}, \binits{G.}}:
\batitle{{Optimal Transport in Competition with Reaction: The
  Hellinger--Kantorovich Distance and Geodesic Curves}}.
\bjtitle{SIAM Journal on Mathematical Analysis}
\bvolume{48}(\bissue{4}),
\bfpage{2869}--\blpage{2911}
(\byear{2016})
\end{barticle}
\endbibitem

%%% 21
\bibitem{chizat2018interpolating}
\begin{barticle}
\bauthor{\bsnm{Chizat}, \binits{L.}},
\bauthor{\bsnm{Peyr{\'e}}, \binits{G.}},
\bauthor{\bsnm{Schmitzer}, \binits{B.}},
\bauthor{\bsnm{Vialard}, \binits{F.-X.}}:
\batitle{{An interpolating distance between optimal transport and Fisher--Rao
  metrics}}.
\bjtitle{Foundations of Computational Mathematics}
\bvolume{18}(\bissue{1}),
\bfpage{1}--\blpage{44}
(\byear{2018})
\end{barticle}
\endbibitem

%%% 22
\bibitem{bronstein2009partial}
\begin{barticle}
\bauthor{\bsnm{Bronstein}, \binits{A.M.}},
\bauthor{\bsnm{Bronstein}, \binits{M.M.}},
\bauthor{\bsnm{Bruckstein}, \binits{A.M.}},
\bauthor{\bsnm{Kimmel}, \binits{R.}}:
\batitle{Partial similarity of objects, or how to compare a centaur to a
  horse}.
\bjtitle{International Journal of Computer Vision}
\bvolume{84}(\bissue{2}),
\bfpage{163}
(\byear{2009})
\end{barticle}
\endbibitem

%%% 23
\bibitem{robinson2012functional}
\begin{botherref}
\oauthor{\bsnm{Robinson}, \binits{D.T.}}:
Functional data analysis and partial shape matching in the square root velocity
  framework
(2012)
\end{botherref}
\endbibitem

%%% 24
\bibitem{rodola2017partial}
\begin{bchapter}
\bauthor{\bsnm{Rodol{\`a}}, \binits{E.}},
\bauthor{\bsnm{Cosmo}, \binits{L.}},
\bauthor{\bsnm{Bronstein}, \binits{M.M.}},
\bauthor{\bsnm{Torsello}, \binits{A.}},
\bauthor{\bsnm{Cremers}, \binits{D.}}:
\bctitle{Partial functional correspondence}.
In: \bbtitle{Computer Graphics Forum},
vol. \bseriesno{36},
pp. \bfpage{222}--\blpage{236}
(\byear{2017}).
\bcomment{Wiley Online Library}
\end{bchapter}
\endbibitem

%%% 25
\bibitem{antonsanti2021partial}
\begin{bchapter}
\bauthor{\bsnm{Antonsanti}, \binits{P.-L.}},
\bauthor{\bsnm{Glaun{\`e}s}, \binits{J.}},
\bauthor{\bsnm{Benseghir}, \binits{T.}},
\bauthor{\bsnm{Jugnon}, \binits{V.}},
\bauthor{\bsnm{Kaltenmark}, \binits{I.}}:
\bctitle{Partial matching in the space of varifolds}.
In: \bbtitle{International Conference on Information Processing in Medical
  Imaging},
pp. \bfpage{123}--\blpage{135}
(\byear{2021}).
\bcomment{Springer}
\end{bchapter}
\endbibitem

%%% 26
\bibitem{sukurdeep2021new}
\begin{botherref}
\oauthor{\bsnm{Sukurdeep}, \binits{Y.}},
\oauthor{\bsnm{Bauer}, \binits{M.}},
\oauthor{\bsnm{Charon}, \binits{N.}}:
A new variational model for the analysis of shape graphs with partial matching
  constraints.
arXiv preprint arXiv:2105.00678
(2021)
\end{botherref}
\endbibitem

%%% 27
\bibitem{attaiki2021dpfm}
\begin{botherref}
\oauthor{\bsnm{Attaiki}, \binits{S.}},
\oauthor{\bsnm{Pai}, \binits{G.}},
\oauthor{\bsnm{Ovsjanikov}, \binits{M.}}:
Dpfm: Deep partial functional maps.
arXiv preprint arXiv:2110.09994
(2021)
\end{botherref}
\endbibitem

%%% 28
\bibitem{hsieh2021diffeomorphic}
\begin{bchapter}
\bauthor{\bsnm{Hsieh}, \binits{H.-W.}},
\bauthor{\bsnm{Charon}, \binits{N.}}:
\bctitle{Diffeomorphic registration with density changes for the analysis of
  imbalanced shapes}.
In: \bbtitle{International Conference on Information Processing in Medical
  Imaging},
pp. \bfpage{31}--\blpage{42}
(\byear{2021}).
\bcomment{Springer}
\end{bchapter}
\endbibitem

%%% 29
\bibitem{Almgren}
\begin{bbook}
\bauthor{\bsnm{Almgren}, \binits{F.}}:
\bbtitle{Plateau's Problem: An Invitation to Varifold Geometry}.
\bpublisher{Student Mathematical Library},
(\byear{1966})
\end{bbook}
\endbibitem

%%% 30
\bibitem{Allard}
\begin{botherref}
\oauthor{\bsnm{Allard}, \binits{W.}}:
On the first variation of a varifold.
Annals of mathematics
\textbf{95}(3)
(1972)
\end{botherref}
\endbibitem

%%% 31
\bibitem{Simon}
\begin{bbook}
\bauthor{\bsnm{Simon}, \binits{L.}}:
\bbtitle{Lecture Notes on Geometric Measure Theory}.
\bpublisher{Australian national university},
(\byear{1983})
\end{bbook}
\endbibitem

%%% 32
\bibitem{Ambrosio2000}
\begin{bbook}
\bauthor{\bsnm{Ambrosio}, \binits{L.}},
\bauthor{\bsnm{Fusco}, \binits{N.}},
\bauthor{\bsnm{Pallara}, \binits{D.}}:
\bbtitle{Functions of Bounded Variation and Free Discontinuity Problems}.
\bpublisher{Oxford : Clarendon Press},
(\byear{2000})
\end{bbook}
\endbibitem

%%% 33
\bibitem{buet2018discretization}
\begin{barticle}
\bauthor{\bsnm{Buet}, \binits{B.}},
\bauthor{\bsnm{Leonardi}, \binits{G.P.}},
\bauthor{\bsnm{Masnou}, \binits{S.}}:
\batitle{Discretization and approximation of surfaces using varifolds}.
\bjtitle{Geometric Flows}
\bvolume{3}(\bissue{1}),
\bfpage{28}--\blpage{56}
(\byear{2018})
\end{barticle}
\endbibitem

%%% 34
\bibitem{buet2020mean}
\begin{botherref}
\oauthor{\bsnm{Buet}, \binits{B.}},
\oauthor{\bsnm{Rumpf}, \binits{M.}}:
Mean curvature motion of point cloud varifolds.
arXiv preprint arXiv:2010.09419
(2020)
\end{botherref}
\endbibitem

%%% 35
\bibitem{Younes2019}
\begin{bbook}
\bauthor{\bsnm{Younes}, \binits{L.}}:
\bbtitle{Shapes and Diffeomorphisms}.
\bpublisher{Springer},
(\byear{2019})
\end{bbook}
\endbibitem

%%% 36
\bibitem{friedrich1991fisher}
\begin{barticle}
\bauthor{\bsnm{Friedrich}, \binits{T.}}:
\batitle{Die fisher-information und symplektische strukturen}.
\bjtitle{Mathematische Nachrichten}
\bvolume{153}(\bissue{1}),
\bfpage{273}--\blpage{296}
(\byear{1991})
\end{barticle}
\endbibitem

%%% 37
\bibitem{bauer2016uniqueness}
\begin{barticle}
\bauthor{\bsnm{Bauer}, \binits{M.}},
\bauthor{\bsnm{Bruveris}, \binits{M.}},
\bauthor{\bsnm{Michor}, \binits{P.W.}}:
\batitle{{Uniqueness of the Fisher--Rao metric on the space of smooth
  densities}}.
\bjtitle{Bulletin of the London Mathematical Society}
\bvolume{48}(\bissue{3}),
\bfpage{499}--\blpage{506}
(\byear{2016})
\end{barticle}
\endbibitem

%%% 38
\bibitem{Glaunes2004}
\begin{barticle}
\bauthor{\bsnm{Glaun\`es}, \binits{J.}},
\bauthor{\bsnm{Trouv\'e}, \binits{A.}},
\bauthor{\bsnm{Younes}, \binits{L.}}:
\batitle{Diffeomorphic matching of distributions: A new approach for unlabelled
  point-sets and sub-manifolds matching}.
\bjtitle{IEEE Computer Society Conference on Computer Vision and Pattern
  Recognition}
\bvolume{2},
\bfpage{712}--\blpage{718}
(\byear{2004})
\end{barticle}
\endbibitem

%%% 39
\bibitem{Roussillon2016}
\begin{barticle}
\bauthor{\bsnm{Roussillon}, \binits{P.}},
\bauthor{\bsnm{Glaunès}, \binits{J.}}:
\batitle{{Kernel Metrics on Normal Cycles and Application to Curve Matching}}.
\bjtitle{SIAM Journal on Imaging Sciences}
\bvolume{9}(\bissue{4}),
\bfpage{1991}--\blpage{2038}
(\byear{2016})
\end{barticle}
\endbibitem

%%% 40
\bibitem{Feydy2017}
\begin{bchapter}
\bauthor{\bsnm{Feydy}, \binits{J.}},
\bauthor{\bsnm{Charlier}, \binits{B.}},
\bauthor{\bsnm{Vialard}, \binits{F.-X.}},
\bauthor{\bsnm{Peyr{\'e}}, \binits{G.}}:
\bctitle{{Optimal Transport for Diffeomorphic Registration}}.
In: \bbtitle{Medical Image Computing and Computer Assisted Intervention},
pp. \bfpage{291}--\blpage{299}
(\byear{2017})
\end{bchapter}
\endbibitem

%%% 41
\bibitem{feydy2019interpolating}
\begin{bchapter}
\bauthor{\bsnm{Feydy}, \binits{J.}},
\bauthor{\bsnm{S{\'e}journ{\'e}}, \binits{T.}},
\bauthor{\bsnm{Vialard}, \binits{F.-X.}},
\bauthor{\bsnm{Amari}, \binits{S.-i.}},
\bauthor{\bsnm{Trouv{\'e}}, \binits{A.}},
\bauthor{\bsnm{Peyr{\'e}}, \binits{G.}}:
\bctitle{{Interpolating between optimal transport and MMD using Sinkhorn
  divergences}}.
In: \bbtitle{The 22nd International Conference on Artificial Intelligence and
  Statistics},
pp. \bfpage{2681}--\blpage{2690}
(\byear{2019}).
\bcomment{PMLR}
\end{bchapter}
\endbibitem

%%% 42
\bibitem{charon2020fidelity}
\begin{bchapter}
\bauthor{\bsnm{Charon}, \binits{N.}},
\bauthor{\bsnm{Charlier}, \binits{B.}},
\bauthor{\bsnm{Glaun{\`e}s}, \binits{J.}},
\bauthor{\bsnm{Gori}, \binits{P.}},
\bauthor{\bsnm{Roussillon}, \binits{P.}}:
\bctitle{Fidelity metrics between curves and surfaces: currents, varifolds, and
  normal cycles}.
In: \bbtitle{Riemannian Geometric Statistics in Medical Image Analysis},
pp. \bfpage{441}--\blpage{477}.
\bpublisher{Elsevier},
(\byear{2020})
\end{bchapter}
\endbibitem

%%% 43
\bibitem{Aronszajn1950}
\begin{barticle}
\bauthor{\bsnm{Aronszajn}, \binits{N.}}:
\batitle{Theory of reproducing kernels}.
\bjtitle{Trans. Amer. Math. Soc.}
\bvolume{68},
\bfpage{337}--\blpage{404}
(\byear{1950})
\end{barticle}
\endbibitem

%%% 44
\bibitem{Sriperumbudur10}
\begin{bchapter}
\bauthor{\bsnm{Sriperumbudur}, \binits{B.K.}},
\bauthor{\bsnm{Fukumizu}, \binits{K.}},
\bauthor{\bsnm{Lanckriet}, \binits{G.}}:
\bctitle{{On the relation between universality, characteristic kernels and RKHS
  embedding of measures}}.
In: \bbtitle{Proceedings of the Thirteenth International Conference on
  Artificial Intelligence and Statistics (AISTATS-10)},
vol. \bseriesno{9},
pp. \bfpage{773}--\blpage{780}
(\byear{2010})
\end{bchapter}
\endbibitem

%%% 45
\bibitem{glaunes2014matrix}
\begin{barticle}
\bauthor{\bsnm{Glaun\`{e}s}, \binits{J.}},
\bauthor{\bsnm{Micheli}, \binits{M.}}:
\batitle{Matrix-valued kernels for shape deformation analysis}.
\bjtitle{Imaging and Computing}
\bvolume{1}(\bissue{1}),
\bfpage{57}--\blpage{139}
(\byear{2014})
\end{barticle}
\endbibitem

%%% 46
\bibitem{Pontryagin1962}
\begin{bbook}
\bauthor{\bsnm{Pontryagin}, \binits{L.}},
\bauthor{\bsnm{Boltyanskii}, \binits{V.}},
\bauthor{\bsnm{Gamkrelidze}, \binits{R.}},
\bauthor{\bsnm{Mishchenko}, \binits{E.}}:
\bbtitle{The Mathematical Theory of Optimal Processes}.
\bpublisher{John Wiley $\&$ Sons},
(\byear{1962})
\end{bbook}
\endbibitem

%%% 47
\bibitem{Glaunes2008}
\begin{barticle}
\bauthor{\bsnm{Glaunès}, \binits{J.}},
\bauthor{\bsnm{Qiu}, \binits{A.}},
\bauthor{\bsnm{Miller}, \binits{M.}},
\bauthor{\bsnm{Younes}, \binits{L.}}:
\batitle{Large deformation diffeomorphic metric curve mapping}.
\bjtitle{International Journal of Computer Vision}
\bvolume{80}(\bissue{3}),
\bfpage{317}--\blpage{336}
(\byear{2008})
\end{barticle}
\endbibitem

%%% 48
\bibitem{charlier2021kernel}
\begin{barticle}
\bauthor{\bsnm{Charlier}, \binits{B.}},
\bauthor{\bsnm{Feydy}, \binits{J.}},
\bauthor{\bsnm{Glaun{\`e}s}, \binits{J.}},
\bauthor{\bsnm{Collin}, \binits{F.-D.}},
\bauthor{\bsnm{Durif}, \binits{G.}}:
\batitle{{Kernel operations on the GPU, with autodiff, without memory
  overflows}}.
\bjtitle{Journal of Machine Learning Research}
\bvolume{22}(\bissue{74}),
\bfpage{1}--\blpage{6}
(\byear{2021})
\end{barticle}
\endbibitem

%%% 49
\bibitem{Miller2015b}
\begin{barticle}
\bauthor{\bsnm{Miller}, \binits{M.}},
\bauthor{\bsnm{Younes}, \binits{L.}},
\bauthor{\bsnm{Ratnanather}, \binits{J.}},
\bauthor{\bsnm{Brown}, \binits{T.}},
\bauthor{\bsnm{Trinh}, \binits{H.}},
\bauthor{\bsnm{Lee}, \binits{D.}},
\bauthor{\bsnm{Tward}, \binits{D.}},
\bauthor{\bsnm{Mahon}, \binits{P.}},
\bauthor{\bsnm{Mori}, \binits{S.}},
\bauthor{\bsnm{Albert}, \binits{M.}}:
\batitle{{Amygdalar atrophy in symptomatic Alzheimer's disease based on
  diffeomorphometry: the BIOCARD cohort}}.
\bjtitle{Neurobiology of Aging}
\bvolume{36 Supplement 1},
\bfpage{3}--\blpage{10}
(\byear{2015})
\end{barticle}
\endbibitem

%%% 50
\bibitem{durrett2019probability}
\begin{bbook}
\bauthor{\bsnm{Durrett}, \binits{R.}}:
\bbtitle{Probability: Theory and Examples}
vol. \bseriesno{49}.
\bpublisher{Cambridge university press},
(\byear{2019})
\end{bbook}
\endbibitem

%%% 51
\bibitem{arguillere14:_shape}
\begin{barticle}
\bauthor{\bsnm{Arguillere}, \binits{S.}},
\bauthor{\bsnm{Tr\'elat}, \binits{E.}},
\bauthor{\bsnm{Trouv\'e}, \binits{A.}},
\bauthor{\bsnm{Younes}, \binits{L.}}:
\batitle{Shape deformation analysis from the optimal control viewpoint}.
\bjtitle{{Journal de Mathématiques Pures et Appliquées}}
\bvolume{104}(\bissue{1}),
\bfpage{139}--\blpage{178}
(\byear{2015})
\end{barticle}
\endbibitem

\end{thebibliography}

\end{document}